\documentclass[12pt]{article}

\usepackage[dvips]{graphicx} 
\usepackage[letterpaper,margin=1in]{geometry}
\usepackage{times}
\usepackage{setspace}       
\usepackage{amsmath,amsthm,comment}
\usepackage[colorlinks=true, pdfstartview=FitV,
linkcolor=blue, citecolor=blue, urlcolor=blue]{hyperref}
\usepackage{algorithm, algorithmic}
\usepackage{caption}
\usepackage{subcaption}
\usepackage{bbm,mathalpha,mathrsfs}
\usepackage{natbib}
\usepackage{color} 
\definecolor{myblue}{rgb}{.8, .8, 1}
\usepackage{relsize}
\usepackage{mathtools}
\usepackage[dvipsnames]{xcolor}
\usepackage{amssymb,fancyhdr}
\usepackage{url}
\usepackage{enumitem}

\usepackage{mathabx}
\usepackage{mathtools}
\usepackage{etoolbox}
\usepackage{accents}
\usepackage{multicol}
\usepackage{booktabs}
\usepackage{upgreek}
\usepackage{empheq}
\usepackage[most]{tcolorbox}
\DeclareFontFamily{U}{matha}{\hyphenchar\font45}
\DeclareFontShape{U}{matha}{m}{n}{
<-6> matha5 <6-7> matha6 <7-8> matha7
<8-9> matha8 <9-10> matha9
<10-12> matha10 <12-> matha12
}{}
\DeclareSymbolFont{matha}{U}{matha}{m}{n}

\DeclareFontFamily{U}{mathx}{\hyphenchar\font45}
\DeclareFontShape{U}{mathx}{m}{n}{
<-6> mathx5 <6-7> mathx6 <7-8> mathx7
<8-9> mathx8 <9-10> mathx9
<10-12> mathx10 <12-> mathx12
}{}
\DeclareSymbolFont{mathx}{U}{mathx}{m}{n}

\DeclareMathDelimiter{\vvvert} {0}{matha}{"7E}{mathx}{"17}%

\DeclarePairedDelimiterX{\mnorm}[1]
{\vvvert}
{\vvvert}
{\ifblank{#1}{\:\cdot\:}{#1}}

\DeclareMathOperator*{\argmin}{arg\,min}

\newcommand{\ip}[1]{\langle#1\rangle}

\newcommand{\Norm}[1]{\left\Vert#1\right\Vert}
\newcommand{\norm}[1]{\Vert#1\Vert}
\newcommand{\Abs}[1]{\left|#1\right|}
\newcommand{\abs}[1]{|#1|}
\newcommand{\Set}[1]{\left\{#1\right\}}
\newcommand{\set}[1]{\{#1\}}
\newcommand{\Real}{\mathbb R}

\def\ev{\mathbf e}

\def\gv{\mathbf g}

\def\gv{\mathbf g}

\def\Av{\mathbf A}

\def\Iv{\mathbf I}

\newcommand\mybluebox{%
    \@ifnextchar[
       {\@mybluebox}%
       {\@mybluebox[0pt]}}

\def\@mybluebox[#1]{%
    \@ifnextchar[
       {\@@mybluebox[#1]}%
       {\@@mybluebox[#1][0pt]}}

\def\@@mybluebox[#1][#2]#3{
    \sbox\mytempbox{#3}%
    \mytemplen\ht\mytempbox
    \advance\mytemplen #1\relax
    \ht\mytempbox\mytemplen
    \mytemplen\dp\mytempbox
    \advance\mytemplen #2\relax
    \dp\mytempbox\mytemplen
    \colorbox{myblue}{\hspace{1em}\usebox{\mytempbox}\hspace{1em}}}

\newtcbox{\mymath}[1][]{%
    nobeforeafter, math upper, tcbox raise base,
    enhanced, colframe=blue!30!black,
    colback=blue!30, boxrule=1pt,
    #1}

\newcommand{\dto}{\xrightarrow{d}}

\newcommand{\Ac}{\mathcal{A}}

\newcommand{\Cc}{\mathcal{C}}
\newcommand{\Dc}{\mathcal{D}}
\newcommand{\Ec}{\mathcal{E}}

\newcommand{\Ic}{\mathcal{I}}

\newcommand{\Lc}{\mathcal{L}}
\newcommand{\Mc}{\mathcal{M}}
\newcommand{\Nc}{\mathcal{N}}

\newcommand{\Pc}{\mathcal{P}}

\newcommand{\Sc}{\mathcal{S}}

\newcommand{\Wc}{\mathcal{W}}

\newcommand{\Yc}{\mathcal{Y}}
\newcommand{\Zc}{\mathcal{Z}}

\newcommand{\Bb}{\mathbb{B}}

\newcommand{\Eb}{\mathbb{E}}

\newcommand{\Pb}{\mathbb{P}}

\newcommand{\Sb}{\mathbb{S}}

\newcommand{\Ck}{\mathfrak{C}}

\newcommand{\Fk}{\mathfrak{F}}

\newcommand{\Rk}{\mathfrak{R}}

\usepackage{authblk}

\renewcommand{\hat}{\widehat}

\newtheorem{theorem}{Theorem}

\numberwithin{intassumption}{assumption}

\newtheorem{lemma}{Lemma}
\newtheorem{remark}{Remark}
\newtheorem{proposition}{Proposition}

\date{\today}

\title{Assumption-Lean Honest Inference for $Z$-functionals}
\author{Woonyoung Chang}
\author{Arun Kumar Kuchibhotla}
\affil{Department of Statistics \& Data Science, Carnegie Mellon University}

\begin{document}

\maketitle
\begin{abstract}
We develop a general assumption-lean framework for constructing uniformly valid confidence sets for functionals defined by moment equalities, referred to as $Z$-functionals. Our approach combines self-normalized statistics with a test inversion principle, enabling honest inference under mild regularity conditions and without explicit variance estimation. To enhance geometric tractability, we propose novel split-normalized and Gateaux-normalized statistics that yield computationally feasible and interpretable confidence sets. A central contribution of this work is a comprehensive non-asymptotic width analysis: we derive high-probability upper bounds on the diameter of the proposed confidence sets, and quantify their proximity to Wald intervals under minimal assumptions. Applications to high-dimensional non-sparse linear and generalized linear regression demonstrate that our procedures achieve valid coverage and near-optimal rate of convergence for the width/diameter, while the classical methods including Wald and bootstrap fail. 
\end{abstract}
\section{Introduction}
This paper studies the problem of inference for parameters defined implicitly through a system of population-level equations or so-called moment equalities. Given a measurable function $\psi : \Zc \times \Theta \to \mathbb{R}^d$, where $\Theta\subseteq\Real^d$, and a distribution $P\in\Pc$, where $\Pc$ is a class of distribution on a sample space $\Wc$, our inferential target is a functional, defined as\begin{equation}\label{eq:target_iid} \theta_0 := \theta_0(P) \quad \text{such that} \quad \mathbb{E}_{Z\sim P}[\psi(Z, \theta_0)] = 0_d. \end{equation} We refer to such a map $P\mapsto\theta_0(P)$ as a $Z$-functional. This formulation includes many classical and modern statistical problems such as smooth $M$-estimation and the generalized method of moments, and further supports model-agnostic perspectives in which the target parameter is defined in the absence of any statistical model \citep{buja2019models1,buja2019models2}. It is conceptually related to partially identified models defined by moment inequalities, an area of large interest in empirical economics \citep{tamer2003incomplete,chernozhukov2007estimation,andrews2010inference,romano2010inference}, where the target may be set-valued. Even with moment equalities (as in our setting), the target $\theta_0(P)$ need not be unique, especially if the moment equation is obtained as the derivative of a non-convex function. All our coverage guarantees are valid with the interpretation $\mathbb{P}_P(\theta_0(P)\notin\widehat{\mathrm{CI}}_{n,\alpha}) = \sup_{\theta^*\in\theta_0(P)}\mathbb{P}_P(\theta^*\notin\widehat{\mathrm{CI}}_{n,\alpha})$. Note that, under such non-uniqueness, the diameters of the proposed confidence sets do not converge to zero in probability. For brevity, throughout, we assume $\theta_0(P)$ is a singleton.

The primary goal of this paper is to construct valid confidence sets for $\tau(P)=\theta_0(P)$ or $\tau(P)=\xi^\top(\theta_0(P))$ (for a known contrast vector $\xi$), under minimal assumptions. We focus on honest confidence sets $\widehat{\mathrm{CI}}_{n,\alpha}$ over a class of distributions $\Pc$, i.e., for all $\alpha\in(0, 1)$, 
\begin{equation*}\label{eq:uniformvalidity}
\limsup_{n\to\infty}\sup_{P\in\Pc}\big(\Pb_P(\tau(P)\notin\widehat{\mathrm{CI}}_{n,\alpha} ) - \alpha\big)_+ = 0,
\end{equation*} where $n$ is the sample size. This notion is stronger than pointwise validity, which requires
\begin{equation*}\label{eq:pointwisevalidity}
\sup_{P\in\Pc}\limsup_{n\to\infty}\big(\Pb_P(\tau(P)\notin\widehat{\mathrm{CI}}_{n,\alpha} ) - \alpha\big)_+ = 0.
\end{equation*} Pointwise-valid procedures can severely under-cover along some distribution sequences \citep{romano2000finite,poetscher2002lower}.

While valid coverage is the primary focus of confidence sets, their utility crucially depends on their shape and size. Our secondary goal is to analyze the worst-case diameter of the proposed confidence sets along with their ``asymptotic’’ shape. Given that the proposed confidence sets are implicit (unlike the traditional Wald or resampling ones), this secondary goal is also of importance. Specifically, we seek a sequence $s_n$ such that
\begin{equation}\label{eq:intro_width}
    \limsup_{n\to\infty}\sup_{P\in\Pc}\Pb_P({\rm diam}_{\norm{\cdot}}(\widehat{\mathrm{CI}}_{n,\alpha})\geq s_n)= 0,
\end{equation} where ${\rm diam}_{\norm{\cdot}}(\widehat{\mathrm{CI}}_{n,\alpha}) = \sup\{\norm{\tau_1-\tau_2}:\tau_1,\tau_2\in\widehat{\mathrm{CI}}_{n,\alpha} \}$ and $\norm{\cdot}$ is an appropriate norm, often the Euclidean norm. Additionally, we seek a set $\widehat{\mathrm{S}}_{n,\alpha}$ with an interpretable shape such as a hyperrectangle or an ellipse, and a sequence $\gamma_n = o(s_n)$ as $n\to\infty$ such that
\begin{equation}\label{eq:asymptotic-shape}
\limsup_{n\to\infty}\sup_{P\in\mathcal{P}}\, \mathbb{P}_P({\rm Haus}_{\|\cdot\|}(\widehat{\mathrm{CI}}_{n,\alpha}, \widehat{\mathrm{S}}_{n,\alpha}) \ge \gamma_n) = 0,
\end{equation}
where ${\rm Haus}_{\|\cdot\|}(A, B)$ is the Hausdorff distance between sets $A$ and $B$ measured in terms of the norm $\norm{\cdot}$.  The convergence in \eqref{eq:intro_width} ensures that the confidence set contracts at rate $s_n$ uniformly over $\mathcal{P}$, while the convergence in~\eqref{eq:asymptotic-shape} ensures that the confidence set has approximately the same shape as that of $\widehat{\mathrm{S}}_{n,\alpha}$. It is important that $\gamma_n$ converges to zero faster than $s_n$; for any sequence of sets $\widehat{\mathrm{E}}_{n}$ with $\mbox{diam}_{\|\cdot\|}(\widehat{\mathrm{E}}_n) = o_P(s_n)$ and $\widehat{\mathrm{E}}_n\cap\widehat{\mathrm{CI}}_{n,\alpha}\neq\emptyset$, we have $\mbox{Haus}_{\|\cdot\|}(\widehat{\mathrm{CI}}_{n,\alpha}, \widehat{\mathrm{E}}_n) = o_P(s_n)$. We refer to such Hausdorff convergence of confidence sets as a second-order analysis of the confidence set. 

We remark that in the context of duality of confidence sets and hypothesis tests, the analysis of the diameter yields a stronger statement than the power analysis. The diameter is more closely related to the separation rate for distinguishable hypotheses. We refer the reader to Section~6 of \cite{robins2006adaptive} and also to \cite{nickl2013confidence}.

\subsection{Literature Review}
\paragraph{Standard Inferences Based on Estimators and Their Limitations.}
A common route in inference for a general functional $\tau:=\tau(P)$ is to rely on pointwise asymptotics of an estimator $\hat\tau_n$, \begin{equation}\label{eq:intro_distappproximation}
    b_n^{-1}(\hat\tau_n -\tau)\dto \Lc(\omega_P)\quad\mbox{for all}\quad P\in\Pc,
\end{equation} for some convergence rate $b_n$, the limit distribution $\Lc$, and a nuisance parameter $\omega_P$. A canonical example is asymptotic normality where $b_n=n^{-1/2}$ and $\mathcal{L}(\omega_P) = \sigma_P\mathcal{N}(0, 1)$. Asymptotic normality with root-n consistency is by no means the only possible setting. For instance, many estimators enjoy cube-root asymptotics with non-Gaussian limits \citep{kim1990cube,Pflug1995,delgado2001subsampling,cattaneo2020practical,einmahl2022cube}.

Two standard approaches to construct a pointwise valid confidence set are: (1) the Wald method and (2) resampling methods. The Wald method, assuming a consistent estimator $\hat\omega_P$, yields a confidence interval $[\hat\tau_n-b_nq_{\alpha/2}(\Lc(\hat \omega_P)), \hat\tau_n- b_nq_{1-\alpha/2}(\Lc(\hat \omega_P)) ]$ where $q_\gamma$ denotes the upper $\gamma$-th quantile functional. In contrast, resampling methods do not require additional estimation procedures for nuisance parameters. For instance, $m_n$-out-of-$n$ bootstrap (or sub-sampling) relies on approximating the distribution of $b_n^{-1}(\hat\tau_n -\tau)$ by the conditional distribution of $b_{m_n}^{-1}(\hat\tau_{m_n}-\hat\tau_n)$, whereas the bootstrap validity typically requires the suitable stability of the functional $\tau(P)$ in $P$, such as Hadamard differentiability \citep{dumbgen1993nondifferentiable, rockafellar1998variational} and sub-sampling requires knowledge of the convergence rate and a careful calibration of the subsample size \citep{politis1999subsampling,bickel1997resampling}.

While these are foundational methods, yielding pointwise valid inferences, they do not guarantee uniform validity in general \citep{AndrewsGuggenberger2009,andrews2010asymptotic}. Among others, one source of complications arises when the convergence rate $b_n$ itself depends on the unknown features of $P$.  For example, the sample median converges at rate $b_n = n^{-1/(2\beta)}$ where $\beta$ denotes the local H\"older smoothness exponent of the underlying distribution function \citep{knight1998limiting,knight2002limiting}. In such cases, the rate $b_n$ should be estimated (see, for e.g., \cite{bertail1999subsampling}), but achieving uniform consistency of a rate estimator is generally difficult. 

In response to such irregular scenarios, recent work has developed robust inference procedures that maintain uniform validity under minimal assumptions \citep{wasserman2020universal,kuchibhotla2023hulc,park2023robust,kim2024dimension,takatsu2025bridging,kim2025locally}. A common strategy in these frameworks is to decouple estimation from inference by sample splitting \citep{wasserman2009high,meinshausen2009p,rinaldo2019bootstrapping,wasserman2020universal,takatsu2025bridging} or cross-fitting \citep{chernozhukov2018double,chernozhukov2022locally,kim2024dimension}, allowing valid inference often at a small price in diameters. Not surprisingly, the uniform validity of these methods sometimes yields conservative inference in regular settings \citep{nguyen2020universal,shi2024universal}, reflecting a standard notion of trade-off between robustness and efficiency. A possible refinement to mitigate the conservativeness was studied in \cite{takatsu2025precise}.

\paragraph{Validity and Width Analysis in Increasing Dimensions.}
Validity of the traditional confidence sets is suspect when the dimension is allowed to grow with the sample size. For example, in the context of linear regression, asymptotic normality of the traditional $Z$-estimator holds only when the dimension $d$ grows slower than $\sqrt{n}$ (i.e., $d = o(\sqrt{n})$) and thus the validity of the Wald interval is limited to such low-dimensional regimes. Similarly, validity guarantees for resampling methods relying on asymptotic distribution could fail for larger dimensions. Only a handful of works provide asymptotically valid confidence sets for the general range of $d = o(n)$ without relying on \emph{well-specified linear model} and/or \emph{sparsity}. \cite{takatsu2025bridging} is one such example. In addition, none of these works provides sharp (root-n) inference for linear contrasts of the projection parameter. Although not related to our work, we mention here that there are numerous inference methods in the literature for high-dimensional $Z$-estimators that rely on sparsity assumptions and derive minimax optimal confidence sets~\citep{Belloni2015Uniform}. In the context of non-sparse fitting of linear regression with increasing dimension and (potential) misspecification, \cite{mammen1993bootstrap} proved the validity of bootstrap confidence intervals when $d = o(n^{3/4})$ assuming a finite 24-th moment on the covariates. \cite{chang2023inference} construct honest intervals for linear contrasts, proving $1/\sqrt{n}$ width and validity when $d = \tilde{o}(n^{2/3})$\footnote{Here, $\tilde{o}$ denotes the small-o notation, which disregards the poly-logarithmic dependency on $n$.}, assuming a finite 12-th moment on the covariates. Recently, \cite{lin2024worthwhilejackknifebreakingquadratic} constructed honest confidence sets for general $Z$-estimators assuming thrice differentiability and sub-Gaussianity of $\psi$ along with a restrictive assumption of compact parameter space. The existence of honest confidence intervals for general $Z$-functionals that remain valid for the range of $d = o(n)$ remains unknown.

\subsection{Main Contributions} This paper adheres to an assumption-lean framework for inference on $Z$-functionals, imposing minimal assumptions and providing non-asymptotic guarantees. Our contributions are as follows.
\paragraph*{Honest, Nuisance-Free, and Computationally Tractable Confidence Sets.}
We propose a family of self-normalized and sample-splitting-based statistics \ref{eq:def.SeNstatistics}, \ref{eq:def.SpNstatistics}, and \ref{eq:def.TNstatistics} for honest inference. Unlike the traditional inference methods, our methods do not depend on the limiting distribution of an estimator and, hence, avoid the estimation of asymptotic variance or quantile. With a purely data-driven quantile estimate using a preliminary estimator, the resulting confidence sets are both computationally efficient and geometrically tractable. We derive non-asymptotic coverage guarantees under weak regularity conditions, moment assumptions, and the consistency of the preliminary estimator. When applied to high-dimensional linear and generalized linear regression, these confidence sets remain uniformly valid provided that $d = \tilde{o}(n)$.

\paragraph*{Rate-Optimal Concentration for $Z$-Estimators.}
By construction, all the proposed confidence sets contain a $Z$-estimator, $\widehat{\theta}$, the solution to an empirical version of~\eqref{eq:target_iid}. Given the asymptotic validity of the confidence sets, we get that the expected diameter of our confidence sets is asymptotically larger than $(1-\alpha)\|\widehat{\theta} - \theta_0\|_2$. Given the lack of results on the rate of convergence of $Z$-estimators under our mild conditions, we provide finite-sample concentration inequality (Theorem~\ref{thm:7}) for $Z$-estimators, significantly generalizing the results of \cite{lin2024worthwhilejackknifebreakingquadratic} and \cite{gauss2024asymptotics} by removing their restrictive compact parameter space assumption, among others. Applied to generalized linear models, Theorem~\ref{thm:11} proves $\sqrt{n/d}$-consistency under finite $(2+\varepsilon)$-moment assumptions and without strong regularity. To our knowledge, this is the first such result that accommodates both heavy-tailed covariates and unbounded parameter spaces. The rate closely mirrors that of least squares estimators in linear regression \citep{oliveira2016lower,kuchibhotla2020berry}, and therefore cannot be substantially improved in general. See also \cite{chardon2024finite} for the $\sqrt{n/d}$-rate optimality results for the MLE in misspecified logistic regression.

\paragraph*{Non-Asymptotic Width Control.}
We provide bounds on the diameter of the confidence set constructed via the SpN statistic for general $Z$-functionals under a mild curvature and a (one-sided quantitative)  stochastic equicontinuity condition. When applied to misspecified generalized linear models, Theorem~\ref{thm:13} shows that the diameter scales as $\sqrt{d \log(ed)/n}$ with probability converging to 1, matching the convergence rate of the $Z$-estimator up to logarithmic factors. This result holds under finite moment assumptions and does not rely on model specification or strong regularity conditions. In misspecified linear regression, we further show that the SeN, SpN, and GaN procedures yield confidence sets with a minimax optimal diameter scaling of $\sqrt{d/n}$ with a high probability.

\paragraph*{Honest confidence intervals for linear contrasts and second order width analysis.} We introduce a general framework for constructing confidence intervals for one-dimensional contrasts via directional projection. These intervals are uniformly valid and, under suitable regularity conditions, converge to the classical Wald interval in Hausdorff distance. In high-dimensional linear regression, the resulting confidence intervals are honest under finite-moment assumptions and dimension growth $d= o(n)$. Moreover, they achieve a width scaling of $1/\sqrt{n} + d/n$ with a high probability. To our knowledge, these are the first confidence intervals to attain this rate under moment assumptions and model misspecification.

\subsection{Outline}
The rest of the paper is organized as follows. In Section~2, we introduce the proposed confidence sets based on self-normalized (SeN), split-normalized (SpN), and Gateaux-normalized (GaN) statistics, and establish their non-asymptotic coverage guarantees under weak moment and regularity conditions. Section~3 focuses on inference for one-dimensional linear contrasts and develops uniformly valid confidence intervals with sharp second-order width bounds. Section~4 applies the framework to high-dimensional non-sparse linear and generalized linear regression, deriving explicit rates for coverage and width. Section~5 presents numerical studies that evaluate the finite-sample performance of the proposed procedures in various data-generating processes. Section~6 concludes with a comprehensive summary and a discussion of potential extensions and open problems. Proofs of theorems, corollaries, and lemmas are included in Appendix.

\paragraph{Notations} For two real numbers $a$ and $b$, let $a\vee b=\max\{a,b\}$, $a\wedge b=\min\{a,b\}$, and $a_+=\max\{a,0\}$. For any integer $m \geq 1$, let $[m] := \set{1,\ldots,m}$. We denote by $I_d$ the $d \times d$ identity matrix, $1_d$ the $d$-dimensional vector of ones, and $0_d$ of zeros. Let $\ev_j$ denote the $j$th canonical basis vector in $\mathbb{R}^d$, for $j \in [d]$. For $x \in \mathbb{R}^d$, let $\norm{x}_p$ denote the usual $\ell_p$-norm for $p \in [0,\infty]$, and let $\mathbb{S}^{d-1} := \set{ \eta \in \mathbb{R}^d : \norm{\eta}_2 = 1 }$ denote the unit sphere in $\mathbb{R}^d$. For a matrix $M \in \mathbb{R}^{d \times d}$, define the operator norm $\norm{M}_{\rm op} := \sup_{\eta \in \mathbb{S}^{d-1}} \norm{M\eta}_2$, and let ${\rm diag}(M)$ denote the diagonal matrix with the same diagonal elements as $M$. For positive definite $M$, define ${\rm Corr}(M) := {\rm diag}(M)^{-\frac{1}{2}} M {\rm diag}(M)^{-\frac{1}{2}}$. The convergence in distribution and probability are denoted by $\xrightarrow{d}$ and $\xrightarrow{P}$, respectively. The smallest and largest singular values of a matrix are denoted by $\sigma_{\min}(\cdot)$ and $\sigma_{\max}(\cdot)$, respectively, while for symmetric matrices we use $\lambda_{\min}(\cdot)$ and $\lambda_{\max}(\cdot)$ to denote the smallest and largest eigenvalues. For two positive sequences $\set{a_n}$ and $\set{b_n}$, we denote $a_n=o(b_n)$ if $a_n/b_n\to0$ as $n\to\infty$, and $a_n=\tilde{o}(b_n)$ if $a_n\log^{k}(en)/b_n\to0$ for some $k<\infty$ as $n\to\infty$.

\section{Inference for high-dimensional \texorpdfstring{$Z$}{Z}-functional}
The motivation for our approach is the reduction from functional inference to high-dimensional mean inference problems via the inversion principle. This principle has widely used in the construction of confidence sets for partially identified models \citep{imbens2004confidence,romano2014practical, andrews2016conditional}. The defining property of the $Z$-functional is that $\psi(Z, \theta_0(P))$ is a $d$-dimensional mean zero random vector under $P$ for all $P\in\Pc$. For a fixed $\theta\in\Theta$, let $\hat\Ck_{n,\alpha}(\theta)\subseteq\Real^d$ be a confidence set constructed from $\psi(Z_1, \theta),\ldots,\psi(Z_n, \theta)$ for $\Eb_P[\psi(Z,\theta_0(P)]$ and define the inversion-based confidence set $\widehat{\mathrm{CI}}_{n, \alpha} := \{ \theta \in \Theta : 0_d \in \hat\Ck_{n,\alpha}(\theta) \}.$ From the construction, it follows that 
\begin{equation*}
    \limsup_{n\to\infty}\sup_{P\in\Pc}\mathbb{P}_P( \theta_0(P) \notin \widehat{\mathrm{CI}}_{n, \alpha} ) \leq  \limsup_{n\to\infty}\sup_{P\in\Pc}\mathbb{P}_P( 0_d \notin \widehat\Ck_{n, \alpha}(\theta_0(P))).
\end{equation*} Thus, honest inference for $Z$-functionals reduces to the construction of an honest confidence set for the high-dimensional mean $\Eb_P[\psi(Z,\theta)]$, evaluated at $\theta=\theta_0(P)$.

\begin{figure}[h]
    \centering
    \includegraphics[width=.6\linewidth]{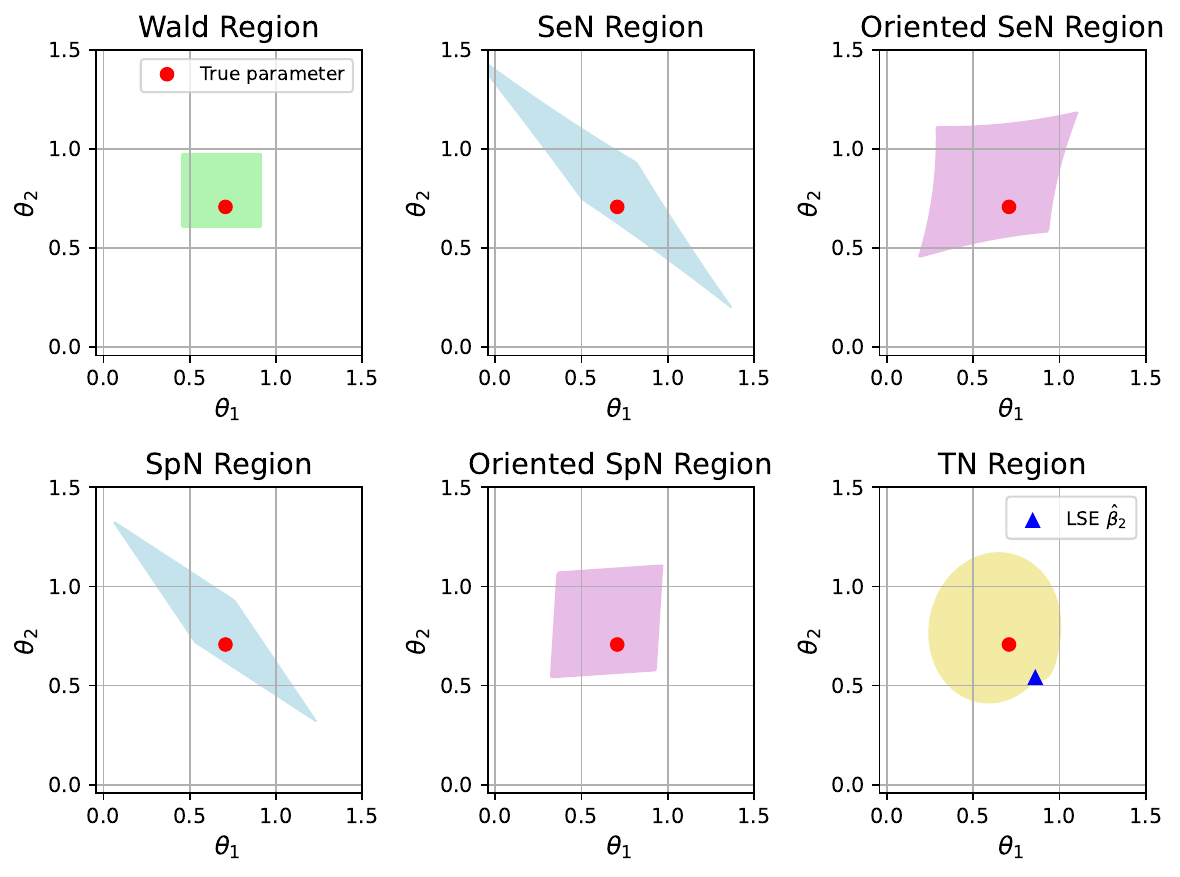}
    \caption{Visualization of confidence regions for the projection parameter~\eqref{eq:def_projection_parameters} in linear regression with two covariates without intercept. Sample size $n = 200$. We split the data into two parts ($\mathcal{D}_1$ and $\mathcal{D}_2$) each with 100 observations. Preliminary estimator is $\widehat{\beta}_2$, the OLSE from $\mathcal{D}_2$. The plots show the comparison of shapes of standard Wald region to our proposed methods with different choices. (i) Wald region: obtained from OLSE (from full sample) and the sandwich variance estimator~\citep{huber1967behavior} with a Bonferroni correction (baseline); (ii) SeN region: the SeN statistic using $Q = I_d$; (iii) Oriented SeN region: the SeN statistic using $Q = \widehat{\Sigma}_2^{-1}$ (the inverted Gram matrix from $\mathcal{D}_2$); (iv) SpN region: the SpN statistic using $Q = I_d$; (v) Oriented SpN region: the SpN statistic using $Q = \widehat{\Sigma}_2^{-1}$; and (iv) GaN region: the GaN statistic with $Q = I_d$. It is clear that orienting with a suitable $Q$ helps shaping the confidence region for interpretability.
}
    \label{fig:1}
\end{figure}

Another property of $Z$-functional is its affine equivariance under linear transformations. For any $Q\in\Real^{d\times d}$, $\Eb_P[\psi(Z,\theta)]=0_d$ implies that $\Eb_P[Q\psi(Z,\theta)]=0_d$. Consequently, an honest confidence set for the transformed vector $\Eb_P[Q\psi(Z,\theta)]$ can be pulled back as the preimage of the map $u\mapsto Qu$, yielding an honest confidence set for $\Eb_P[\psi(Z,\theta)]$. This is valid even if $Q$ is rank-deficient; in such cases, the resulting preimage is an affine cylinder in $\Real^d$ and does not constrain the directions orthogonal to the row space of $Q$. Nevertheless, introducing a linear transformation can be beneficial in shaping the resulting confidence set. Suitable choices of $Q$ may allow the construction of approximately rectangular confidence sets, which are both interpretable and efficient for coordinate-wise inference. See Figure~\ref{fig:1} for an illustration of how different choices of $Q$ yield different shapes; see \ref{dgp:0} for details on the data generating scheme for this figure. 

\subsection{Construction based on Self-Normalization}

Given a sample of $n$ i.i.d. observations $Z_1,\ldots,Z_n\sim P$, we adopt a sample-splitting scheme with subsets $\Dc_1=\set{Z_i:i\in\Ic_1}$ and $\Dc_2=\set{Z_i:i\in\Ic_2}$ where $\Ic_1$ and $\Ic_2$ form a partition of indices $[n]$. Let $n_1=|\Ic_1|$ and $n_2=|\Ic_2|$, and define, for $\theta\in\Theta$,\begin{align*}
    \mu(\theta) =\Eb[\psi(Z,\theta)],&\quad V(\theta)={\rm Var}(\psi(Z,\theta)),\\
    \bar\psi_1(\theta) = \frac{1}{n_1}\textstyle\sum_{i\in\Ic_1}\psi(Z_i,\theta),&\quad \hat V_1(\theta)=\frac{1}{n_1}\textstyle\sum_{i\in\Ic_1}(\psi(Z_i,\theta)-\bar\psi_1(\theta))(\psi(Z_i,\theta)-\bar\psi_1(\theta))^\top.
\end{align*} We use the subscript $\cdot_1$ to denote data-driven quantities computed using $\Dc_1$, and $\cdot_2$ to denote those computed using $\Dc_2$. Moreover, we write $\hat \cdot$ to denote a specific estimator, for e.g., $Z$-estimator, while we write $\widetilde \cdot$ to denote an arbitrary estimator.

For any $\Dc_2$-measurable matrix $\widetilde Q_2 \in \Real^{d \times d}$, the self-normalized (SeN) statistic is defined as

\begin{tcolorbox}
\begin{equation}\label{eq:def.SeNstatistics}
    T^{\rm SeN}(\theta)=T^{\rm SeN}(\theta; \widetilde Q_2) := \max_{j\in[d]}\frac{\sqrt{n_1}\abs{\ev_j^\top \widetilde Q_2 \bar\psi_1(\theta)}}{\sqrt{\ev_j^\top \widetilde Q_2 \hat V_1(\theta) \widetilde Q_2^\top\ev_j}}. \tag{SeN Statistic}
\end{equation} 
\end{tcolorbox}

This construction can be used without sample splitting by setting $\Ic_2 = \emptyset$ and choosing a non-stochastic matrix $Q$. 
In the univariate case ($d=1$), the SeN statistics converge in distribution to the standard normal law under weak conditions \citep{Gine1997}. In multivariate settings ($d\geq2$), the distribution can be approximated with the $\ell_{\infty}$-norm of a $d$-dimensional mean zero Gaussian whose variance is a correlation matrix.

Due to non-linear dependence of the SeN statistic on $\theta$, the geometry of the resulting confidence set is generally intractable. To improve tractability, we introduce a split-normalized (SpN) statistic. Given any estimator $\widetilde\theta_2\in\Theta$ obtained from $\Dc_2$, we define
\begin{tcolorbox}
    \begin{equation}\label{eq:def.SpNstatistics}
    T^{\rm SpN}(\theta)=T^{\rm SpN}(\theta;\widetilde\theta_2, \widetilde Q_2) := \max_{j\in[d]}\frac{\sqrt{n_1}\abs{\ev_j^\top \widetilde Q_2 \bar\psi_1(\theta)}}{\sqrt{\ev_j^\top \widetilde Q_2 \hat V_1(\widetilde\theta_2) \widetilde Q_2^\top\ev_j}}. \tag{SpN Statistic}
\end{equation}
\end{tcolorbox}
\noindent Unlike the SeN statistic, the SpN statistic has an explicit dependence on a preliminary estimator computed from $\Dc_2$. This preliminary estimator is used to remove the dependence on $\theta$ in the denominator of the SeN statistic.

Distributional approximations of the SeN and SpN statistics with $\theta = \theta_0(P)$ depend heavily on the relative growth of dimension with respect to the sample size, and also on the moment assumptions. For example, if $\psi(Z, \theta_0)$ is sub-Gaussian, then $d$ can be as large as $\exp(n^{1/3 - \varepsilon})$, but if $\psi(Z, \theta_0)$ has only a finite fourth moment, then $d$ can at most be $n^{1-\varepsilon}$ ~\citep[Corollary 2.1]{chernozhukov2023nearly}. Such a phase transition in the Gaussian approximation depending on moment conditions has been studied in \cite{kock2024remark}. Rather surprisingly, it is possible to consider a statistic whose distributional approximation is independent of the dimension $d$. We finally introduce the GaN statistic: 
\begin{tcolorbox}
    \begin{equation}\label{eq:def.TNstatistics}
    T^{\rm GaN}(\theta) := T^{\rm GaN}(\theta;\widetilde\theta_2, \widetilde Q_2)=\frac{\sqrt{n_1}\abs{(\theta - \widetilde\theta_2)^\top \widetilde Q_2\bar\psi_1(\theta)}}{\sqrt{(\theta - \widetilde\theta_2)^\top \widetilde Q_2\hat V_1(\widetilde\theta_2) \widetilde Q_2^\top(\theta - \widetilde\theta_2)}}. \tag{GaN Statistic}
\end{equation}
\end{tcolorbox}
\noindent When $\psi$ corresponds to the gradient of a loss $\theta\mapsto\Lc(\cdot,\theta)$, the GaN statistic is a normalized Gateaux directional derivative of $\theta\mapsto n_1^{-1}\sum_{i\in\Ic_1}\Lc(X_i,\theta)$ along the direction $\widetilde Q_2^\top(\theta-\widetilde\theta_2)$, motivating its name. The associated confidence set is then defined as
\begin{equation*}
    \widehat{\mathrm{CI}}^{\rm GaN}_\alpha:=\set{\theta\in\Theta:T^{\rm GaN}(\theta)\leq z_{\alpha/2}},
\end{equation*} where $z_{\alpha/2}$ denotes the upper $\alpha/2$-th quantile of standard normal distribution.

Self-normalization combined with and without sample splitting has been widely used in the hypothesis testing literature. \cite{Chernozhukov2019inference} use an analogue of the SeN statistic for testing in moment inequality models. \cite{kim2024dimension} developed a cross-fitted studentized statistic (with sample splitting) for dimension-agnostic inference. Although there is a duality between hypothesis tests and confidence sets, these methods yield computationally and geometrically intractable confidence sets when applied directly. To elaborate, when testing the null hypothesis $H_0: \theta_0 = \theta^*$ using a non-negative statistic $T(\theta)$, the level $\alpha$-test would reject $H_0$ if $T(\theta^*) > K_{\alpha}(\theta^*)$, for some data-driven (usually bootstrap-based) quantile $K_{\alpha}(\theta^*)$. The test-inversion based confidence set would be $\{\theta\in\Theta:\, T(\theta) \le K_{\alpha}(\theta)\}$. The computation of such a confidence set requires the computation of the quantile for each $\theta\in\Theta$, while for hypothesis testing, the computation of the quantile is only required at the hypothesized value of $\theta_0$. Hence, for tractability, we use sample splitting to compute a preliminary estimator $\tilde{\theta}_2$ from $\Dc_2$ and use that for the computation of the quantile. Note that such a substitution with a preliminary estimator is not needed for the GaN statistic because $T^{\mathrm{GaN}}(\theta_0)$ can be approximated by the folded standard Gaussian random variable. For the SeN and SpN statistics, the quantile is computed as follows: Let
\begin{equation*}\label{eq:limit_cov}
\Gamma(\theta;Q) := \mathrm{Corr}(Q V(\theta)Q^\top) \quad\mbox{for}\quad\theta\in\theta~\mbox{and}~Q\in\Real^{d\times d}.
\end{equation*} Define the method of moment estimator for this correlation matrix, obtained from $\Dc_1$, as
\begin{equation*}\label{eq:def.gammahat.iid}
    \widehat\Gamma_1(\theta;Q) ={\rm Corr}(Q\widehat V_1(\theta)Q^\top)\quad\mbox{for}\quad\theta\in\theta~\mbox{and}~Q\in\Real^{d\times d}
\end{equation*} For $\alpha\in(0,1)$, define the Gaussian multiplier bootstrap \citep{chen2018gaussian,rinaldo2019bootstrapping} quantile function $K^{\rm b}_\alpha:\Theta\to\Real$ as
\begin{equation*}\label{eq:bootstrap_quantile}
    K^{\rm b}_\alpha(\theta)=K^{\rm b}_\alpha(\theta;\widetilde Q_2):=\inf\big\{t\geq0:n^{-1}{\textstyle\sum}_{i=1}^n\mathbbm{1}(\norm{W_i^{\rm b}}_\infty>t)\leq\alpha\big\},
\end{equation*} where $W_i^{\rm b}\mid Z_1,\ldots,Z_n\sim\Nc(0_d,\widehat\Gamma_1(\theta;\widetilde Q_2))$ for $i\in[n]$. Let $\hat K^{\rm b}_\alpha :=  K^{\rm b}_{\alpha}(\widetilde\theta_2)$ evaluated only at a preliminary estimator. We then define the SeN and SpN confidence sets as
\begin{equation*}
    \widehat{\mathrm{CI}}_\alpha^{\rm SeN}:=\set{\theta\in\Theta:T^{\rm SeN}(\theta)\leq \hat K^{\rm b}_\alpha}\quad\mbox{and}\quad \widehat{\mathrm{CI}}_\alpha^{\rm SpN}:=\set{\theta\in\Theta:T^{\rm SpN}(\theta)\leq \hat K^{\rm b}_\alpha}.
\end{equation*} In the following subsections, we provide results on the miscoverage errors and diameter of the confidence sets. Throughout, we assume that $n$ is even and $n_1 = n_2 = n/2$ for convenience. This supposition is not required for proof techniques and is only made to reduce the notational burden.
\subsection{Non-asymptotic Coverage Guarantees}
In this section, we establish finite-sample coverage guarantees for the proposed confidence sets under general moment conditions and weak regularity conditions.

To describe the required moment conditions, we introduce Lyapunov ratios. For any $\theta \in \Theta$, recall $\mu(\theta)=\Eb[\psi(Z,\theta)]$, and define
\begin{equation*}
    L_q(\theta):= \sup_{u \in \Sb^{d-1}} \frac{(\Eb[|u^\top (\psi(Z, \theta) - \mu(\theta))|^q])^{1/q}}{(\Eb[|u^\top (\psi(Z, \theta) - \mu(\theta))|^2])^{1/2}}\quad\mbox{for}~q\geq2.
\end{equation*}
We have $1 \leq L_q(\theta)$ for all $\theta$ and $q\geq 2$, and we write $L_q := L_q(\theta_0)$ for the shorthand. Lyapunov ratio is closely related to $\Lc^q$–$\Lc^2$ moment equivalence for one-dimensional projections, a significant weakening of exponentially decaying tail assumptions. Moreover, this ratio is invariant under affine transformations of $\psi$, aligning with the invariance of the target $\theta_0$.


To characterize the regularity conditions, let $\Bb(\theta_0,\delta)=\set{\theta\in\Theta:\norm{\theta-\theta_0}_2\leq \delta}$ be a $\delta$-neighborhood of $\theta_0$. For $V(\theta)={\rm Var}(\psi(Z,\theta))$, we define, for $\delta>0$,
\begin{equation*}
    \phi(\delta):=\sup_{\theta\in \Bb(\theta_0,\delta)}\Norm{V(\theta_0)^{-1/2}(V(\theta)-V(\theta_0))V(\theta_0)^{-1/2}}_{\rm op}.
\end{equation*}
Moreover, we define $\gamma_q(\delta):=\sup_{\theta\in \Bb(\theta_0,\delta)}L_q(\theta)/L_q(\theta_0)$ for $q\geq2$.
\begin{remark}\label{rmk:2}
    The function $\phi$ measures the local smoothness of the variance map $\theta \mapsto V(\theta)$ around $\theta_0$. If $\psi$ is local Lipschitz in $\theta$, then $\phi(\delta) = O(\delta)$ as $\delta \to 0$. In particular, $\phi$ encodes moment-level regularity and can be controlled without the continuity of $\psi$ itself. Similarly, $\gamma_q$ captures the variation in the moment level and remains uniformly bounded around $\theta_0$, that is, $\gamma_q(\delta)=O(1)$ as $\delta\to0$, under mild conditions, such as the local continuity of $\theta\mapsto L_q(\theta)$.
\end{remark}

Theorem~\ref{thm:1} provides non-asymptotic bounds on the coverage of the SeN confidence set under minimal assumptions.

\begin{theorem}\label{thm:1}
    For $q>2$, there exists a constant $C=C(q)$, depending only on $q$, such that
    \begin{align}\label{eq:thm1:1}
        \sup_{\alpha\in(0,1)}\left(\Pb(\theta_0\notin\widehat{\mathrm{CI}}_\alpha^{\rm SeN})-\alpha\right)_+&\leq C\frac{L_4\log^{5/4}(ed)}{n^{1/4}}+ C\gamma_q(\delta)\left(\frac{L^2_{q}\log^{3}(e d) d^{2/q}}{ n^{1-2/q}}\right)^{\frac{q}{2q+2}}\nonumber\\
        &+C\log(en)\sqrt{\phi(\delta)}+\Pb(\norm{\widetilde\theta_2-\theta_0}_2\geq \delta),
    \end{align}for any $\delta>0$. Suppose further that there exists $\lambda_\circ>0$ such that
    \begin{equation}\label{eq:min_eigen}
        \Pb(\lambda_{\rm min}(\Gamma(\theta_0,\widetilde Q_2))\geq \lambda_\circ)=1.
    \end{equation} For $q\geq4$, there exists a constant $C=C(q,\lambda_\circ)$ such that for any $\delta>0$,
    \begin{align}\label{eq:thm1:2}
        \sup_{\alpha\in(0,1)}\left(\Pb(\theta_0\notin\widehat{\mathrm{CI}}_\alpha^{\rm SeN})-\alpha\right)_+&\leq  C\log^2(en)\gamma_q(\delta)\left(\frac{L^2_{q}\log^{3}(e d) d^{2/q}}{ n^{1-2/q}}\right)^{\frac{q}{q+2}\wedge\frac{q}{2(q-2)}}\nonumber\\
        &+ C\log^2(en)\phi(\delta) +\Pb(\norm{\widetilde\theta_2-\theta_0}_2\geq \delta).
    \end{align}
\end{theorem}

The general bound \eqref{eq:thm1:1} requires finite fourth moments, but is particularly useful when $L_q$ $(q\in(2,4))$ is sufficiently smaller than $L_4$. Although implicit, the quantities on the right hand side of \eqref{eq:thm1:1} and \eqref{eq:thm1:2} depend on $P$, the data generating distribution. Hence, Theorem~\ref{thm:1} implies uniform validity of SeN confidence set under mild regularity conditions by taking the supremum on both sides of \eqref{eq:thm1:1} and \eqref{eq:thm1:2} over $P\in\mathcal{P}$.

As highlighted in Remark~\ref{rmk:2}, under mild regularity conditions, we have $\phi(\delta)\to 0$ and $\gamma_q(\delta)=O(1)$ as $\delta\to 0$. Additionally, if $L_4$ is uniformly bounded and a preliminary estimator $\widetilde\theta_2$ is uniformly consistent over a distribution class $\Pc$, then the SeN confidence set is honest over $\Pc$ provided $d\log^6(en)=o(n)$. While uniform consistency may appear strong, Theorem~\ref{thm:4} establishes such a property under mild conditions; see also Theorem~\ref{thm:11}. For distribution classes allowing for higher moments, the dimension condition can be relaxed further.

The bounds in Theorem~\ref{thm:1} constitute two sources of approximation error: (i) one-sided Berry–Esséen bounds for the self-normalized sums, and (ii) quantile estimation error from the estimated correlation matrix at a preliminary estimator. The former relies on recent quantitative high-dimensional CLTs for standardized averages where the minimal eigenvalue condition \eqref{eq:min_eigen} is crucial. When this condition is satisfied, sharper rates follow from CLTs tailored to non-degenerate Gaussian approximations \citep{lopes2022central,chernozhukov2023nearly}; otherwise, the general rate in \eqref{eq:thm1:1} applies, relying only on control of the diagonal entries of the covariance matrix \citep{chernozhuokov2022improved}. In our case, the diagonals are ones, so no further assumption is needed. The quantile approximation error is controlled by Gaussian comparison inequalities \citep{lopes2022central,chernozhukov2023nearly}, and the applicability of sharper results again depends on the validity of the condition \eqref{eq:min_eigen}.

The dimension dependence in \eqref{eq:thm1:1} and \eqref{eq:thm1:2} arises from both sources of approximations. From the approximation error in Gaussian approximation (high-dimensional CLT), we use the bound: 
\begin{equation}\label{eq:1}
    \Eb^{1/q}[\max_{j\in[d]}|U_j|^q]\leq \Eb^{1/q}[{\textstyle\sum}_{j=1}^d|U_j|^q]\leq d^{1/q}L_q~\mbox{where}~U_j = \ev_j^\top V(\theta_0)^{-1/2}\psi(Z, \theta_0),~j\in[d].
\end{equation} 
While this appears to be crude, the dimension factor cannot be improved under only the finite $q$-th moment assumptions. Under a stronger distributional assumption, this can be relaxed to logarithmic dependence; for e.g., if $U_j$ admits a sub-Weibull($\beta$) tail, then $d^{1/q}$ in \eqref{eq:1} improves to $\log^{1/\beta}(ed)$ up to constants in $\beta$ and $q$; see \cite{vladimirova2020subweibull} and \cite{kuchibhotla2022moving} for details on the sub-Weibull random variables. On the other hand, the approximation error from the quantile appears via the estimation error of $\tilde{\theta}_2$ which often restricts the growth rate of $d$, unless additional structural assumptions are available. For example, with sparsity, a consistent estimator can be obtained for $\theta_0$ even if $d$ grows faster than $n$. In a similar vein, if the correlation matrix $\Gamma(\theta; Q) = \mbox{Corr}(QV(\theta)Q^{\top})$ is sparse or banded, then this structure can be used to reduce the dimension requirement of this approximation error.


Next, we establish the coverage bounds for the SpN confidence sets.

\begin{theorem}\label{thm:2} The miscoverage bound in \eqref{eq:thm1:1} applies to the confidence set $\widehat{\mathrm{CI}}_\alpha^{\rm SpN}$. Under minimum eigenvalue condition in \eqref{eq:min_eigen}, the sharper bound \eqref{eq:thm1:2} applies.
\end{theorem}

The SpN statistic shares the same ``limiting" distribution as the SeN statistic. The only modification lies in normalization, and this, in turn, is absorbed into the existing bounds in Theorem~\ref{thm:1} up to constants. Therefore, the SpN confidence set inherits the same uniform validity guarantees as the SeN confidence set, without requiring additional conditions.

We now turn to the GaN statistic. A key feature of the GaN construction is that it reduces inference to a one-dimensional problem. This simplification removes the explicit dependence on the dimension $d$ in the distributional approximation. The following result provides a two-sided coverage bound for the GaN confidence set under weaker moment conditions.
\begin{theorem}\label{thm:3}
    For $q>2$, let $\widebar q = q\wedge 4$. There exists a constant $C:=C(q)>0$ such that
    \begin{align*}
        \sup_{\alpha\in(0,1)}\Abs{\Pb(\theta_0\notin\widehat{\Cc}_{\alpha}^{\rm GaN})-\alpha}&\leq
        C\left(\frac{L_{\widebar q}^2\gamma_{\widebar q}^2(\delta)}{n^{1-2/\widebar q}}\right)^{\frac{\widebar q}{\widebar q+2}}+C\phi(\delta)+\Pb(\norm{\widetilde\theta_2-\theta_0}_2>\delta)~\mbox{for any }\delta>0.
    \end{align*}  
\end{theorem} Although the bound in Theorem~\ref{thm:3} does not explicitly depend on the dimension $d$, the tail behavior of the preliminary estimator $\widetilde\theta_2$ may still induce a dimension dependence.

\subsection{Rate of Convergence of $Z$-estimator and Width Analysis of Confidence Sets}
This section analyzes the width properties of honest confidence sets, with a focus on the SpN confidence set, which is designed to facilitate geometric tractability. While our results are stated in general form, explicit rates can be derived in more structured settings; see Section~\ref{sec:glm} for an application to generalized linear regression.

Defining $\widehat{\theta}_1$ as the $Z$-estimator from $\Dc_1$, i.e., $\bar{\psi}_1(\hat\theta_1) = 0_d$, it is clear that $\widehat{\theta}_1$ belongs to all our confidence sets. Because each of them contains $\theta_0$ with an asymptotic probability of at least $1-\alpha$, we get that, under mild regularity conditions,
\[
\liminf_{n\to\infty}\,\inf_{P\in\mathcal{P}}\,\mathbb{P}\big(\mbox{diam}_{\|\cdot\|_2}(\widehat{\mathrm{CI}}_{n,\alpha}) \ge \|\widehat{\theta}_1 - \theta_0\|_2\big) \ge 1 - \alpha.
\]
Because of this, as a preliminary but instructive step, we first analyze the rate of convergence of the $Z$-estimator.

The study of $Z$-estimators in increasing dimensions dates back to the seminal work of \cite{huber1973robust}, who analyzed $M$-estimators in linear models. This was extended by \cite{yohai1979asymptotic, portnoy1984asymptotic,portnoy1985asymptotic, portnoy1986asymptotic, portnoy1988asymptotic,mammen1993bootstrap}, who established consistency and asymptotic normality under various conditions. An approach based on stochastic equicontinuity was introduced by \cite{Pakes1989simulation} and later refined by \cite{he1996general,he2000parameters} for the convex $M$-estimation, although the conditions presented are often implicit. Recently, \cite{gauss2024asymptotics} showed consistency for $d \log(en) = o(n)$ and asymptotic normality for $d^2\log(en)=o(n)$, under weak moment assumptions and compact parameter space.
All of these works are primarily qualitative.

In the modern high-dimensional setting, the finite-sample guarantee is essential, and several works fall within the scope of a model-agnostic framework. In linear regression, \cite{oliveira2016lower} proved the $\sqrt{n/d}$-consistency of least squares under the fourth moment conditions
and this rate is minimax optimal \citep{mourtada2022exact}. \cite{kuchibhotla2020berry} established the Berry-Esséen bound for the studentized least squares, where their bound vanishes if $d^2 = \tilde{o}(n)$. For general $Z$-estimators, the key technical challenge lies in controlling local curvature, often tied to the Jacobian of the estimating equation. A relatively weak condition is self-concordance \citep{bach2010self}, which has been leveraged in finite-sample risk bounds \citep{bach2010self,ostrovskii2021finite} and deterministic approximations of estimators \citep{kuchibhotla2018deterministic}. In particular, \cite{ostrovskii2021finite} showed $\sqrt{n/d}$-consistency of estimators that minimize self-concordant loss under sub-Gaussianity and $d^2\log(en)=o(n)$. Recently, \cite{lin2024worthwhilejackknifebreakingquadratic} showed the $\sqrt{n/d}$-consistency of smooth $Z$-estimators under sub-Gaussianity and compact parameter space.


We propose the following assumptions.
\begin{enumerate}[label=\textbf{(M)},leftmargin=1cm]
    \item For $V(\theta_0)={\rm Var}(\psi(Z,\theta_0))$, there exists $\overline{\lambda}_V>0$ such that $\lambda_{\rm max}(V(\theta_0))\leq \overline{\lambda}_V<\infty$.\label{asmp:m}
\end{enumerate}



\begin{enumerate}[label=\textbf{(C\arabic*)},start=1,leftmargin=1cm]
\setlength\topsep{-1em}
    \item There exists a real-valued function $\Lc:\Zc\times \Theta\to\Real$ such that $\theta\mapsto\Lc(\cdot,\theta)$ is convex and differentiable, and $\psi(\cdot,\theta)=\nabla_\theta\Lc(\cdot,\theta)$. For $\mu(\theta)=\Eb[\psi(Z,\theta)]$, there exist $\delta_0>0$ and non-decreasing function $\underline{\mu}(\cdot):(0,\delta_0)\to\Real^+$ such that $\underline{\mu}(\delta)\leq \inf_{\eta\in\Sb^{d-1}}\eta^\top\mu(\theta_0+\delta\eta)$.\label{asmp:C.1}
    \item For $n\geq 1$, there exists $u_{1,n}:\Real^+\times(0,1)\to\Real^+$ such that, for all $\delta\in(0,\delta_0)$ and $\varepsilon\in(0,1)$,
    \begin{equation*}
        \Pb^*\bigg(\sqrt{n/2}\sup_{\eta\in\Sb^{d-1}}-\eta^\top\left(\bar\psi_1(\theta_0+\delta\eta)-\mu(\theta_0+\delta\eta)-\bar\psi_1(\theta_0)\right)\geq u_{1,n}(\delta;\varepsilon) \bigg)\leq \varepsilon.
    \end{equation*}\label{asmp:C.2}
\end{enumerate} 
\vspace{-0.5cm}
We adopt the outer probability $\Pb^*$ to account for possible measurability issues; see Section~1.2 of \cite{van1996weak} for a comprehensive discussion. 

Assumption~\ref{asmp:C.1} links the problem of $Z$-estimation with the convex $M$-estimation. We note that for $\eta = (\theta-\theta_0)/\norm{\theta-\theta_0}_2$, the difference in population objectives satisfies
\begin{align*}
    &\Eb[\Lc(X,\theta)]-\Eb[\Lc(X,\theta_0)]=\int_0^{\norm{\theta-\theta_0}_2}\eta^\top \mu(\theta_0+\delta\eta)\,d\delta\geq \int_0^{\norm{\theta-\theta_0}_2}\underline{\mu}(\delta)\,d\delta.
\end{align*} Hence, \ref{asmp:C.1} ensures that the parameter $\theta_0$ is the unique minimizer of the population objective $\theta \mapsto \Eb[\Lc(X, \theta)]$, as well as the unique root of the moment equation $\Eb[\psi(Z, \theta)] = 0_d$. Assumption \ref{asmp:C.2} ensures that the local empirical process remains controlled, providing a finite-sample analogue of the stochastic equicontinuity conditions studied in \cite{Pakes1989simulation}, \cite{he1996general}, and \cite{he2000parameters}. 

The following theorem establishes a non-asymptotic concentration of the $Z$-estimator $\hat\theta_1$. 
\begin{theorem}\label{thm:4} Assume \ref{asmp:m}, \ref{asmp:C.1}, and \ref{asmp:C.2}. For $q>2$, there exists a constant $C=C(q)>0$ such that: for any $\varepsilon>0$, let $r_n:=r_n(\varepsilon)$ satisfy the inequalities: $0<r_n<\delta_0$ and 
    \begin{equation}\label{eq:thm:4:1}
        \underline{\mu}(r_n)> \frac{u_{1,n}(r_n;\varepsilon)}{\sqrt{n}}+C\overline{\lambda}_V\bigg[\sqrt{\frac{d+\log(2/
        \varepsilon)}{n}}+\frac{L_qd^{1/2}}{n^{1-1/q}\varepsilon^{1/q}}\bigg].
    \end{equation}  Then, $\Pb(\norm{\hat\theta_1-\theta_0}_2> r_n)\leq2\varepsilon$.
\end{theorem}

Hence, the convergence rate is determined implicitly by the local curvature, via $\underline\mu(\cdot)$, and the local complexity of empirical processes through $u_{1,n}(\cdot;\cdot)$. 


To analyze the width of the SpN confidence set, we impose a spectral condition on $\widetilde Q_2$.

\begin{enumerate}[label=\textbf{(S)},start=3,leftmargin=1cm]
    \item There exist constants $\underline\sigma,\overline\sigma>0$ such that $\Pb(\underline\sigma\leq\sigma_{\rm min}(\widetilde Q_2)\leq\sigma_{\rm max}(\widetilde Q_2)\leq \overline\sigma)=1$.\label{asmp:S}
\end{enumerate}
Let $\hat\theta_2$ be a $Z$-estimator obtained from $\Dc_2$ and we denote
\begin{equation}\label{eq:SpN_Kn}
    \widehat{\mathrm{CI}}^{\rm SpN}(K_n):=\set{\theta:T^{\rm SpN}(\theta;\hat\theta_2,\widetilde Q_2)\leq K_n}\quad\mbox{for}\quad\{K_n>0:n\geq 1\}.
\end{equation} We now establish a non-asymptotic bound for the diameter of SpN confidence sets.

\begin{theorem}\label{thm:5} Suppose the assumptions made in Theorem~\ref{thm:4} are met, and assume \ref{asmp:S}. For $q\geq4$, there exists a constant $C=C(q)>0$ such that: for any $\varepsilon>0$, let $r_n:=r_n(\varepsilon)$ be such that $\Pb(\norm{\hat\theta_2-\theta_0}_2>r_n)\leq 2\varepsilon$, and let $s_n:=s_n(\varepsilon)$ satisfy that $0<s_n<\delta_0$ and 
    \begin{equation}\label{eq:thm:5:1}
        \underline{\mu}(s_n)> \frac{u_{1,n}(s_n;\varepsilon)}{\sqrt{n}}+C\frac{\overline{\sigma}^2\overline{\lambda}_V(1\vee\phi(r_n))}{\underline{\sigma}^2}\bigg[\sqrt{\frac{K_n^2d+\log(2/\varepsilon)}{n}}+\frac{\gamma_q(r_n)L_qK_nd^{1/2}}{n^{1-1/q}\varepsilon^{1/q}}\bigg].
    \end{equation}Then, $\Pb({\rm diam}_{\norm{\cdot}_2}(\widehat{\mathrm{CI}}^{\rm SpN}(K_n))> 2s_n)\leq4\varepsilon$.
\end{theorem}

The SpN confidence set contracts at a rate of $s_n$, depending on the estimation rate $r_{n}$, the quantile $K_n$ and the complexity of empirical processes $u_{1,n}$, all against curvature $\underline{\mu}$. It is easy to show that the bootstrap quantile used for our construction scales as $\sqrt{\log(ed/\alpha)}$.

\begin{lemma}\label{lem:1}
    For $\alpha\in(0,1)$, if $\sqrt{6\log(en)/n}\leq \alpha$, then for any $\Dc_1\cup\Dc_2$-measurable quantities $\widetilde\theta\in\Theta$ and $\widetilde Q\in\Real^{d\times d}$, $\Pb(K_\alpha^{\rm b}(\widetilde\theta;\widetilde Q)\leq \sqrt{2\log(4d/\alpha)})\geq 1-1/n$.
\end{lemma} 

\section{Inference for One-Dimensional Contrasts}\label{sec:contrast}
We now consider a different inferential target: $\xi^\top \theta_0$ for a non-stochastic unit vector $\xi \in \Sb^{d-1}$. This includes inference on marginal coefficients and non-sparse linear contrasts. The goal of this section is to construct honest confidence intervals for a scalar $\xi^\top\theta_0$. From the confidence set for $\theta_0$, a confidence set for $\xi^{\top}\theta_0$ readily follows by the idea of projection, i.e., 
\begin{equation*}
    \{\xi^{\top}\theta:\, \theta\in\widehat{\mathrm{CI}}_{n,\alpha}\},
\end{equation*}is an asymptotically $(1-\alpha)$ confidence set for $\xi^{\top}\theta_0$ whenever $\widehat{\mathrm{CI}}_{n,\alpha}$ is an asymptotically $(1-\alpha)$ confidence set for $\theta_0$. This projection set is, in general, extremely conservative and wide because it is simultaneously valid for all $\xi\in\Sb^{d-1}$. Hence, it is conceivable that a much shorter confidence set can be constructed for $\xi^{\top}\theta_0$.

Our constructions are based on directional self-normalization and split-normalization, respectively. Let the $Z$-estimator $\hat\theta_2$ be our preliminary estimator and consider two sets:
\begin{align}\label{eq:directional_SeN_SpN}
    \widehat{\mathscr{ C}}_{\alpha,\xi}^{\rm SeN}=\widehat{\mathscr{ C}}_{\alpha,\xi}^{\rm SeN}(\widetilde Q_2)&=\Bigg\{\theta:\frac{\sqrt{n_1}|\xi^\top \widetilde Q_2\bar\psi_1(\theta)|}{\sqrt{\xi^\top\widetilde Q_2\hat V_1(\theta)\widetilde Q_2^\top \xi}}\leq z_{\alpha/2}\Bigg\},\nonumber\\
    \widehat{\mathscr{ C}}_{\alpha,\xi}^{\rm SpN}=\widehat{\mathscr{ C}}_{\alpha,\xi}^{\rm SpN}(\hat\theta_2,\widetilde Q_2)&=\Bigg\{\theta:\frac{\sqrt{n_1}|\xi^\top \widetilde Q_2\bar\psi_1(\theta)|}{\sqrt{\xi^\top\widetilde Q_2\hat V_1(\hat\theta_2)\widetilde Q_2^\top \xi}}\leq z_{\alpha/2}\Bigg\}.
\end{align} These sets define valid confidence sets for $\theta_0$, but are, in general, unbounded. Hence, we intersect them with the SpN confidence set whose target coverage $1-\alpha/n$ converges to 1 as $n\to\infty$. Our confidence intervals are given as projections of intersected sets:
\begin{align*}
    \widehat{\rm CI}_{\alpha,\xi}^{\rm SeN}=\set{\xi^\top\theta: \theta \in \widehat{\mathscr{ C}}^{\rm SeN}_{\alpha,\xi}\cap \widehat{\mathrm{CI}}^{\rm SpN}_{\alpha/n}},\qquad\widehat{\rm CI}_{\alpha,\xi}^{\rm SpN}=\set{\xi^\top\theta: \theta \in \widehat{\mathscr{ C}}_{\alpha,\xi}^{\rm SpN}\cap \widehat{\mathrm{CI}}^{\rm SpN}_{\alpha/n}}.
\end{align*}

Theorem~\ref{thm:6} establishes the validity of these confidence intervals under mild conditions. 

\begin{theorem}\label{thm:6} Suppose assumptions in Theorem~\ref{thm:4} are met and assume \ref{asmp:S} and \eqref{eq:min_eigen}. As per Theorem~\ref{thm:4}, let $r_{n}$ be such that $\Pb(\norm{\hat\theta_2-\theta_0}_2\geq r_{n})\leq 1/n$. Define
\begin{equation*}
    \Ec_{n,q} =\gamma_q(r_{n})\log^2(en)\left(\frac{L^2_{q}\log^{3}(e d) d^{2/q}}{ n^{1-2/q}}\right)^{\frac{q}{q+2}\wedge\frac{q}{2(q-2)}}+ \log^2(en)\phi(r_{n})+\frac{1}{n}.
\end{equation*} For $q\geq 4$, there exists a constant $C=C(q)>0$ such that
\begin{equation*}
    \sup_{\xi\in\Sb^{d-1}}\sup_{\alpha\in(0,1)}\big(\Pb(\xi^\top\theta_0\notin \widehat{\rm CI}_{\alpha,\xi}^{\rm SeN})-\alpha\big)_+\leq C\Ec_{n,q}.
\end{equation*} The same miscoverage guarantee holds for $\widehat{\rm CI}_{\alpha,\xi}^{\rm SpN}$.
\end{theorem}

Thus, the SeN and SpN confidence intervals are valid under finite fourth moments, provided that $\gamma_4(r_{n})=O(1)$, $\phi(r_{n})=\tilde{o}(1)$ and $d=\tilde{o}(n)$, as $n\to\infty$. 


In the context of linear regression, we demonstrate that the proposed confidence intervals exhibit favorable width properties: they scale as $1/\sqrt{n} + d/n$, which recovers the parametric rate when $d \ll \sqrt{n}$. The width analyses of the SeN and SpN confidence intervals for the general setting are provided in Appendix A.

\section{Applications}
In this section, we apply the general results from the previous sections to linear and generalized linear regression.
\subsection{High-dimensional Misspecified Linear Model}\label{sec:lm}
Linear regression is among the most widely used statistical tools in science. Define the projection parameter \citep{buja2019models1} as the minimizer of the population least squares criterion: 
\begin{equation}\label{eq:def_projection_parameters}
    \beta_0:=~\argmin_{\theta\in\Real^d}\Eb[(Y-X^\top\theta)^2],
\end{equation} where $Y \in \Real$ is a scalar response and $X \in \Real^d$ is a random covariate vector.  Under finite second moments, $\beta_0$ solves the moment condition $\Eb[\psi((X,Y),\theta)]=0_d$ where $\psi((X,Y),\theta)=-X(Y-X^\top\theta)$, and is uniquely given by $\beta_0=\Sigma^{-1}\Eb[XY]$, if population Gram matrix $\Sigma=\Eb[XX^\top]$ is invertible. The ordinary least squares estimator (OLSE) is the conventional estimator for $\beta_0$ and takes the closed-form $\hat\beta = \hat\Sigma^{-1}n^{-1}\sum_{i=1}^nX_iY_i$ when the sample Gram matrix $\hat\Sigma = n^{-1}\sum_{i=1}^nX_iX_i^\top$ is almost surely invertible.

Under model misspecification and moment conditions, the OLSE remains consistent in $\norm{\cdot}_2$ if $d=o(n)$ \citep{oliveira2016lower}, and linear contrasts are asymptotically normal as long as $d=\tilde{o}(\sqrt{n})$ \citep{kuchibhotla2020berry}. The dimension requirement for asymptotic normality is sharp: when $d\gg \sqrt{n}$, the OLSE incurs non-negligible bias, invalidating OLSE-based Wald inference \citep{mammen1993bootstrap,lin2024worthwhilejackknifebreakingquadratic}. Recent works \citep{cattaneo2019two,chang2023inference,lin2024worthwhilejackknifebreakingquadratic} have proposed bias correction procedures to extend the validity range of Wald inference to $d=\tilde{o}(n^{2/3})$.


We consider the following assumptions throughout this section.

\begin{enumerate}[label=\textbf{(LM\arabic*)}, left = 0cm]
\item Suppose that there exists $q_x,q_y,K_x,K_{y}\geq1$ such that $q_{xy}:=(1/q_x+1/q_y)^{-1}> 2$, $\sup_{u\in\Sb^{d-1}}\Eb[\abs{u^\top\Sigma^{-1/2}X}^{q_x}]\leq K_x^{q_x}$, and $\Eb[\abs{Y-X^\top\beta_0}^{q_{y}}]\leq K_y^{q_y}$.\label{asmp:lm.1}
\item For $\Sigma = \Eb[XX^\top]$ and $V=\Eb[XX^\top(Y-X^\top\beta_0)^2]$, there exist constants $\underline{\lambda}_\Sigma$, $\overline{\lambda}_\Sigma$, and $\underline{\lambda}_V$ such that $0<\underline{\lambda}_\Sigma\leq \lambda_{\rm min}(\Sigma)\leq \lambda_{\rm max}(\Sigma)\leq \overline{\lambda}_\Sigma<\infty$ and $0<\underline{\lambda}_V\leq \lambda_{\rm min}(V)$. \label{asmp:lm.2}
\end{enumerate}

Assumption~\ref{asmp:lm.1} imposes only finite-moment assumptions on the covariates and the population residual (termed by \cite{buja2019models1}), allowing for discrete and heavy-tailed distribution. This relaxes the exponentially decaying tail assumptions, fairly common in the high-dimensional linear regression literature. Assumption~\ref{asmp:lm.2} serves as a technical condition to ensure the well-conditioned Gram matrix and non-degeneracy of the ``meat part" of the sandwich variance \citep{white1980heteroskedasticity,buja2019models1}. 

\subsubsection{Inference for Parameter Vector}
We begin by validating the confidence sets for the parameter vector. We set a preliminary estimator to $\hat\beta_2$, OLSE form $\Dc_2$, and assume $d\leq n/2$, so that the OLSE is well-defined on each split.

\begin{theorem}\label{thm:7}
    \begin{enumerate}[left=0cm]
        \item\label{thm:7.1} Assume \ref{asmp:lm.1} with $q:=\min\set{q_x/2,q_{xy}}\geq 4$ and \ref{asmp:lm.2}. 
        \begin{enumerate}[left=-.2cm]
            \item Assume \ref{asmp:S}. Then, there exists a constant $C>0$ such that
    \begin{equation}\label{eq:thm:7:1}
        \sup_{\alpha\in(0,1)}\left(\Pb(\theta_0\not\in\widehat{\mathrm{CI}}_{\alpha}^{\rm SeN})-\alpha\right)_+\leq C\Set{\frac{\log^2(en)d^{1/2}}{n^{1/2}} + \left(\frac{\log^3(e n) d^{2/q}}{n^{1-2/q}}\right)^{\frac{q}{q+2}}}.
    \end{equation} The same miscoverage bound applies to $\widehat{\mathrm{CI}}_{\alpha}^{\rm SpN}$.
    \item Let $\widetilde Q_2=\hat\Sigma_2^{-1}$ where $\hat\Sigma_2:=n_2^{-1}\sum_{i\in\Dc_2}X_iX_i^\top$. If \ref{asmp:lm.1} holds further with $q\geq 6$, then there exists a constant $C>0$ such that \begin{equation}\label{eq:thm:7:2}
        \sup_{\alpha\in(0,1)}\left(\Pb(\theta_0\not\in\widehat{\mathrm{CI}}_{\alpha}^{\rm SeN})-\alpha\right)_+\leq C\frac{\log^2(en)d^{1/2}}{n^{1/2}}.
    \end{equation}The same miscoverage bound applies to $\widehat{\mathrm{CI}}_{\alpha}^{\rm SpN}$.
    
        \end{enumerate}
    \item Assume \ref{asmp:lm.1} with $q>2$ and \ref{asmp:lm.2}. Then, there exists a constant $C>0$ such that
    \begin{equation}\label{eq:thm:7:3}        \sup_{\alpha\in(0,1)}\Abs{\Pb(\theta_0\not\in\widehat{\mathrm{CI}}_{\alpha}^{\rm GaN})-\alpha}\leq C\Set{\sqrt{\frac{d}{n}}+\left(\frac{1}{n}\right)^{\frac{(q\wedge 4)-2}{(q\wedge 4)+2}}}.
    \end{equation}
    \end{enumerate}
\end{theorem}
The bound in \eqref{eq:thm:7:1} vanishes when $d \log^6 (ed) = o(n)$ for $q = 4$ and $d \log^4 (ed) = o(n)$ for $q > 4$, as $n\to\infty$, and thereby the SeN and SpN confidence sets are honest over a class of joint distributions that admit finite moments under this dimension condition. Assumption \ref{asmp:S} that $\widetilde Q_2$ has an almost surely bounded spectrum can be relaxed to a high probability condition: for example, if $X$ is sub-Gaussian, then a matrix $\Sigma^{-\frac{1}{2}}\hat\Sigma_2\Sigma^{-\frac{1}{2}}$ has a bounded spectrum with a probability converging to 1 as long as $d = o(n_2)$ \citep{vershynin2018high}. This property fails for heavy-tailed covariates \citep{tikhomirov2018sample}, but the validity still holds, as in \eqref{eq:thm:7:2}, under a higher moment condition. The second part of Theorem~\ref{thm:7} proves a sharper bound for coverage of the GaN confidence sets under a weaker moment assumption, and the bound vanishes as long as $d=o(n)$ as $n\to\infty$.


To study width, we recall the SpN confidence set from \eqref{eq:SpN_Kn} and further consider, for $K_n\geq1$,
\begin{equation*}
    \widehat{\mathrm{CI}}^{\rm SeN}(K_n):=\set{\beta: T^{\rm SeN}(\beta;\widetilde Q_2)\leq K_n}\quad\mbox{and}\quad \widehat{\mathrm{CI}}^{\rm GaN}(K_n):=\set{\beta: T^{\rm GaN}(\beta;\hat\beta_2,\widetilde Q_2)\leq K_n}.
\end{equation*} 

\begin{theorem}\label{thm:8}
    Assume \ref{asmp:lm.1} with $q\geq 4$, \ref{asmp:lm.2}, and \ref{asmp:S}. Then, there exist constants $C_1$ and $C_2$, independent of $\underline{\sigma}, \overline{\sigma},n,$ and $d$, such that 
    \begin{equation}\label{eq:thm:8:1}
        \Pb\left({\rm diam}_{\norm{\cdot}_2}(\widehat{\mathrm{CI}}^{\rm SeN}(K_n))\geq \frac{C_1\overline{\sigma}}{\underline{\sigma}}\sqrt{\frac{d+\log(en)}{n}} \right)\leq C_2\bigg[\frac{1}{n}+\left(\frac{\overline{\sigma}^2K_n^2d\log^{3/2}(en)}{\underline{\sigma}^2n^{1-2/q_x}}\right)^{q_x/2}\bigg].
    \end{equation} The bound \eqref{eq:thm:8:1} applies to ${\rm diam}_{\norm{\cdot}_2}(\widehat{\mathrm{CI}}^{\rm SeN}(K_n))$. Moreover, if $\widetilde Q_2$ is positive definite, then the same bound \eqref{eq:thm:8:1} holds for ${\rm diam}_{\norm{\cdot}_2}(\widehat{\mathrm{CI}}^{\rm GaN}(K_n))$.
\end{theorem}

Therefore, the widths of confidence sets scale as $\sqrt{d/n}$, with a high probability, which matches the minimax estimation rate \citep{mourtada2022exact}, and is therefore unimprovable up to constants. A confidence set with a similar diameter bound was recently developed in \cite{takatsu2025bridging}. 

\subsubsection{Inference for Linear Contrasts}
In this section, we apply the procedure described in Section~\ref{sec:contrast} to construct confidence intervals for $\xi^\top\beta_0$ for a fixed $\xi\in\Sb^{d-1}$. Section~\ref{sec:contrast} employs the SpN confidence set as a bounding set due to its favorable width guarantees in general settings. In linear regression, we use the SeN confidence set with the Bonferroni quantile as our bounding set, since it attains the desired width with a better coverage guarantee. We recall $\hat{\mathscr{C}}_\alpha^{\rm SeN}$ and $\hat{\mathscr{C}}_\alpha^{\rm SpN}$ from \eqref{eq:directional_SeN_SpN}, and define
\begin{align*}
    &\widehat{\rm CI}_\alpha^{\rm SeN}(\xi^\top\beta_0)=\Set{\xi^\top\theta: \theta \in \hat{\mathscr{C}}_\alpha^{\rm SeN}(\hat\Sigma_2^{-1})~\mbox{and}~ T^{\rm SeN}(\theta;\hat\Sigma_2^{-1})\leq z_{\frac{1}{2nd}}},\\
    &\widehat{\rm CI}_\alpha^{\rm SpN}(\xi^\top\beta_0)=\Set{\xi^\top\theta: \theta \in \hat{\mathscr{C}}_\alpha^{\rm SpN}(\hat\beta_2,\hat\Sigma_2^{-1})~\mbox{and}~ T^{\rm SeN}(\theta;\hat\Sigma_2^{-1})\leq z_{\frac{1}{2nd}}},
\end{align*} where the sample precision matrix $\hat\Sigma_2^{-1}$ serves as the rotation matrix. We note that the construction of the SeN confidence interval here does not require preliminary estimators.

\begin{theorem}\label{thm:9}
    Suppose that Assumption~\ref{asmp:lm.2} holds. 
    \begin{enumerate}[left=0cm]
        \item Assume \ref{asmp:lm.1} and let $s=q_{xy}\wedge 3$, then there exists a constant $C>0$ such that
        \begin{equation}\label{eq:thm:9:1}
        \sup_{\alpha\in(0,1)}\left(\Pb\left(\xi^\top\beta_0\notin\widehat{\rm CI}^{\rm SeN}_\alpha(\xi^\top\beta_0)\right)-\alpha\right)_+\leq C\frac{\log^{s/2}(end)}{n^{s/2-1}}.        
    \end{equation}
    \item Assume \ref{asmp:lm.1} with $q>2$, then there exists a constant $C>0$ such that
    \begin{equation}\label{eq:thm:9:2}        \sup_{\alpha\in(0,1)}\left(\Pb\left(\xi^\top\beta_0\notin\widehat{\rm CI}_\alpha^{\rm SpN}(\xi^\top\beta_0)\right)-\alpha\right)_+\leq C\left[\sqrt{\frac{d}{n}}+\left(\frac{\log(en)}{n}\right)^{1-2/(q\wedge 4)}\right].
    \end{equation}
    \end{enumerate}
\end{theorem}
The SeN confidence interval is honest over a distribution class that admits finite $(2+\delta)$-th moments as long as $d\leq n/2$ and $n\to\infty$. The SpN confidence interval is honest over a class that further allows a higher moment of covariates, provided that $d=o(n)$. 

We next analyze their width properties. The standard distribution approximation for OLSE under model misspecification \citep{buja2019models1} suggests $\sqrt{n}(\hat\beta_1-\beta_0)\approx \Nc(0_d,\Sigma^{-1}V\Sigma^{-1})$. This approximation leads to a canonical Wald interval for $\xi^\top\theta_0$:
\begin{equation*}
    \widehat{\rm CI}^{\rm Wald}_\alpha(\xi^\top\beta_0) = \big[\xi^\top\hat\beta_1\pm z_{\alpha/2}\sqrt{\xi^\top\Sigma^{-1}V\Sigma^{-1}\xi/n}\big].
\end{equation*}
To formalize this geometric similarity, we consider the one-side Hausdorff distance between intervals $A,B\subseteq\Real$, defined as
\begin{equation*}
    {\rm Haus}_{|\cdot|}(A\to B):=\sup_{a\in A}\inf_{b\in B} |a-b|.
\end{equation*}
Theorem~\ref{thm:10} reports the second-order width analyses for the SeN and SpN confidence intervals toward the Wald's, in terms of the Hausdorff distance.

\begin{theorem}\label{thm:10} Assume \ref{asmp:lm.1} with $q\geq 4$ and \ref{asmp:lm.2}. There exists constants $C_1,C_2>0$ such that
\begin{equation*}
    \Pb\left(\sqrt{n/2}~{\rm Haus}_{|\cdot|}\left(\widehat{\rm CI}^{\rm SeN}_\alpha(\xi^\top\beta_0)\to\widehat{\rm CI}^{\rm Wald}_\alpha(\xi^\top\beta_0)\right)\geq \frac{C_1d}{\sqrt{n}}\right)\leq C_2\bigg[\frac{1}{n}+\left(\frac{d\log^3(en)}{n^{1-4/q_x}}\right)^{q_x/4}\bigg].
\end{equation*} The same bound applies to ${\rm Haus}_{|\cdot|}(\widehat{\rm CI}^{\rm SpN}_\alpha(\xi^\top\beta_0)\to\widehat{\rm CI}^{\rm Wald}_\alpha(\xi^\top\beta_0))$.
\end{theorem} Thus, the widths of the SeN and SpN confidence intervals scale as $1/\sqrt{n} + d/n$ with high probability, matching the scaling of $|\xi^\top(\hat\beta_1 - \beta_0)|$ under model misspecification \citep{mammen1993bootstrap,chang2023inference}. Consequently, any valid confidence interval containing OLSE cannot improve on this rate up to constants. When $d = o(\sqrt{n})$, the SeN and SpN intervals converge to the Wald interval faster than the Wald interval shrinks to a point, with widths comparable to the full-sample Wald interval up to a factor of $\sqrt{2}$. Notably, while the Wald interval loses validity when $d \gg \sqrt{n}$ \citep{mammen1993bootstrap,kuchibhotla2020berry}, the SeN and SpN intervals remain valid under the weaker conditions $d \leq n/2$ and $d = o(n)$ as $n\to\infty$, respectively.

\subsection{High-dimensional Misspecified Generalized Linear Model}\label{sec:glm}

We now extend to generalized linear models (GLMs). Consider the canonical form of a loss:
\begin{equation*}
    \Lc((y,x),\theta)=\ell(x^\top\theta)-yx^\top\theta,
\end{equation*} where $\ell: \mathbb{R} \to \mathbb{R}$ is convex and thrice differentiable. The gradient and Hessian of the loss function are given by $\psi((y,x),\theta)=\nabla_\theta\Lc((y,x),\theta)=-x^\top(y- \ell'(x^\top\theta))$ and $\nabla_\theta^2\Lc((y,x),\theta)=xx^\top \ell''(x^\top\theta)\succeq 0$. In this setup, the target $\theta_0$ is the solution to the following moment equation:
\begin{equation}\label{eq:$Z$-functional_glm}
    \theta_0~:~\Eb[X(Y-\ell'(X^\top\theta_0))]=0_d.
\end{equation}
We reemphasize that $\theta_0$ in \eqref{eq:$Z$-functional_glm} is defined as a statistical functional and does not require the data-generating process to be a parametric GLM.


We impose the following assumptions throughout this section.
\begin{enumerate}[label=\textbf{(GLM\arabic*)}, leftmargin=2cm]
    \item There exists a positive constant $b_1>0$ such that $\sup_{x}\abs{\ell''(x)}\vee\sup_{x}\abs{\ell'''(x)}\leq b_1$.\label{asmp:glm1}
    \item Suppose that there exists $q_x,q_y,K_x,K_{y}\geq1$ such that $q_{xy}:=(1/q_x+1/q_y)^{-1}> 2$, $\sup_{u\in\Sb^{d-1}}\Eb[\abs{u^\top\Sigma^{-1/2}X}^{q_x}]\leq K_x^{q_x}$, and $\Eb[\abs{Y-\ell'(X^\top\theta_0)}^{q_{y}}]\leq K_y^{q_y}$.\label{asmp:glm2}
    \item For $\Sigma = \Eb[XX^\top]$ and $V = \Eb[XX^\top\ell''(X^\top\theta_0)]$, there exist constants $\underline{\lambda}_\Sigma,\overline{\lambda}_\Sigma, \underline\lambda_V\in(0,\infty)$ such that $0<\underline{\lambda}_\Sigma\leq \lambda_{\rm min}(\Sigma)\leq \lambda_{\rm max}(\Sigma)\leq \overline{\lambda}_\Sigma<\infty$ and $\lambda_{\rm min}(V)\geq \underline\lambda_V>0$.\label{asmp:glm3}
\end{enumerate} 

Assumption~\ref{asmp:glm1} encompasses a wide range of commonly used loss functions, including squared loss, logistic loss, probit loss, and several smoothed robust proxies, such as smoothed hinge or bisquare loss. Assumption~\ref{asmp:glm2} generalizes \ref{asmp:lm.1} to the GLM setting, requiring finite moments of standardized covariates and population residuals. Assumption~\ref{asmp:glm3} ensures the well-conditioned Gram matrix and the Hessian matrix. Importantly, \ref{asmp:glm3} is a pointwise condition at $\theta_0$ and does not impose uniform strong convexity in a neighborhood.



We first establish a non-asymptotic concentration bound for $\hat\theta_1$ under mild assumptions.
\begin{theorem}\label{thm:11} Assume \ref{asmp:glm1}, \ref{asmp:glm2}, and \ref{asmp:glm3}. Then, there exist constants $c>0$ and $C>0$ such that for any $\varepsilon\in(0,1)$, with probability at least $1-\varepsilon$,
\begin{equation}\label{eq:thm:10:1}
   C^{-1}\norm{\hat\theta_1-\theta_0}_2\leq \sqrt{\frac{d+\log(2/\varepsilon)}{n}}+\frac{d^{1/2}}{n^{1-1/q_{xy}}\varepsilon^{1/q_{xy}}},
\end{equation}provided that the right-hand side is no greater than $c$.
\end{theorem}
Therefore, $Z$-estimators for misspecified GLMs are $\sqrt{n/d}$-consistent under $(2+\epsilon)$-th moments if $d = o(n)$. To the best of our knowledge, this is the first non-asymptotic result that accommodates heavy-tailed covariates and responses, unbounded parameter spaces, and model misspecification. The rate comprises a sub-Gaussian term with a polynomially decaying remainder, reflecting a structure analogous to classical Fuk–Nagaev inequalities \citep{einmahl2008characterization, rio2017constants}. As per Theorem~\ref{thm:11}, the $Z$-estimator $\hat\theta_2$ from $\Dc_2$ also has the convergence rate of $\sqrt{d/n}.$ The next result establishes that $\hat\theta_2$ serves as a valid preliminary estimator for constructing uniformly valid confidence sets.
\begin{theorem}\label{thm:12}
Assume \ref{asmp:glm1} and \ref{asmp:glm3}. Let $q=\min\set{q_x/2,q_{xy}}$.
    \begin{enumerate}[left=0cm]
        \item Assume \ref{asmp:glm2} with $q \geq 4$ and \ref{asmp:S}. Then, there exists a constant $C>0$ such that
    \begin{equation*}
        \sup_{\alpha\in(0,1)}\left(\Pb(\theta_0\not\in\widehat{\mathrm{CI}}_{\alpha}^{\rm SeN})-\alpha\right)_+\leq C\bigg\{\frac{\log^2(en)d^{1/2}}{n^{1/2}} + \left(\frac{\log^3(e n) d^{2/q}}{n^{1-2/q}}\right)^{\frac{q}{q+2}}\bigg\}.
    \end{equation*} Same bound applies to the uniform miscoverage rate of $\widehat{\mathrm{CI}}_{\alpha}^{\rm SpN}$.
    \item\label{thm:11.2} Assume \ref{asmp:glm2} with $q>2$, then there exists a constant $C>0$ such that
    \begin{equation*}      \sup_{\alpha\in(0,1)}\abs{\Pb(\theta_0\not\in\widehat{\mathrm{CI}}_{\alpha}^{\rm GaN})-\alpha}\leq C\bigg\{\sqrt{\frac{d}{n}}+\left(\frac{1}{n}\right)^{\frac{(q\wedge 4)-2}{(q\wedge 4)+2}}\bigg\}.
    \end{equation*}
    \end{enumerate}
\end{theorem}
Therefore, the SeN, SpN, and GaN confidence sets are honest over a distribution class that admits finite moments as long as $d=\tilde{o}(n)$. For SeN and SpN, Assumption~\ref{asmp:S} can be relaxed to a high-probability condition.

We recall $\widehat{\mathrm{CI}}^{\rm SpN}(K_n)$ from \eqref{eq:SpN_Kn} and Theorem~\ref{thm:13} establishes the bound for its diameter.

\begin{theorem}\label{thm:13}
    Assume \ref{asmp:glm1}, \ref{asmp:glm2} with $q\geq 4$, \ref{asmp:glm3}, and \ref{asmp:S}. Then, there exist constants $c,C>0$ such that for any $\varepsilon\in(0,1)$, 
    \begin{equation*}
        C^{-1}{\rm diam}_{\norm{\cdot}_2}\left(\widehat{\mathrm{CI}}^{\rm SpN}(K_n)\right)\leq \sqrt{\frac{K_n^2d+\log(2/\varepsilon)}{n}}+\frac{d^{1/2}}{n^{1-1/q_{xy}}\varepsilon^{1/q_{xy}}}\,
    \end{equation*} holds with probability at least $1-\varepsilon$, provided that the right-hand side is no greater than $c$.
\end{theorem}

Since the quantile used in the construction of the SpN set satisfies $\hat K_\alpha^{\rm b} \asymp \sqrt{\log(ed/\alpha)}$, the diameter of the set scales as $\sqrt{d \log(ed)/n}$ across a family of generalized linear regressions.

\section{Numerical Studies}

This section evaluates the finite sample performance of the proposed confidence sets in linear regression across a range of data-generating processes (dgps), chosen to reflect both well-specified and misspecified scenarios.

\begin{enumerate}[label=\textbf{(dgp\arabic*)},start=1,left=0cm]
    \item \textbf{Well-specified linear model with independent covariates.} Let $(X_i,Y_i)\sim(X,Y)$ for $i\in[n]$ where $Y= X^\top\beta_0+\varepsilon_i$. Here, the covariate $X$ and the error $\varepsilon$ are independently drawn from $X\sim \Nc(0_d, I_d)$ and $\varepsilon\sim\Nc(0,1)$. We let $\beta_0 = 1_d/\sqrt{d}$.
    \label{dgp:0}
    \item \textbf{Mis-specified linear model with dependent covariates.} The covariates are defined as $\ev_j^\top X_i=U_{ij}Z_j$ for $i\in[n]$ and $j\in[d]$ where $U_{ij}\sim\Nc(0,1)$ and $Z_i\sim{\rm unif}[1,3]$. The responses were put as $Y_i = Q_i +0.25Q_iV_i+V_i$ where $Q_i=\norm{X_i}_2^2/d$ and $V_i$ is drawn from a mixture of Gaussian: $1/2\Nc(1/2, (1.2)^2)+1/2\Nc(-1/2, (0.7)^2)$ for $i\in[n]$. The variables $U_{ij},Z_i,V_i$ are independent.\label{dgp:2}
\end{enumerate}
The sample size is fixed at $n=2000$ with a split ratio of $1/2$, and the number of covariates varies as $d = \lfloor n^\gamma \rfloor$ for $\gamma \in [0.15, 0.8]$. For \textbf{(dgp2)}, the true projection parameters are approximated by Monte Carlo simulation with $10^5$ replicates.
The effective width of a confidence set $\mathcal{A} \subseteq \mathbb{R}^d$ is defined as ${\rm vol}(\mathcal{A})^{1/d}$, serving as an optimistic proxy for its geometric diameter. Closed-form expressions are available for all methods except \texttt{SeN}, \texttt{O-SeN}, and \texttt{GaN}, whose volumes are approximated numerically. The GaN set remains bounded with probability one, whereas the SeN sets may be unbounded with a positive probability. Therefore, the latter sets are intersected with the SpN set whose target coverage is $1 - 1/n$.

\begin{figure}[h]
    \centering
    \begin{subfigure}[b]{0.49\textwidth}
        \centering
        \includegraphics[width=\textwidth]{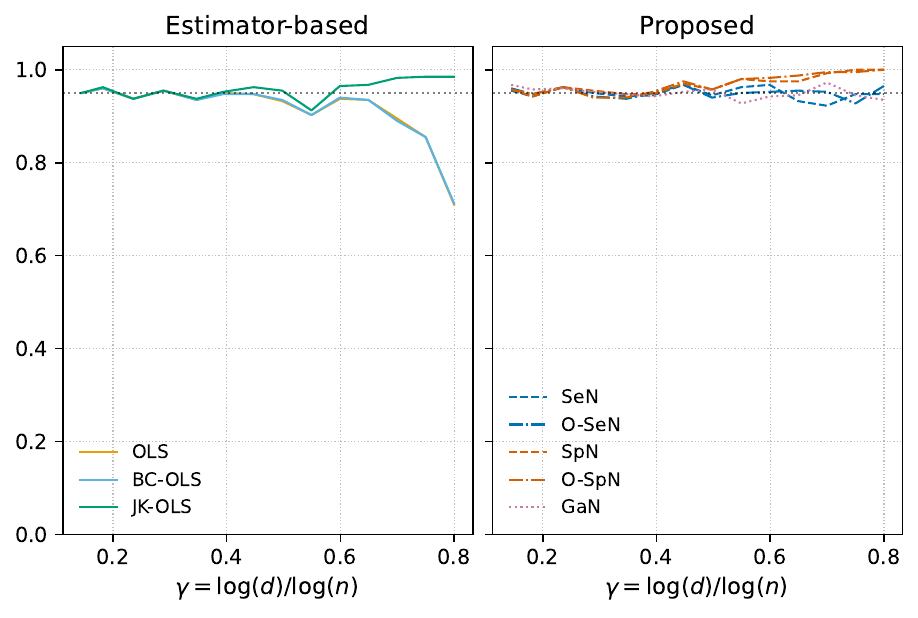}
        \caption{Empirical Coverage for \textbf{(dgp1)}}
    \end{subfigure}%
    ~ 
    \begin{subfigure}[b]{0.49\textwidth}
        \centering
        \includegraphics[width=\textwidth]{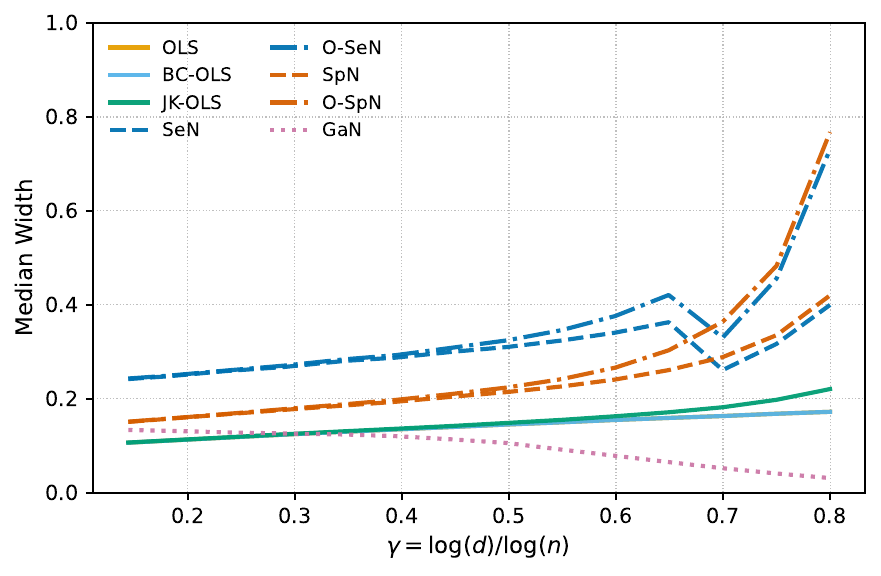}
        \caption{Effective width for \textbf{(dgp1)}}
    \end{subfigure}
\\
    \centering
    \begin{subfigure}[b]{0.49\textwidth}
        \centering
        \includegraphics[width=\textwidth]{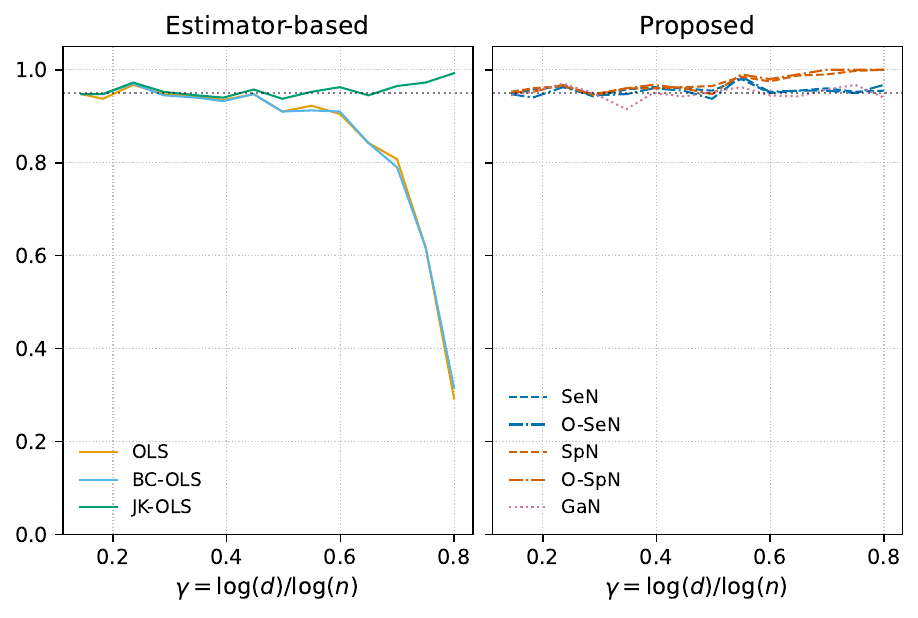}
        \caption{Empirical Coverage for \textbf{(dgp2)}}
    \end{subfigure}%
    ~ 
    \begin{subfigure}[b]{0.49\textwidth}
        \centering
        \includegraphics[width=\textwidth]{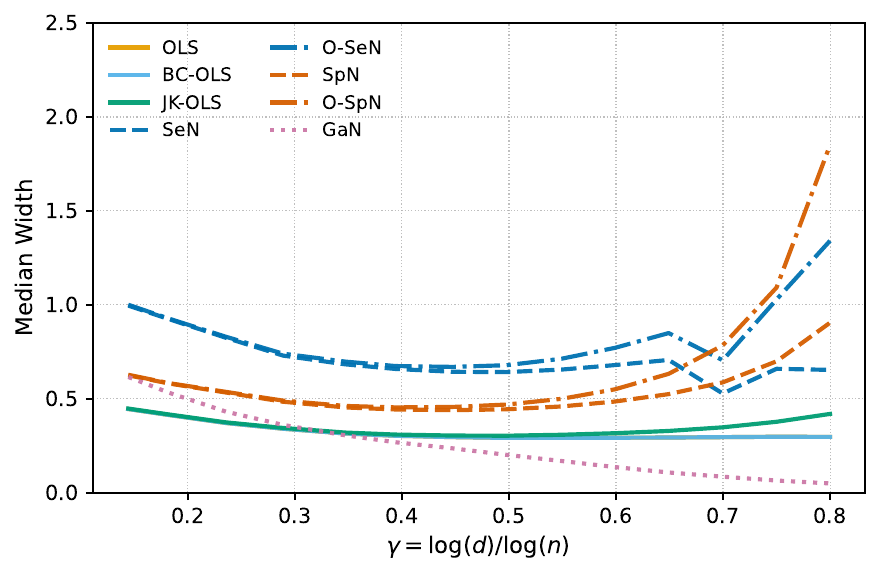}
        \caption{Effective width for \textbf{(dgp2)}}
    \end{subfigure}
    \caption{Coverage probabilities and median effective width of 95\% confidence sets for the full vector of projection parameter under three data-generating processes, plotted against $\log(d)/\log(n)$. Each row corresponds to a different data-generating process. Methods include classical Wald-type intervals, bias-corrected OLS estimators, and procedures based on self-normalization (SN), split-normalization (SpN), and Gateaux-Normalization (GaN). See Table~\ref{tab:abbreviations} for method abbreviations.}
    \label{fig.2}
\end{figure}

The empirical coverages and median effective width of the confidence sets for the parameter vector are reported in Figure~\ref{fig.2}. As anticipated, the classical Wald procedure suffers from severe under-coverage when $d \gg \sqrt{n}$, consistent with the theoretical threshold \citep{mammen1993bootstrap, chang2023inference, lin2024worthwhilejackknifebreakingquadratic}. Jackknife correction shows better coverage, but tend to be conservative. In contrast, the SeN and GaN sets attain near-nominal levels uniformly across settings, while the SpN sets tend to be conservative.

\begin{table}[h]
\centering
\caption{Abbreviations for inference methods.}
\begin{tabular}{ll}
\toprule
\textbf{Abbreviation} & \textbf{Method Description} \\
\midrule
OLS          & Ordinary Least Squares with Sandwich Variance \\
BC-OLS       & Bias-corrected OLSE \citep{chang2023inference} with Sandwich Variance\\
JK-OLS       & Jackknife Correction for both OLSE and Variance \citep{cattaneo2019two} \\
SeN        & Self-Normalization with $\widetilde Q_2 = I_d$ \\
O-SeN      & Oriented Self-Normalization with $\widetilde Q = \hat\Sigma_2^{-1}$ \\
SpN        & Split-Normalization with $\widetilde Q_2 = I_d$ \\
O-SpN      & Oriented Split-Normalization with $\widetilde Q = \hat\Sigma_2^{-1}$\\
GaN      & Gateaux-Normalization with $\widetilde Q_2 = I_d$\\
\bottomrule
\end{tabular}
\label{tab:abbreviations}
\end{table}

In terms of set size, sample splitting introduces scaling $\sqrt{2}$ and the effective radii generally increase with dimension. Since \texttt{SeN} and \texttt{O-SeN} sets are intersected with the superset, their volumes tend to be inflated, particularly in lower dimensions, although an apparent ``elbow" phenomenon suggests that the measured volumes eventually approximate those of the original constructions as the dimension increases. Notably, the \texttt{GaN} confidence sets exhibit substantially smaller effective radii than other methods, including the Wald procedure (baseline). However, this should be interpreted with care. The \texttt{GaN} sets are elliptical, and the volume of elliptical (or spherical) regions shrinks rapidly with dimension compared to axis-aligned hyperrectangles. Thus, this behavior is attributable to the geometry of the set rather than to the inferential efficiency alone.

\begin{figure}
    \centering
    \begin{subfigure}[b]{0.48\textwidth}
        \centering
        \includegraphics[width=\textwidth]{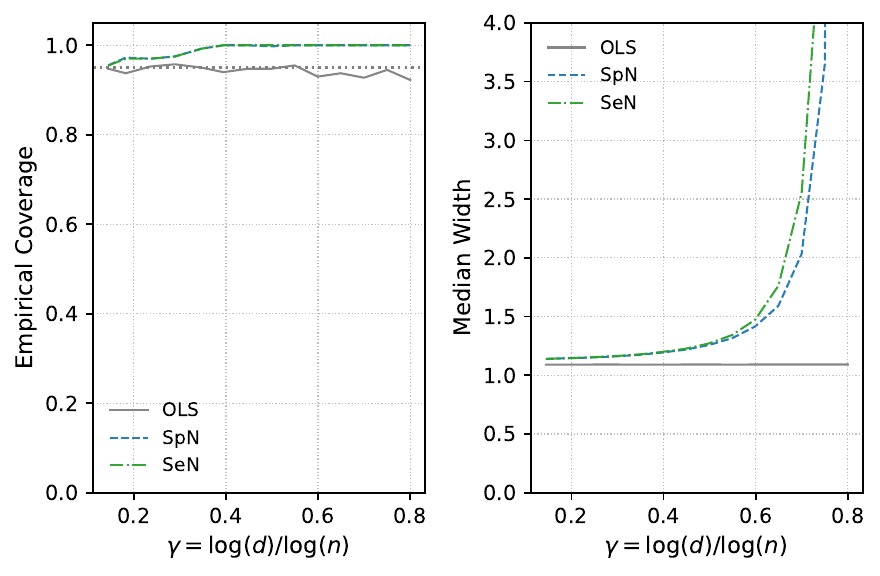}
        \caption{Empirical Coverage and Width for \textbf{(dgp1)}}
    \end{subfigure}
    \begin{subfigure}[b]{0.48\textwidth}
        \centering
        \includegraphics[width=\textwidth]{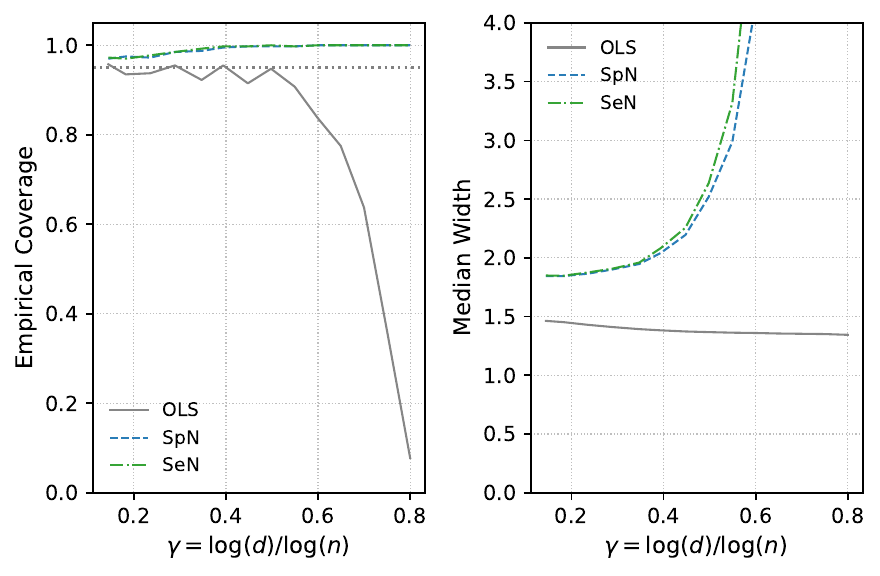}
        \caption{Empirical Coverage and Width for \textbf{(dgp2)}}
    \end{subfigure}

    \caption{Empirical coverage probabilities and median widths of 95\% confidence intervals for the linear contrast $\ev_1^\top\beta_0$. Methods compared include the classical Wald intervals based on OLS and sandwich variance estimation, along with the self-normalization and split-normalization procedures proposed in Section~4.3. Results are reported across a range of dimension-to-sample-size log ratios $\gamma = \log(d)/\log(n)\in[0.15, 0.8]$ for each data-generating process.}

    \label{fig.3}
\end{figure}

We now evaluate the finite sample performance of the proposed confidence intervals for one-dimensional linear functionals of the regression coefficient vector, focusing on the contrast $\xi^\top\beta_0$ with $\xi = 1_d/\sqrt{d}$. Figure~\ref{fig.3} reports the empirical coverage and median width of the constructed intervals across a range of dimension-to-sample-size ratios, parametrized by $\gamma = \log(d)/\log(n) \in [0.15, 0.8]$, for two dgps introduced earlier.

Across all settings, the proposed methods exhibit valid, though conservative, coverage performance throughout the examined dimensional regimes. In contrast, the classical Wald intervals display varying behavior depending on the structure of the data-generating process. Under \textbf{(dgp1)}, where both linearity and covariate independence hold, the Wald procedure maintains nominal coverage up to moderate dimensions, consistent with previous findings in \cite{cattaneo2018alternative, cattaneo2018inference}. However, in \textbf{(dgp2)}, where both linearity and independence are absent, the Wald intervals show severe undercoverage when $\gamma \gtrsim 0.5$. In terms of inferential efficiency, the proposed procedures yield widths comparable to those of Wald intervals for moderate dimensions ($\gamma \lesssim 0.5$). This behavior aligns with our theoretical findings in Theorem~\ref{thm:10}. However, beyond this, the width of our methods grows more noticeably, reflecting the high price for maintaining uniform validity in high dimensions.

\section{Conclusions}
This paper develops a general framework for constructing uniformly valid confidence sets for functionals defined by moment equations, adopting a model-agnostic, assumption-lean approach that requires only weak moment and regularity conditions.

The methodological contribution lies in the introduction of self-normalized (SeN), split-normalized (SpN), and Gateaux-normalized (GaN) statistics. These are constructed through a combination of test inversion and sample splitting, enabling both statistical robustness and computational tractability. We show that the resulting confidence sets are honest over broad distribution classes and admit worst-case width guarantees. In particular, the proposed procedures remain valid under heavy-tailed distributions, model misspecification, and unbounded parameter spaces.

When applied to high-dimensional linear and generalized linear regression, our theory yields the non-asymptotic width and coverage guarantees for confidence sets that are robust to model and distribution misspecification. We establish that the SpN sets attain uniform coverage with diameters scaling as $\sqrt{d\log(ed)/n}$ for generalized linear regression. In linear regression, the SeN, SpN, and GaN sets achieve honesty with the minimax optimal diameter scaling $\sqrt{d/n}$. Furthermore, we show that directional inference procedures based on these sets yield valid confidence intervals for linear contrasts, with second-order width analysis toward the classical Wald interval. Numerical studies support our theoretical findings, demonstrating that the proposed methods achieve reliable coverage across a range of nonlinear data-generating processes. Compared with classical procedures, the SeN, SpN, and GaN constructions exhibit better robustness and comparable width behavior in high-dimensional settings.

Future directions include extending the framework to settings with dependent data, such as time series or spatial models, as well as to partially identified parameters defined through moment inequalities. The combination of geometric tractability, finite-sample validity, and distributional robustness makes the proposed approach a promising foundation for honest inference in modern high-dimensional applications.
\bibliographystyle{abbrvnatnew}
\bibliography{bib}

\begin{thebibliography}{96}
\providecommand{\natexlab}[1]{#1}
\providecommand{\url}[1]{\texttt{#1}}
\expandafter\ifx\csname urlstyle\endcsname\relax
  \providecommand{\doi}[1]{doi: #1}\else
  \providecommand{\doi}{doi: \begingroup \urlstyle{rm}\Url}\fi

\bibitem[Andrews and Guggenberger(2009)]{AndrewsGuggenberger2009}
D.~W. Andrews and P.~Guggenberger.
\newblock Validity of subsampling and plug-in asymptotic inference for
  parameters defined by moment inequalities.
\newblock \emph{Econometric Theory}, 25\penalty0 (3):\penalty0 669--709, 2009.

\bibitem[Andrews and Guggenberger(2010)]{andrews2010asymptotic}
D.~W.~K. Andrews and P.~Guggenberger.
\newblock Asymptotic size and a problem with subsampling and with the \(m\) out
  of \(n\) bootstrap.
\newblock \emph{Econometric Theory}, 26\penalty0 (2):\penalty0 426--468, 2010.

\bibitem[Andrews and Mikusheva(2016)]{andrews2016conditional}
D.~W.~K. Andrews and A.~Mikusheva.
\newblock Conditional inference with a functional nuisance parameter.
\newblock \emph{Econometrica}, 84\penalty0 (4):\penalty0 1571--1612, 2016.

\bibitem[Andrews and Soares(2010)]{andrews2010inference}
D.~W.~K. Andrews and G.~Soares.
\newblock Inference for parameters defined by moment inequalities using
  generalized moment selection.
\newblock \emph{Econometrica}, 78\penalty0 (1):\penalty0 119--157, 2010.

\bibitem[Bach(2010)]{bach2010self}
F.~Bach.
\newblock Self-concordant analysis for logistic regression.
\newblock \emph{Electronic Journal of Statistics}, 4:\penalty0 384--414, 2010.

\bibitem[Belloni et~al.(2015)Belloni, Chernozhukov, and
  Kato]{Belloni2015Uniform}
A.~Belloni, V.~Chernozhukov, and K.~Kato.
\newblock Uniform post-selection inference for least absolute deviation
  regression and other \({Z}\)-estimation problems.
\newblock \emph{Biometrika}, 102\penalty0 (1):\penalty0 77--94, 2015.

\bibitem[Bentkus and G{\"o}tze(1996)]{bentkus1996berry}
V.~Bentkus and F.~G{\"o}tze.
\newblock The {Berry-Esseen} bound for {Student's} statistic.
\newblock \emph{Annals of Probability}, 24\penalty0 (1):\penalty0 491--503,
  1996.

\bibitem[Bertail et~al.(1999)Bertail, Politis, and
  Romano]{bertail1999subsampling}
P.~Bertail, D.~N. Politis, and J.~P. Romano.
\newblock On subsampling estimators with unknown rate of convergence.
\newblock \emph{Journal of the American Statistical Association}, 94\penalty0
  (446):\penalty0 569--579, 1999.

\bibitem[Berthet and Einmahl(2022)]{einmahl2022cube}
P.~Berthet and J.~H.~J. Einmahl.
\newblock Cube root weak convergence of empirical estimators of a density level
  set.
\newblock \emph{Annals of Statistics}, 50\penalty0 (3):\penalty0 1423--1446,
  2022.

\bibitem[Bickel et~al.(1997)Bickel, G{\"o}tze, and van
  Zwet]{bickel1997resampling}
P.~J. Bickel, F.~G{\"o}tze, and W.~R. van Zwet.
\newblock Resampling fewer than \(n\) observations: Gains, losses, and remedies
  for losses.
\newblock \emph{Statistica Sinica}, 7\penalty0 (1):\penalty0 1--31, 1997.

\bibitem[Boucheron et~al.(2013)Boucheron, Lugosi, and
  Massart]{Boucheron2013concentration}
S.~Boucheron, G.~Lugosi, and P.~Massart.
\newblock \emph{Concentration inequalities}.
\newblock Oxford University Press, Oxford, 2013.
\newblock A nonasymptotic theory of independence, With a foreword by Michel
  Ledoux.

\bibitem[Brailovskaya and van Handel(2024)]{brailovskaya2023universality}
T.~Brailovskaya and R.~van Handel.
\newblock Universality and sharp matrix concentration inequalities.
\newblock \emph{Geometric and Functional Analysis}, 34\penalty0 (4):\penalty0
  1101--1186, 2024.

\bibitem[Buja et~al.(2019{\natexlab{a}})Buja, Brown, Berk, George, Pitkin,
  Traskin, Zhang, and Zhao]{buja2019models1}
A.~Buja, L.~Brown, R.~Berk, E.~George, E.~Pitkin, M.~Traskin, K.~Zhang, and
  L.~Zhao.
\newblock Models as approximations {I}: Consequences illustrated with linear
  regression.
\newblock \emph{Statistical Science}, 34\penalty0 (4):\penalty0 523--544,
  2019{\natexlab{a}}.

\bibitem[Buja et~al.(2019{\natexlab{b}})Buja, Brown, Kuchibhotla, Berk, George,
  and Zhao]{buja2019models2}
A.~Buja, L.~Brown, A.~K. Kuchibhotla, R.~Berk, E.~George, and L.~Zhao.
\newblock Models as approximations {II}: A model-free theory of parametric
  regression.
\newblock \emph{Statistical Science}, 34\penalty0 (4):\penalty0 545--565,
  2019{\natexlab{b}}.

\bibitem[Cattaneo et~al.(2018{\natexlab{a}})Cattaneo, Jansson, and
  Newey]{cattaneo2018alternative}
M.~D. Cattaneo, M.~Jansson, and W.~K. Newey.
\newblock Alternative asymptotics and the partially linear model with many
  regressors.
\newblock \emph{Econometric Theory}, 34\penalty0 (2):\penalty0 277--301,
  2018{\natexlab{a}}.

\bibitem[Cattaneo et~al.(2018{\natexlab{b}})Cattaneo, Jansson, and
  Newey]{cattaneo2018inference}
M.~D. Cattaneo, M.~Jansson, and W.~K. Newey.
\newblock Inference in linear regression models with many covariates and
  heteroscedasticity.
\newblock \emph{Journal of the American Statistical Association}, 113\penalty0
  (523):\penalty0 1350--1361, 2018{\natexlab{b}}.

\bibitem[Cattaneo et~al.(2019)Cattaneo, Jansson, and Ma]{cattaneo2019two}
M.~D. Cattaneo, M.~Jansson, and X.~Ma.
\newblock Two-step estimation and inference with possibly many included
  covariates.
\newblock \emph{Review of Economic Studies}, 86\penalty0 (3):\penalty0
  1095--1122, 2019.

\bibitem[Cattaneo et~al.(2020)Cattaneo, Jansson, and Ma]{cattaneo2020practical}
M.~D. Cattaneo, M.~Jansson, and X.~Ma.
\newblock Bootstrap-based inference for cube root consistent estimators.
\newblock \emph{Econometrica}, 88\penalty0 (3):\penalty0 1159--1194, 2020.

\bibitem[Chang et~al.(2024)Chang, Kuchibhotla, and Rinaldo]{chang2023inference}
W.~Chang, A.~K. Kuchibhotla, and A.~Rinaldo.
\newblock Inference for projection parameters in linear regression: Beyond \(d
  = o(n^{1/2})\).
\newblock \emph{arXiv preprint arXiv:2307.00795}, 2024.

\bibitem[Chardon et~al.(2024)Chardon, Lerasle, and Mourtada]{chardon2024finite}
H.~Chardon, M.~Lerasle, and J.~Mourtada.
\newblock Finite-sample performance of the maximum likelihood estimator in
  logistic regression.
\newblock \emph{arXiv preprint arXiv:2411.02137}, 2024.

\bibitem[Chen(2018)]{chen2018gaussian}
X.~Chen.
\newblock Gaussian and bootstrap approximations for high-dimensional
  \({U}\)-statistics and their applications.
\newblock \emph{Annals of Statistics}, 46\penalty0 (2):\penalty0 642--678,
  2018.

\bibitem[Chernozhukov et~al.(2007)Chernozhukov, Hong, and
  Tamer]{chernozhukov2007estimation}
V.~Chernozhukov, H.~Hong, and E.~Tamer.
\newblock Estimation and confidence regions for parameter sets in econometric
  models.
\newblock \emph{Econometrica}, 75\penalty0 (5):\penalty0 1243--1284, 2007.

\bibitem[Chernozhukov et~al.(2017{\natexlab{a}})Chernozhukov, Chetverikov, and
  Kato]{chernozhukov2017clt}
V.~Chernozhukov, D.~Chetverikov, and K.~Kato.
\newblock Central limit theorems and bootstrap in high dimensions.
\newblock \emph{Annals of Probability}, 45\penalty0 (4):\penalty0 2309--2352,
  2017{\natexlab{a}}.

\bibitem[Chernozhukov et~al.(2017{\natexlab{b}})Chernozhukov, Chetverikov, and
  Kato]{chernozhukov2017detailed}
V.~Chernozhukov, D.~Chetverikov, and K.~Kato.
\newblock Detailed proof of {Nazarov's} inequality.
\newblock \emph{arXiv preprint arXiv:1711.10696}, 2017{\natexlab{b}}.

\bibitem[Chernozhukov et~al.(2018)Chernozhukov, Chetverikov, Demirer, Duflo,
  Hansen, Newey, and Robins]{chernozhukov2018double}
V.~Chernozhukov, D.~Chetverikov, M.~Demirer, E.~Duflo, C.~Hansen, W.~Newey, and
  J.~Robins.
\newblock Double/debiased machine learning for treatment and structural
  parameters.
\newblock \emph{The Econometrics Journal}, 21\penalty0 (1):\penalty0 C1--C68,
  2018.

\bibitem[Chernozhukov et~al.(2019)Chernozhukov, Chetverikov, and
  Kato]{Chernozhukov2019inference}
V.~Chernozhukov, D.~Chetverikov, and K.~Kato.
\newblock Inference on causal and structural parameters using many moment
  inequalities.
\newblock \emph{Review of Economic Studies}, 86\penalty0 (5):\penalty0
  1867--1900, 2019.

\bibitem[Chernozhukov et~al.(2022)Chernozhukov, Escanciano, Ichimura, Newey,
  and Robins]{chernozhukov2022locally}
V.~Chernozhukov, J.~C. Escanciano, H.~Ichimura, W.~K. Newey, and J.~M. Robins.
\newblock Locally robust semiparametric estimation.
\newblock \emph{Econometrica}, 90\penalty0 (1):\penalty0 327--354, 2022.

\bibitem[Chernozhukov et~al.(2023)Chernozhukov, Chetverikov, and
  Koike]{chernozhukov2023nearly}
V.~Chernozhukov, D.~Chetverikov, and Y.~Koike.
\newblock Nearly optimal central limit theorem and bootstrap approximations in
  high dimensions.
\newblock \emph{Annals of Applied Probability}, 33\penalty0 (3):\penalty0
  2374--2425, 2023.

\bibitem[Chernozhuokov et~al.(2022)Chernozhuokov, Chetverikov, Kato, and
  Koike]{chernozhuokov2022improved}
V.~Chernozhuokov, D.~Chetverikov, K.~Kato, and Y.~Koike.
\newblock Improved central limit theorem and bootstrap approximations in high
  dimensions.
\newblock \emph{Annals of Statistics}, 50\penalty0 (5):\penalty0 2562--2586,
  2022.

\bibitem[de~la Peña et~al.(2009)de~la Peña, Lai, and Shao]{pena2009self}
V.~H. de~la Peña, T.~L. Lai, and Q.-M. Shao.
\newblock \emph{Self-Normalized Processes: Limit Theory and Statistical
  Applications}, volume~35 of \emph{Probability and Its Applications}.
\newblock Springer, Berlin, 2009.

\bibitem[Delgado et~al.(2001)Delgado, Rodr{\'\i}guez-Poo, and
  Wolf]{delgado2001subsampling}
M.~A. Delgado, J.~M. Rodr{\'\i}guez-Poo, and M.~Wolf.
\newblock Subsampling inference in cube root asymptotics with an application to
  {Manski's} maximum score estimator.
\newblock \emph{Econometrica}, 69\penalty0 (3):\penalty0 675--712, 2001.

\bibitem[D{\"u}mbgen(1993)]{dumbgen1993nondifferentiable}
L.~D{\"u}mbgen.
\newblock On nondifferentiable functions and the bootstrap.
\newblock \emph{Probability Theory and Related Fields}, 95\penalty0
  (1):\penalty0 125--140, 1993.

\bibitem[Einmahl and Li(2008)]{einmahl2008characterization}
U.~Einmahl and D.~Li.
\newblock Characterization of {LIL} behavior in {B}anach space.
\newblock \emph{Transactions of the American Mathematical Society},
  360\penalty0 (12):\penalty0 6677--6693, 2008.

\bibitem[Fan et~al.(2018)Fan, Liu, Sun, and Zhang]{Fan2018ILAMM}
J.~Fan, H.~Liu, Q.~Sun, and T.~Zhang.
\newblock I-{LAMM} for sparse learning: simultaneous control of algorithmic
  complexity and statistical error.
\newblock \emph{Annals of Statistics}, 46\penalty0 (2):\penalty0 814--841,
  2018.

\bibitem[Gauss and Nagler(2024)]{gauss2024asymptotics}
J.~Gauss and T.~Nagler.
\newblock Asymptotics for estimating a diverging number of parameters-with and
  without sparsity.
\newblock \emph{arXiv preprint arXiv:2411.17395}, 2024.

\bibitem[Gin\'e et~al.(1997)Gin\'e, G\"otze, and Mason]{Gine1997}
E.~Gin\'e, F.~G\"otze, and D.~M. Mason.
\newblock When is the {S}tudent {$t$}-statistic asymptotically standard normal?
\newblock \emph{Annals of Probability}, 25\penalty0 (3):\penalty0 1514--1531,
  1997.

\bibitem[He and Shao(1996)]{he1996general}
X.~He and Q.-M. Shao.
\newblock A general {Bahadur} representation of \({M}\)-estimators and its
  application to linear regression with nonstochastic designs.
\newblock \emph{Annals of Statistics}, 24\penalty0 (6):\penalty0 2608--2630,
  1996.

\bibitem[He and Shao(2000)]{he2000parameters}
X.~He and Q.-M. Shao.
\newblock On parameters of increasing dimensions.
\newblock \emph{Journal of Multivariate Analysis}, 73\penalty0 (1):\penalty0
  120--135, 2000.

\bibitem[Huber(1967)]{huber1967behavior}
P.~J. Huber.
\newblock The behavior of maximum likelihood estimates under nonstandard
  conditions.
\newblock \emph{Proceedings of the Fifth Berkeley Symposium on Mathematical
  Statistics and Probability}, 1:\penalty0 221--233, 1967.

\bibitem[Huber(1973)]{huber1973robust}
P.~J. Huber.
\newblock Robust regression: asymptotics, conjectures and {Monte Carlo}.
\newblock \emph{Annals of Statistics}, pages 799--821, 1973.

\bibitem[Imbens and Manski(2004)]{imbens2004confidence}
G.~W. Imbens and C.~F. Manski.
\newblock Confidence intervals for partially identified parameters.
\newblock \emph{Econometrica}, 72\penalty0 (6):\penalty0 1845--1857, 2004.

\bibitem[Kim and Ramdas(2024)]{kim2024dimension}
I.~Kim and A.~Ramdas.
\newblock Dimension-agnostic inference using cross \({U}\)-statistics.
\newblock \emph{Bernoulli}, 30\penalty0 (1):\penalty0 1--30, 2024.

\bibitem[Kim and Ramdas(2025)]{kim2025locally}
I.~Kim and A.~Ramdas.
\newblock Locally minimax optimal and dimension-agnostic discrete argmin
  inference.
\newblock \emph{arXiv preprint arXiv:2503.21639}, 2025.

\bibitem[Kim and Pollard(1990)]{kim1990cube}
J.~Kim and D.~Pollard.
\newblock Cube root asymptotics.
\newblock \emph{Annals of Statistics}, 18\penalty0 (1):\penalty0 191--219,
  1990.

\bibitem[Knight(1998)]{knight1998limiting}
K.~Knight.
\newblock Limiting distributions for \({L}_1\) regression estimators under
  general conditions.
\newblock \emph{Annals of statistics}, pages 755--770, 1998.

\bibitem[Knight(2002)]{knight2002limiting}
K.~Knight.
\newblock What are the limiting distributions of quantile estimators?
\newblock In \emph{Statistical data analysis based on the {$L_1$}-norm and
  related methods ({N}euch\^atel, 2002)}, Stat. Ind. Technol., pages 47--65.
  Birkh\"auser, Basel, 2002.

\bibitem[Kock and Preinerstorfer(2024)]{kock2024remark}
A.~B. Kock and D.~Preinerstorfer.
\newblock A remark on moment-dependent phase transitions in high-dimensional
  {Gaussian} approximations.
\newblock \emph{Statistics \& Probability Letters}, 211:\penalty0 110149, 2024.

\bibitem[Kuchibhotla(2018)]{kuchibhotla2018deterministic}
A.~K. Kuchibhotla.
\newblock Deterministic inequalities for smooth \({M}\)-estimators.
\newblock \emph{arXiv preprint arXiv:1809.05172}, 2018.

\bibitem[Kuchibhotla and Chakrabortty(2022)]{kuchibhotla2022moving}
A.~K. Kuchibhotla and A.~Chakrabortty.
\newblock Moving beyond {sub-Gaussianity} in high-dimensional statistics:
  Applications in covariance estimation and linear regression.
\newblock \emph{Information and Inference: A Journal of the IMA}, 11\penalty0
  (4):\penalty0 1389--1456, 2022.

\bibitem[Kuchibhotla and Patra(2022)]{kuchibhotla2022least}
A.~K. Kuchibhotla and R.~K. Patra.
\newblock On least squares estimation under heteroscedastic and heavy-tailed
  errors.
\newblock \emph{Annals of Statistics}, 50\penalty0 (1):\penalty0 277--302,
  2022.

\bibitem[Kuchibhotla et~al.(2020)Kuchibhotla, Rinaldo, and
  Wasserman]{kuchibhotla2020berry}
A.~K. Kuchibhotla, A.~Rinaldo, and L.~Wasserman.
\newblock {Berry-Esseen} bounds for projection parameters and partial
  correlations with increasing dimension.
\newblock \emph{arXiv preprint arXiv:2007.09751}, 2020.

\bibitem[Kuchibhotla et~al.(2023)Kuchibhotla, Balakrishnan, and
  Wasserman]{kuchibhotla2023hulc}
A.~K. Kuchibhotla, S.~Balakrishnan, and L.~Wasserman.
\newblock The {HulC}: Confidence regions from convex hulls.
\newblock \emph{Journal of the Royal Statistical Society: Series B},
  86\penalty0 (3):\penalty0 586--622, 2023.

\bibitem[Ledoux and Talagrand(2011)]{Ledoux2011Isoperimetry}
M.~Ledoux and M.~Talagrand.
\newblock \emph{Probability in {B}anach spaces}.
\newblock Classics in Mathematics. Springer-Verlag, Berlin, 2011.
\newblock Isoperimetry and processes, Reprint of the 1991 edition.

\bibitem[Lin et~al.(2024)Lin, Su, Mou, Ding, and
  Wainwright]{lin2024worthwhilejackknifebreakingquadratic}
L.~Lin, F.~Su, W.~Mou, P.~Ding, and M.~J. Wainwright.
\newblock When is it worthwhile to jackknife? breaking the quadratic barrier
  for \({Z}\)-estimators.
\newblock \emph{arXiv preprint arXiv:2411.02909}, 2024.

\bibitem[Lopes(2022)]{lopes2022central}
M.~E. Lopes.
\newblock Central limit theorem and bootstrap approximation in high dimensions:
  Near \(1/\sqrt{n}\) rates via implicit smoothing.
\newblock \emph{Annals of Statistics}, 50\penalty0 (5):\penalty0 2492--2513,
  2022.

\bibitem[Mammen(1993)]{mammen1993bootstrap}
E.~Mammen.
\newblock Bootstrap and wild bootstrap for high dimensional linear models.
\newblock \emph{Annals of Statistics}, 21\penalty0 (1):\penalty0 255--285,
  1993.

\bibitem[Massart(1990)]{massart1990}
P.~Massart.
\newblock The tight constant in the {Dvoretzky-Kiefer-Wolfowitz} inequality.
\newblock \emph{Annals of Probability}, 18\penalty0 (3):\penalty0 1269 -- 1283,
  1990.

\bibitem[Meinshausen et~al.(2009)Meinshausen, Meier, and
  B{\"u}hlmann]{meinshausen2009p}
N.~Meinshausen, L.~Meier, and P.~B{\"u}hlmann.
\newblock $p$-values for high-dimensional regression.
\newblock \emph{Journal of the American Statistical Association}, 104\penalty0
  (488):\penalty0 1671--1681, 2009.

\bibitem[Mourtada(2022)]{mourtada2022exact}
J.~Mourtada.
\newblock Exact minimax risk for linear least squares, and the lower tail of
  sample covariance matrices.
\newblock \emph{Annals of Statistics}, 50\penalty0 (4):\penalty0 2157--2178,
  2022.

\bibitem[Nazarov(2003)]{nazarov2003maximal}
F.~Nazarov.
\newblock On the maximal perimeter of a convex set in {${\mathbb R}^n$} with
  respect to a {G}aussian measure.
\newblock In \emph{Geometric aspects of functional analysis}, volume 1807 of
  \emph{Lecture Notes in Math.}, pages 169--187. Springer, Berlin, 2003.

\bibitem[Nguyen et~al.(2020)Nguyen, Bagnall-Guerreiro, and
  Jones]{nguyen2020universal}
H.~D. Nguyen, J.~Bagnall-Guerreiro, and A.~T. Jones.
\newblock Universal inference with composite likelihoods.
\newblock \emph{arXiv preprint arXiv:2009.00848}, 2020.

\bibitem[Nickl and van~de Geer(2013)]{nickl2013confidence}
R.~Nickl and S.~van~de Geer.
\newblock Confidence sets in sparse regression.
\newblock \emph{Annals of Statistics}, 41\penalty0 (6):\penalty0 2852--2876,
  2013.

\bibitem[Oliveira(2016)]{oliveira2016lower}
R.~I. Oliveira.
\newblock The lower tail of random quadratic forms with applications to
  ordinary least squares.
\newblock \emph{Probability Theory and Related Fields}, 166:\penalty0
  1175--1194, 2016.

\bibitem[Ostrovskii and Bach(2021)]{ostrovskii2021finite}
D.~M. Ostrovskii and F.~Bach.
\newblock Finite-sample analysis of \({M}\)-estimators using self-concordance.
\newblock \emph{Electronic Journal of Statistics}, 15\penalty0 (1):\penalty0
  326--391, 2021.

\bibitem[Pakes and Pollard(1989)]{Pakes1989simulation}
A.~Pakes and D.~Pollard.
\newblock Simulation and the asymptotics of optimization estimators.
\newblock \emph{Econometrica}, 57\penalty0 (5):\penalty0 1027--1057, 1989.

\bibitem[Park et~al.(2023)Park, Balakrishnan, and Wasserman]{park2023robust}
B.~Park, S.~Balakrishnan, and L.~Wasserman.
\newblock Robust universal inference for misspecified models.
\newblock \emph{arXiv preprint arXiv:2307.04034}, 2023.

\bibitem[Petrov(1975)]{petrov1975sums}
V.~V. Petrov.
\newblock \emph{Sums of Independent Random Variables}, volume~82 of
  \emph{Ergebnisse der Mathematik und ihrer Grenzgebiete}.
\newblock Springer, Berlin, 1975.
\newblock Translated from the Russian by A. A. Brown.

\bibitem[Pflug(1995)]{Pflug1995}
G.~C. Pflug.
\newblock Asymptotic stochastic programs.
\newblock \emph{Mathematics of Operations Research}, 20\penalty0 (4):\penalty0
  769--789, 1995.

\bibitem[Politis et~al.(1999)Politis, Romano, and Wolf]{politis1999subsampling}
D.~N. Politis, J.~P. Romano, and M.~Wolf.
\newblock \emph{Subsampling}.
\newblock Springer Series in Statistics. Springer-Verlag, New York, 1999.

\bibitem[Portnoy(1984)]{portnoy1984asymptotic}
S.~Portnoy.
\newblock Asymptotic behavior of \({M}\)-estimators of $ p $ regression
  parameters when $ p^2/n $ is large. {I. Consistency}.
\newblock \emph{Annals of Statistics}, 12\penalty0 (4):\penalty0 1298--1309,
  1984.

\bibitem[Portnoy(1985)]{portnoy1985asymptotic}
S.~Portnoy.
\newblock Asymptotic behavior of \({M}\)-estimators of $ p $ regression
  parameters when $ p^2/n $ is large; {II. Normal} approximation.
\newblock \emph{Annals of Statistics}, 13\penalty0 (4):\penalty0 1403--1417,
  1985.

\bibitem[Portnoy(1986)]{portnoy1986asymptotic}
S.~Portnoy.
\newblock Asymptotic behavior of the empiric distribution of \({M}\)-estimated
  residuals from a regression model with many parameters.
\newblock \emph{Annals of Statistics}, 14\penalty0 (3):\penalty0 1152--1170,
  1986.

\bibitem[Portnoy(1988)]{portnoy1988asymptotic}
S.~Portnoy.
\newblock Asymptotic behavior of likelihood methods for exponential families
  when the number of parameters tends to infinity.
\newblock \emph{Annals of Statistics}, 16\penalty0 (1):\penalty0 356--366,
  1988.

\bibitem[P{\"o}tscher(2002)]{poetscher2002lower}
B.~M. P{\"o}tscher.
\newblock Lower risk bounds and properties of confidence sets for ill-posed
  estimation problems with applications to spectral density and persistence
  estimation, unit roots, and estimation of long memory parameters.
\newblock \emph{Econometrica}, 70\penalty0 (3):\penalty0 1035--1065, 2002.

\bibitem[Rinaldo et~al.(2019)Rinaldo, Wasserman, and
  G’Sell]{rinaldo2019bootstrapping}
A.~Rinaldo, L.~Wasserman, and M.~G’Sell.
\newblock Bootstrapping and sample splitting for high-dimensional,
  assumption-lean inference.
\newblock \emph{Annals of Statistics}, 47\penalty0 (6):\penalty0 3438--3469,
  2019.

\bibitem[Rio(2017)]{rio2017constants}
E.~Rio.
\newblock About the constants in the {F}uk-{N}agaev inequalities.
\newblock \emph{Electronic Communications in Probability}, 22:\penalty0 Paper
  No. 28, 12, 2017.

\bibitem[Robins et~al.(2006)Robins, van~der Vaart, and
  Ventura]{robins2006adaptive}
J.~Robins, A.~van~der Vaart, and V.~Ventura.
\newblock Adaptive nonparametric confidence sets.
\newblock \emph{Annals of Statistics}, 34\penalty0 (1):\penalty0 229--253,
  2006.

\bibitem[Rockafellar and Wets(1998)]{rockafellar1998variational}
R.~T. Rockafellar and R.~J.-B. Wets.
\newblock \emph{Variational Analysis}, volume 317 of \emph{Grundlehren der
  mathematischen Wissenschaften}.
\newblock Springer, Berlin, 1998.

\bibitem[Romano and Shaikh(2010)]{romano2010inference}
J.~P. Romano and A.~M. Shaikh.
\newblock Inference for the identified set in partially identified econometric
  models.
\newblock \emph{Econometrica}, 78\penalty0 (1):\penalty0 169--211, 2010.

\bibitem[Romano and Shaikh(2014)]{romano2014practical}
J.~P. Romano and A.~M. Shaikh.
\newblock Practical confidence regions for parameter sets in econometric
  models.
\newblock \emph{Econometrica}, 82\penalty0 (3):\penalty0 1697--1714, 2014.

\bibitem[Romano and Wolf(2000)]{romano2000finite}
J.~P. Romano and M.~Wolf.
\newblock Finite sample nonparametric inference and large sample efficiency.
\newblock \emph{Annals of Statistics}, 28\penalty0 (3):\penalty0 820--847,
  2000.

\bibitem[Shi and Drton(2024)]{shi2024universal}
H.~Shi and M.~Drton.
\newblock On universal inference in {Gaussian} mixture models.
\newblock \emph{arXiv preprint arXiv:2407.19361}, 2024.

\bibitem[Sun et~al.(2020)Sun, Zhou, and Fan]{Sun2020adaptivehuber}
Q.~Sun, W.-X. Zhou, and J.~Fan.
\newblock Adaptive {H}uber regression.
\newblock \emph{Journal of the American Statistical Association}, 115\penalty0
  (529):\penalty0 254--265, 2020.

\bibitem[Takatsu(2025)]{takatsu2025precise}
K.~Takatsu.
\newblock On the precise asymptotics of universal inference.
\newblock \emph{arXiv preprint arXiv:2503.14717}, 2025.

\bibitem[Takatsu and Kuchibhotla(2025)]{takatsu2025bridging}
K.~Takatsu and A.~K. Kuchibhotla.
\newblock Bridging root-$n$ and non-standard asymptotics: Dimension-agnostic
  adaptive inference in \({M}\)-estimation.
\newblock \emph{arXiv preprint arXiv:2501.07772}, 2025.

\bibitem[Tamer(2003)]{tamer2003incomplete}
E.~Tamer.
\newblock Incomplete simultaneous discrete response model with multiple
  equilibria.
\newblock \emph{Review of Economic Studies}, 70\penalty0 (1):\penalty0
  147--165, 2003.

\bibitem[Tikhomirov(2018)]{tikhomirov2018sample}
K.~Tikhomirov.
\newblock Sample covariance matrices of heavy-tailed distributions.
\newblock \emph{International Mathematics Research Notices}, 2018\penalty0
  (20):\penalty0 6254--6289, 2018.

\bibitem[Tropp et~al.(2015)]{tropp2015introduction}
J.~A. Tropp et~al.
\newblock An introduction to matrix concentration inequalities.
\newblock \emph{Foundations and Trends in Machine Learning}, 8\penalty0
  (1-2):\penalty0 1--230, 2015.

\bibitem[van~der Vaart and Wellner(1996)]{van1996weak}
A.~W. van~der Vaart and J.~A. Wellner.
\newblock \emph{Weak Convergence and Empirical Processes: With Applications to
  Statistics}.
\newblock Springer Series in Statistics. Springer, 1996.

\bibitem[Vershynin(2018)]{vershynin2018high}
R.~Vershynin.
\newblock \emph{High-Dimensional Probability: An Introduction with Applications
  in Data Science}, volume~47 of \emph{Cambridge Series in Statistical and
  Probabilistic Mathematics}.
\newblock Cambridge University Press, 2018.

\bibitem[Vladimirova et~al.(2020)Vladimirova, Girard, Nguyen, and
  Arbel]{vladimirova2020subweibull}
M.~Vladimirova, S.~Girard, H.~D. Nguyen, and J.~Arbel.
\newblock {Sub-Weibull} distributions: generalizing {sub-Gaussian} and
  {sub-Exponential} properties to heavier-tailed distributions.
\newblock \emph{Stat}, 9\penalty0 (1):\penalty0 e318, 2020.

\bibitem[von Bahr and Esseen(1965)]{vonbahr1965inequalities}
B.~von Bahr and C.-G. Esseen.
\newblock Inequalities for the $r$th absolute moment of a sum of random
  variables, $1 \leq r \leq 2$.
\newblock \emph{Annals of Mathematical Statistics}, 36\penalty0 (1):\penalty0
  299--303, 1965.

\bibitem[Wasserman and Roeder(2009)]{wasserman2009high}
L.~Wasserman and K.~Roeder.
\newblock High-dimensional variable selection.
\newblock \emph{Annals of Statistics}, 37\penalty0 (5A):\penalty0 2178--2201,
  2009.

\bibitem[Wasserman et~al.(2020)Wasserman, Ramdas, and
  Balakrishnan]{wasserman2020universal}
L.~Wasserman, A.~Ramdas, and S.~Balakrishnan.
\newblock Universal inference.
\newblock \emph{Proceedings of the National Academy of Sciences}, 117\penalty0
  (29):\penalty0 16880--16890, 2020.

\bibitem[White(1980)]{white1980heteroskedasticity}
H.~White.
\newblock A heteroskedasticity-consistent covariance matrix estimator and a
  direct test for heteroskedasticity.
\newblock \emph{Econometrica}, 48\penalty0 (4):\penalty0 817--838, 1980.

\bibitem[Yohai and Maronna(1979)]{yohai1979asymptotic}
V.~J. Yohai and R.~A. Maronna.
\newblock Asymptotic behavior of {$M$}-estimators for the linear model.
\newblock \emph{Ann. Statist.}, 7\penalty0 (2):\penalty0 258--268, 1979.

\end{thebibliography}

\appendix
\numberwithin{equation}{section}
\numberwithin{lemma}{section}
\numberwithin{proposition}{section}
\numberwithin{remark}{section}
\numberwithin{theorem}{section}
\numberwithin{figure}{section}
\numberwithin{table}{section}
\renewcommand{\theequation}{\thesection.\arabic{equation}}
\renewcommand{\theassumption}{\thesection.\arabic{assumption}}
\renewcommand{\thelemma}{\thesection.\arabic{lemma}}
\renewcommand{\theproposition}{\thesection.\arabic{proposition}}
\renewcommand{\thetheorem}{\thesection.\arabic{theorem}}
\renewcommand{\thecorollary}{\thesection.\arabic{corollary}}
\renewcommand{\theremark}
{\thesection.\arabic{remark}}
\renewcommand{\thefigure}{\thesection.\arabic{figure}}
\renewcommand{\thetable}{\thesection.\arabic{table}}

\section{Width Analysis on CIs for Linear Contrast of General \texorpdfstring{$Z$}{Z}-functionals}
Classical asymptotic approximation for $Z$-estimators (e.g., Theorem~3.3.1 of \citep{van1996weak}) suggests $\sqrt{n}(\hat\theta_1-\theta)\approx \Nc(0_d,J^{-1}VJ^{-1})$ where $\hat\theta_1$ is the $Z$-estimator from $\Dc_1$ and $J = \Eb[\nabla_\theta\psi(Z,\theta_0)]$ is the Jacobian matrix at $\theta_0$. This approximation leads to a canonical Wald interval for $\xi^\top\theta_0$: $[\xi^\top\hat\theta_{\Dc_1}\pm z_{\alpha/2}\sqrt{\xi^\top J^{-1}VJ^{-1}\xi/n}]$. The following theorem establishes the second order width analysis for the confidence sets described in Section~3.

\begin{theorem}\label{thm:A1}
    Grant assumptions in Theorem~\ref{thm:5} and assume $\sigma_{\rm min}^{-1}(J)\vee \sigma_{\rm max}(J)=O(1)$. and let $\set{r_{n_i}:n_i\geq 1}$  $(i=1,2)$, $\set{s_n\gtrsim n^{-1/2}}$, $\set{A_n},\set{B_n},\set{C_n}$ be such that, as $n,n\to\infty$, $\Pb(\norm{\hat\theta_i-\theta_0}_2\geq r_{n_i})\to0$ ($i=1,2$), $\Pb(\sup\set{\norm{\theta-\theta_0}_2: \theta\in\widehat{\mathrm{CI}}^{\rm SpN}_{\alpha/n}}\geq s_n)\to0$, $\sup_{\theta:\norm{\theta-\theta_0}\leq s_n}\norm{\nabla\bar\psi_1(\theta)-J}_{\rm op}=O_\Pb(A_n),$ $\norm{J\widetilde Q_2^\top-I_d}_{\rm op} = O_\Pb(B_n)$, and $\sup_{\theta:\norm{\theta-\theta_0}\leq s_n}\norm{V(\theta)-\hat V_1(\theta)}_{\rm op}=O_\Pb(C_n)$. Let $D_n^{(1)}=\sqrt{n} s_n\left(A_n+B_n\right) +C_n+\phi(s_n)$, and $D_n^{(2)}=\sqrt{n} s_n\left(A_n+B_n\right) +\phi(r_{n})$. Then,
    \begin{equation*}        
        \sqrt{n}~{\rm d}_{\rm H}(\widehat{\rm CI}_\alpha^{\rm SeN}\to\widehat{\rm CI}_\alpha^{\rm Wald})= O_\Pb(D_n^{(1)}),\mbox{ and }\sqrt{n}~{\rm d}_{\rm H}(\widehat{\rm CI}_\alpha^{\rm SpN}\to\widehat{\rm CI}_\alpha^{\rm Wald})= O_\Pb(D_n^{(2)}).
    \end{equation*}
\end{theorem}
\section{Proofs of Theorems}
\subsection{Proof of Theorem~\ref{thm:1}} Theorem~\ref{thm:1} follows from Theorems \ref{thm:B1}, \ref{thm:B2}, and Lemma~\ref{lem:B1}. Proofs of Theorem~\ref{thm:B1}, \ref{thm:B2}, and Lemma~\ref{lem:B1} are deferred to the end of this section.
\begin{theorem}\label{thm:B1} For an independent $g\sim \Nc(0_d,I_d)$, let $S(t)=\Pb(\norm{[\Gamma(\theta_0,\widetilde Q_2)]^{1/2}g}_\infty>t|\Dc_2).$ Define $\Ec^{\rm SeN}:=\sup_{t\geq0}[\Pb(T^{\rm SeN}(\theta_0)>t|\Dc_2 )-S(t)]_+$. There exists a universal constant $C>0$,
\begin{equation}\label{eq:thm:B1:1}
    \Ec^{\rm SeN}\leq CL_q\bigg[\frac{\log^{5/4}(ed)}{n^{1/4}}+\frac{\log^{3/2}(ed) d^{1/q}}{n^{1/2-1/q}}\bigg]:=\varepsilon_{n,q}^{(1)},
\end{equation} with probability $1$, for $q>2$. Moreover, for any $\lambda\geq 0$, define $\Dc_2$-measurable event $\Ec_\lambda:=\set{\lambda_{\rm min}(\Gamma(\theta_0,\widetilde Q_2))\geq\lambda}$. For $q\geq 4$, one has, with probability 1,
\begin{equation}\label{eq:thm:B1:2}
    \Ec^{\rm SeN}\leq C\bigg\{\log^2(en)\bigg(\frac{L_q^2\log^{3}(ed) d^{2/q}}{n^{1-2/q}}\bigg)^{\frac{q}{2(q-2)}}\mathbbm{1}_{\Ec_\lambda}+\varepsilon_{n,q}^{(1)}\mathbbm{1}_{\Ec_\lambda^\complement}\bigg\}:=\varepsilon_{n,q}^{(2)}(\lambda),
\end{equation}
\end{theorem}

For a matrix $M\in\Real^{d\times d}$, let $\mnorm{M}_{\infty}=\max_{j,k\in[d]}|\ev_j^\top M\ev_k|$.

\begin{theorem}\label{thm:B2}For a (stochastic) nonnegative-definite $\Gamma^\dagger\in\Real^{d\times d}$ and $\alpha\in(0,1)$, define $K_{\alpha,B}(\Gamma^\dagger):=\inf\{t\geq 0:\frac{1}{B}\sum_{i=1}^B\mathbbm{1}(\norm{[\Gamma^\dagger]^{1/2}W_i^{\rm b}}_\infty>t)\leq\alpha\}$ where $W_i^{\rm b}\sim \Nc(0_d,I_d)$ are independent of everything. Let $\Upsilon(\Gamma^\dagger) = \mnorm{\Gamma^\dagger-\Gamma(\theta_0,\widetilde Q_2)}_\infty$ and $\Ec^{\rm SeN}(\Gamma^\dagger)=\sup_{\alpha\in(0,1)}[\Pb(T^{\rm SeN}(\theta_0)>K_{\alpha,B}(\Gamma^\dagger))-\alpha]_+$. Then, there exists a universal constant $c>0$ such that, for $q>2$,
\begin{equation}\label{eq:thm:B2:1}
    \Ec^{\rm SeN}(\Gamma^\dagger)\leq \varepsilon_{n,q}^{(1)}+\sqrt{\frac{\log(en)}{2B}}+\inf_{t>0}\bigg[c\log(ed)t^{1/2}+\Pb(\Upsilon(\Gamma^\dagger)>t)\bigg].
\end{equation}Moreover, let $p_\lambda = \Pb(\Ec_\lambda^\complement)$. For any $\lambda>0$ and $q\geq 4$,
\begin{equation}\label{eq:thm:B2:2}
    \Ec^{\rm SeN}(\Gamma^\dagger)\leq (p_\lambda \varepsilon_{n,q}^{(1)} + \varepsilon_{n,q}^{(2)}(\lambda) )+ \sqrt{\frac{\log(en)}{2B}}+\inf_{t>0}\bigg[cp_\lambda\log(ed)t^{1/2}+c\lambda^{-1}\log(ed)t\log(e/t)+\Pb(\Upsilon(\Gamma^\dagger)>t)\bigg].
\end{equation}
\end{theorem}

\begin{lemma}\label{lem:B1}Let $$w_{n,q}(\delta,\varepsilon):=CL_q^2\gamma_q^2(\delta)\frac{d^{2/q}}{n^{1-2/q}}(\log^{1-2/q}(ed)+\varepsilon^{-2/q})+4\phi(\delta),$$ for some constant $C=C(q)\geq 1$. For any $\varepsilon\in(0,1)$ and $\delta\in(0,\infty)$, it holds w.p.1 that $\Pb(\mnorm{\hat\Delta}_\infty\geq 2\wedge w_{n,q}(\delta,\varepsilon)|\Dc_2)\leq\varepsilon+\mathbbm{1}(\norm{\widetilde\theta_2-\theta_0}_2>\delta)$ where $\hat\Delta = \widehat\Gamma_1(\widetilde\theta_2;\widetilde Q_2)-\Gamma(\theta_0,\widetilde Q_2)$.
\end{lemma}

\noindent We denote $w_1(\delta) = CL_q^2\gamma_q^2(\delta){d^{2/q}}/{n^{1-2/q}}$ and $w_2(\delta)= w_1(\delta)\log^{1-2/q}(ed)+4\phi(\delta).$ It follows from the construction of $w_{n,q}$ in Lemma~\ref{lem:B1} that $w_{n,q}(\delta,\varepsilon) =  w_1(\delta)\varepsilon^{-2/q}+ w_2(\delta)$ for $\delta>0$. 
\paragraph{Proof of \eqref{eq:thm1:1}} Take $\Gamma^\dagger=\widehat\Gamma_1(\widetilde\theta_2;\widetilde Q_2)$ in Theorem~\ref{thm:B2}, then it suffices to control $\inf_{t> 0}\Eb[\Rk_1(t)]$ where $\Rk_1(t)=c\log(ed)\sqrt{t}+\Pb(\mnorm{\hat\Delta}_\infty >t|\Dc_2)$. Lemma~\ref{lem:B1} with $t=w_{n,q}(\delta,\varepsilon)$ implies that
\begin{align*}
    &\inf_{t\geq 0}\Eb[\Rk_1(t)]\leq \inf_{\epsilon\in(0,1),\delta>0}\big\{c\log(ed)w^{1/2}_{n,q}(\delta,\varepsilon)+\varepsilon+\Pb(\norm{\widetilde\theta_2-\theta_0}_2>\delta)\big\}\\
    &\quad\leq\inf_{\delta>0}\big\{\Pb(\norm{\widetilde\theta_2-\theta_0}_2>\delta)+c\log(ed)w_2^{1/2}(\delta)+\inf_{\varepsilon\in(0,1)}[\varepsilon+c\log(ed)w_1^{1/2}(\delta)\varepsilon^{-1/q}]\big\}\\
    &\qquad\leq \inf_{\delta>0}\big\{\Pb(\norm{\widetilde\theta_2-\theta_0}_2>\delta)+c\log(ed)w_2^{1/2}(\delta)+c\log^{q/(q+1)}(ed)w_1^{q/(2q+2)}(\delta)\big\}.
\end{align*} Last inequality follows from taking $\varepsilon = [\log(ed)w_1^{1/2}(\delta)]^{\frac{q}{q+1}}$ , and simplifying this yields \eqref{eq:thm1:1}.
\paragraph{Proof of \eqref{eq:thm1:2}} Again, it suffices to control $\inf_{t\geq 0}\Eb[\Rk_2(t)]$ where $\Rk_2(t)=c\lambda^{-1}\log(ed)t\log(e/t)+cp_\lambda\log(ed)t^{1/2}+\Pb(\mnorm{\hat\Delta}_\infty >t\big|\Dc_2)$. Take $t=w_{n,q}(\delta,\varepsilon)$ in Lemma~\ref{lem:B1}, and we observe that $\log(e/w_{n,q}(\delta,\varepsilon))\leq \log(e n)$, for all $\delta,\varepsilon$. This implies that for any $\delta,\lambda>0$,
\begin{align}\label{eq:pf_of27.1}
    &\inf_{t\geq 0}\Eb[\Rk_2(t)]\leq ~\Pb(\norm{\widetilde\theta_2-\theta_0}_2>\delta)+c[p_\lambda\log(ed)w_2^{1/2}(\delta)+\lambda^{-1}\log(ed)\log(en)w_2(\delta)]\nonumber\\
    &\quad+c\inf_{\epsilon>0}[\epsilon + cp_\lambda\log(ed)w_1^{1/2}(\delta)\varepsilon^{-1/q} + c\lambda^{-1}\log(ed)\log(en)w_1(\delta)\varepsilon^{-2/q}].
\end{align} Since $p_{\lambda_\circ}=0$, taking $\varepsilon = [\log(ed)\log(en)w_1(\delta)/\lambda_\circ]^{q/(q+2)}$ and simplifying yields \eqref{eq:thm1:2}.

\paragraph{proof of Theorem~\ref{thm:B1}}
This result is a consequence of Theorem~\ref{thm:B3}. Applying the first part of Theorem~\ref{thm:B3} conditional on $\Dc_2$ yields the bound in \eqref{eq:thm:B1:1}. Furthermore, combining an application of the first part of Theorem~\ref{thm:B3} on the complement of the event $\Ec_\lambda$ with the second part of the same theorem establishes the remaining claim.

The following results do not involve sample splitting. 
\begin{theorem}\label{thm:B3}
    For a non-stochasitc $Q\in\Real^{d\times d}$, let $\check G_\theta\sim\Nc(0_{d},\Gamma(\theta;Q))$. For any $\theta\in\Theta$ and $q>2$, there exists a universal constant $C>0$ such that
    \begin{equation}\label{eq:thmB3:1}
    \sup_{t\geq 0}[\Pb(T^{\rm SeN}(\theta;Q))-\Pb(\norm{\check G_\theta}_\infty\geq t)]\leq C \left[\frac{L_4(\theta)\log^{5/4}(ed)}{n^{1/4}}+\frac{L_q(\theta)\log^{3/2}(ed)d^{1/q}}{n^{1/2-1/q}}\right],
\end{equation} Moreover, for $\theta\in\Theta$ such that there further exists $\lambda(\theta)>0$ where $\lambda_{\rm min}({\rm Corr}(\Gamma(\theta;Q)))\geq \lambda(\theta)$ for all $n$, then
\begin{equation}\label{eq:thmB3:2}
    \sup_{t\geq 0}[\Pb(T^{\rm SeN}(\theta;Q))-\Pb(\norm{\check G_\theta}_\infty\geq t)]\leq C\left[\frac{L_q^{q/(q-2)}(\theta)\log^2(ed)\log(en)\log(end)d^{1/(q-2)}}{n^{1/2}\lambda^{1/(q-2)}(\theta)}\right],
\end{equation} for any $q\geq 4$, where $C>0$ is a universal constant.
\end{theorem}

For a given $\theta\in\Theta$ and non-stochastic $Q\in\Real^{d\times d}$, define a $d$-dimensional random vectors $V_{i,\theta}=(V_{i1,\theta},\ldots,V_{id,\theta})^\top$ for $i\in\Ic_1$, where
\begin{equation*}
    V_{ij,\theta} = \frac{\ev_j^\top Q(\psi(Z_i,\theta)-\mu(\theta))}{(\ev_j^\top Q V(\theta)Q^\top\ev_j)^{1/2}}\quad\mbox{for}\quad j\in[d].
\end{equation*} From the construction it follows that $V_{i,\theta}$ is mean zero and ${\rm Var}(V_{i,\theta})=\Gamma(\theta;Q)$. Recall $\check G_\theta\sim\Nc(0_{d},\Gamma(\theta;Q))$, and define
\begin{eqnarray*}
    \varrho_n^\theta:=\sup_{t\geq0}\Abs{\Pb\left(\Norm{n_1^{-1/2}\textstyle\sum_{i\in\Ic_1}V_{i,\theta}}_\infty\geq t\right)-\Pb\left(\Norm{\check G_\theta}_\infty\geq t\right)}.
\end{eqnarray*} 
The proof of Theorem~\ref{thm:B3} is based on Lemma~\ref{lem:B2} and Proposition~\ref{prop:B1}, whose proofs are deferred to the end of this subsection. Lemma~\ref{lem:B2} establishes intermediate bounds on the one-sided Kolmogorov–Smirnov distance in terms of moments and the quantity $\varrho_n^\theta$, leveraging an exponential inequality for self-normalized sums stated in Lemma~\ref{lem:C3}. Proposition~\ref{prop:B1} applies recent high-dimensional central limit theorems from \citet{chernozhuokov2022improved} and \citet{chernozhukov2023nearly} to provide an upper bound for $\varrho_n^\theta$.

\begin{lemma}\label{lem:B2} There exists a universal constant $C>0$ such that
\begin{align}\label{eq:lem:B2}
    &\qquad\sup_{t\geq0}\left[\Pb\left(T^{\rm SeN}(\theta;Q)>t\right)-\Pb\left(\Norm{\check G_\theta}_\infty>t\right)\right]_+\nonumber\\
    &\leq \varrho_n^\theta+\frac{d}{n^{q/2-1}}+Cq^{3/2-2/q}\frac{L_q(\theta)^2\log^{3/2-2/q}(en)\log^{1/2}(ed)}{n^{1-2/q}},
\end{align} for any $\theta\in\Theta$ and non-stochastic $Q\in\Real^{d\times d}$.
    
\end{lemma}

\begin{proposition}\label{prop:B1} For any $q>2$, there exists a universal constant $C>0$ such that
    \begin{equation}\label{eq:prop:B1:1}
    \varrho_n^\theta\leq C \left[\frac{L_4(\theta)\log^{5/4}(ed)}{n^{1/4}}+\frac{L_q(\theta)\log^{3/2}(ed)d^{1/q}}{n^{1/2-1/q}}\right],
\end{equation} Moreover, for $\theta\in\Theta$ such that there further exists $\lambda(\theta)>0$ where $\lambda_{\rm min}({\rm Corr}(\Gamma(\theta;Q)))\geq \lambda(\theta)$ for all $n$, then
\begin{equation}\label{eq:prop:B1:2}
    \varrho_n^\theta\leq C\left[\frac{L_q^{q/(q-2)}(\theta)\log^2(ed)\log(en)\log(end)d^{1/(q-2)}}{n^{1/2}\lambda^{1/(q-2)}(\theta)}\right],
\end{equation} for any $q\geq 4$, where $C>0$ is a universal constant.
\end{proposition}

We may assume that $d^{2/q}\leq n^{1-2/q}$, otherwise the inequalities in Theorem~\ref{thm:B3} hold by taking sufficiently large constant. The rightmost term in \eqref{eq:lem:B2} is bounded as
\begin{equation*}
    \mbox{Rightmost term in \eqref{eq:lem:B2}}\leq C\left[\frac{L^2_q(\theta)\log^{3/2-2/q}(en)\log^{1/2}(ed)}{n^{1-2/q}}\right].
\end{equation*} Since for $\epsilon>0$ there exists $c_\epsilon>0$ such that $\log(en)\leq c_\epsilon n^\epsilon$ for all $n\geq 1$, we have, for $q>2$
\begin{align*}
    &\frac{L^2_q(\theta)\log^{3/2-2/q}(en)\log^{1/2}(ed)}{n^{1-2/q}}\\
    &\leq C_q \min\Set{\frac{L_q(\theta)\log^{3/2}(ed)d^{1/q}}{n^{1/2-1/q}}, \frac{L_q^{q/(q-2)}(\theta)\log^2(ed)\log(en)\log(end)d^{1/(q-2)}}{n^{1/2}}}.
\end{align*} Therefore, combining \eqref{eq:lem:B2} and \eqref{eq:prop:B1:1} proves \eqref{eq:thmB3:1}, and combining \eqref{eq:lem:B2} and \eqref{eq:prop:B1:2} proves \eqref{eq:thmB3:2}.

\begin{proof}[proof of Lemma~\ref{lem:B2}] We begin by noting that
\begin{equation*}
    T^{\rm SeN}(\theta)=\max_{j\in[d]}\frac{n_1\abs{\ev_j^\top Q(\bar\psi_1(\theta)-\mu(\theta))}}{\sqrt{\sum_{i\in\Ic_1}(\ev_j^\top Q(\psi(Z_i,\theta)-\mu(\theta))^2-n_1\abs{\ev_j^\top Q(\bar\psi_1(\theta)-\mu(\theta))}^2}}.
\end{equation*} Therefore, we have
\begin{equation}\label{eq:lemA1:1}
    \Set{T^{\rm SeN}(\theta)>t}=\Set{S^{\rm SeN}(\theta)>\frac{t}{\sqrt{1+t^2/n_1}}}~\mbox{where}~S^{\rm SeN}(\theta)=\max_{j\in[d]}\frac{n_1\abs{\ev_j^\top Q(\bar\psi_1(\theta)-\mu(\theta))}}{\sqrt{\sum_{i\in\Ic_1}(\ev_j^\top Q(\psi(Z_i,\theta)-\mu(\theta))^2}}.
\end{equation} Define
\begin{equation}\label{eq:def.c_nq}
    c_{n,q} = L_q(\theta)^2\left(\frac{q\log(en_1)}{2n_1}\right)^{1-2/q},
\end{equation}and assume that $c_n\leq 1/3$, otherwise the inequality in Lemma~\ref{lem:B2} holds by taking a sufficiently large universal constant. Further denote
\begin{equation*}
    R_{n,\Ac} = \min_{j\in[d]}\sqrt{\frac{\sum_{i\in\Ic_1}(\ev_j^\top Q(\psi(Z_i,\theta)-\mu(\theta))^2}{\sum_{i\in\Ic_1}\Eb[(\ev_j^\top Q(\psi(Z_i,\theta)-\mu(\theta))^2]}}.
\end{equation*}
From \eqref{eq:lemA1:1}, it follows from the simple implication that
\begin{align}\label{eq:pf_lemA.1:0}
    &\Pb\left(T^{\rm SeN}(\theta)>t\right)=\Pb\left(S^{\rm SeN}(\theta)>\frac{t}{\sqrt{1+t^2/n_1}}\right)\nonumber\\
    &~=\Pb\left(S^{\rm SeN}(\theta)>\frac{t}{\sqrt{1+t^2/n_1}},~R_{n,\Ac}>\sqrt{1-c_{n,q}}\right)+\Pb\left(S^{\rm SeN}(\theta)>\frac{t}{\sqrt{1+t^2/n_1}},~R_{n,\Ac}\leq\sqrt{1-c_{n,q}}\right)\nonumber\\
    &~\leq \Pb\left(R_{n,\Ac}S^{\rm SeN}(\theta)>\frac{t\sqrt{1-c_{n,q}}}{\sqrt{1+t^2/n_1}}\right)+\Pb\left(R_{n,\Ac}\leq\sqrt{1-c_{n,q}}\right)\nonumber\\
    &~\leq\Pb\left(\max_{j\in[d]}\frac{n_1\abs{\ev_j^\top Q(\bar\psi_1(\theta)-\mu(\theta))}}{\sqrt{\sum_{i\in\Ic_1}\Eb[(\ev_j^\top Q(\psi(Z_i,\theta)-\mu(\theta))^2]}}>\frac{t\sqrt{1-c_{n,q}}}{\sqrt{1+t^2/n_1}}\right)+\Pb\left(R_{n,\Ac}\leq\sqrt{1-c_{n,q}}\right)\nonumber\\
    &~=\Pb\left(\Norm{\frac{1}{\sqrt{n_1}}\sum_{i\in\Ic_1}V_{i,\theta}}_\infty> \frac{t\sqrt{1-c_{n,q}}}{\sqrt{1+t^2/n_1}}\right)+\Pb\left(R_{n,\Ac}\leq\sqrt{1-c_{n,q}}\right).
\end{align} The quantity $R_{n,\Ac}$ can be controlled using Lemma~\ref{lem:C4}. In particular, we check that the $q$:th Lyapunov's ratio of $a^\top(\psi(Z_i,\theta)-\mu(\theta))$ can be bounded as
\begin{eqnarray*}
    &&\max_{j\in[d]}\frac{( \Eb\abs{\ev_j^\top Q(\psi(Z_i,\theta)-\mu(\theta)}^q])^{1/q}}{(\Eb(\ev_j^\top Q(\psi(Z_i,\theta)-\mu(\theta))^2])^{1/2}}\leq L_q(\theta).
\end{eqnarray*} Now Lemma~\ref{lem:C4} applies and yields
\begin{equation}\label{eq:pf_lemA.1:0.1}
    \Pb\left(R_{n,\Ac}\leq \sqrt{1-c_{n,q}}\right)\leq \frac{d}{n^{q/2-1}}.
\end{equation} In order to control the leading term in \eqref{eq:pf_lemA.1:0}, we recall $\check G_\theta\sim\Nc(0_{d},\Gamma(\theta;Q))$ and set
\begin{equation*}
    \Phi_{\rm AC}^\theta:=\sup_{t\geq 0, \varepsilon>0}\frac{1}{\varepsilon}\left[\Pb\left(\norm{\check G_\theta}_\infty>t\right)-\Pb\left(\norm{\check G_\theta}_\infty>t+\varepsilon\right)\right],
\end{equation*} as the anti-concentration constant. The anti-concentration constant for maximum Gaussian is controlled using Nazarov's inequality \citep{nazarov2003maximal}. Lemma~A.1 of \cite{chernozhukov2017clt} implies that there exists a universal constant $C>0$ such that
\begin{equation*}
    \Phi_{\rm AC}^\theta\leq C\sqrt{\log(ed)}.
\end{equation*} Meanwhile, since $(\Gamma(\theta;Q))_{jj}= 1$ for all $j\in[d]$, Gaussian maximal concentration inequality in Lemma~\ref{lem:C1} implies that for any $q>2$,
\begin{equation}\label{eq:gauss_conc}
    \Pb\left(\norm{\check G_\theta}_\infty\geq \sqrt{(q-2)\log(en_1)}\right)\leq\frac{d}{n^{q/2-1}}.
\end{equation} 

Now, we are ready to control the leading term of  \eqref{eq:pf_lemA.1:0}. First, for $t\geq0$ such that $t/\sqrt{1+t^2/n_1}\geq \sqrt{(q-2)\log(en_1)/(1-c_{n,q})}$, we have
\begin{align}\label{eq:pf_lemA.1:1}
&\Pb\left(\Norm{\frac{1}{\sqrt{n_1}}\sum_{i\in\Ic_1}V_{i,\theta}}_\infty> \frac{t\sqrt{1-c_{n,q}}}{\sqrt{1+t^2/n_1}}\right)\leq\Pb\left(\norm{\check G_\theta}_\infty>\frac{t\sqrt{1-c_{n,q}}}{\sqrt{1+t^2/n_1}}\right) +\varrho_n^\theta\nonumber\\
&\qquad\leq\Pb\left(\norm{\check G_\theta}_\infty\geq \sqrt{(q-2)\log(en_1)}\right)+\varrho_n^\theta\leq \frac{d}{n^{q/2-1}}+\varrho_n^\theta,
\end{align} where the first inequality follows from the definition of $\varrho_n^\theta$ and the last from \eqref{eq:gauss_conc}. Meanwhile, for $t\geq0$ such that $t/\sqrt{1+t^2/n_1}<\sqrt{(q-2)\log(en_1)/(1-c_{n,q})}$, we observe that
\begin{align}\label{eq:pf_lemA.1:1.5}
    &\frac{t}{\sqrt{1+t^2/n_1}}<\sqrt{\frac{(q-2)\log(en_1)}{1-c_{n,q}}}~\Longrightarrow~t^2\left(1-\frac{(q-2)\log(en_1)}{n(1-c_{n,q})}\right)<\frac{(q-2)\log(en_1)}{1-c_{n,q}}\nonumber\\
    &\qquad~\Longrightarrow~t^2<\left(1-\frac{(q-2)\log(en_1)}{n(1-c_{n,q})}\right)_+^{-1}\frac{(q-2)\log(en_1)}{1-c_{n,q}}.
\end{align} This implies that
\begin{align}\label{pf_lemA.1:1.6}
    &t-\frac{t\sqrt{1-c_{n,q}}}{\sqrt{1+t^2/n_1}}=\frac{t}{\sqrt{1+t^2/n_1}}\left(\sqrt{1+t^2/n_1}-\sqrt{1-c_{n,q}}\right)\nonumber\\
    &\quad=\frac{t}{\sqrt{1+t^2/n_1}}\frac{t^2/n_1+c_{n,q}}{\sqrt{1+t^2/n_1}+\sqrt{1-c_{n,q}}}\leq \frac{t}{\sqrt{1+t^2/n_1}}\left(t^2/n_1+c_{n,q}\right)\nonumber\\
    &\quad\leq\sqrt{\frac{(q-2)\log(en_1)}{1-c_{n,q}}}\left(c_{n,q}+\left(1-\frac{(q-2)\log(en_1)}{n_1(1-c_{n,q})}\right)_+^{-1}\frac{(q-2)\log(en_1)}{n_1(1-c_{n,q})}\right),
\end{align} where the last inequality follows from \eqref{eq:pf_lemA.1:1.5} and the others from simple algebraic manipulation. We note from the definition of $c_{n,q}$ in \eqref{eq:def.c_nq} that
\begin{equation*}
    c_{n,q}\geq\left(\frac{q\log(en_1)}{2n_1}\right)^{1-2/q}\Longrightarrow~\frac{(q-2)\log(en_1)}{n_1}\leq \frac{2(q-2)c_{n,q}^{q/(q-2)}}{q}\leq 2c_{n,q},
\end{equation*} where the last inequality follows from $c_{n,q}\leq 1$. Since $x\mapsto x/(1-x)_+$ is increasing in positive reals, the right-hand side of \eqref{pf_lemA.1:1.6} can be bounded as
\begin{align*}
    &(\mbox{RHS of }\eqref{pf_lemA.1:1.6})\leq \sqrt{\frac{(q-2)\log(en_1)}{1-c_{n,q}}}\left(c_{n,q}+\frac{2c_{n,q}/(1-c_{n,q})}{1-2c_{n,q}/(1-c_{n,q})}\right)\\
    &\quad=\sqrt{\frac{(q-2)\log(en_1)}{1-c_{n,q}}}\left(c_{n,q}+\frac{2c_{n,q}}{1-3c_{n,q}}\right)\leq 11c_{n,q}\sqrt{q\log(en)},
\end{align*}as long as $c_{n,q}\leq 1/4$. Hence, we have
\begin{align}\label{eq:pf_lemA.1:2}
&\Pb\left(\Norm{\frac{1}{\sqrt{n_1}}\sum_{i\in\Ic_1}V_{i,\theta}}_\infty> \frac{t\sqrt{1-c_{n,q}}}{\sqrt{1+t^2/n_1}}\right)-\Pb\left(\norm{\check G}_\infty>t\right)\nonumber\\
&\qquad\leq \Pb\left(\norm{\check G}_\infty>\frac{t\sqrt{1-c_{n,q}}}{\sqrt{1+t^2/n_1}}\right)-\Pb\left(\norm{\check G}_\infty>t\right)+\varrho_n^\theta\nonumber\\
&\qquad\leq\Phi_{\rm AC}^\theta\left(t-\frac{t\sqrt{1-c_{n,q}}}{\sqrt{1+t^2/n_1}}\right)+\varrho_n^\theta\leq 11\sqrt{qc_{n,q}^2\log(ed)\log(en)}+\varrho_n^\theta.
\end{align} Combining \eqref{eq:pf_lemA.1:1} and \eqref{eq:pf_lemA.1:2}, we get
\begin{align*}
    &\sup_{t\geq 0}\left[\Pb\left(\Norm{\frac{1}{\sqrt{n_1}}\sum_{i\in\Ic_1}V_{i,\theta}}_\infty> \frac{t\sqrt{1-c_{n,q}}}{\sqrt{1+t^2/n_1}}\right)-\Pb\left(\norm{\Av\check G}_\infty>t\right)\right]_+\\
    &\qquad\leq\varrho_n^\theta+\frac{d}{n^{q/2-1}}+Cq^{3/2-2/q}L_q(\theta)^2\left(\frac{\log^{2-2/q}(en)\log(ed)}{n^{1-2/q}}\right). 
\end{align*} Combining this with \eqref{eq:pf_lemA.1:0.1} and \eqref{eq:pf_lemA.1:0} proves the lemma.
\end{proof}

\begin{proof}[proof of Proposition~\ref{prop:B1}]
To prove \eqref{eq:prop:B1:1}, we use Proposition~\ref{prop:B1} which refines Theorem~2.5 of \cite{chernozhuokov2022improved}. For the application, we need to verify the conditions for the moments. First, we have
\begin{align}\label{eq:pf.prop.A.1:2}
    &\max_{j\in[d]}\left(\frac{1}{n_1}\sum_{i\in\Ic_1}\Eb[|V_{ij,\theta}|^q]\right)^q \leq L_q(\theta)
\end{align} Therefore, \eqref{eq:pf.prop.A.1:2} (with $q=4$) verifies the condition \eqref{eq:propB.1:1} for $B_n=L_4(\theta)$. Furthermore, for $q>2$,
\begin{align}\label{eq:pf.prop.A.1:3}
    &\max_{i\in\Ic_1}\left(\Eb\left[\max_{j\in[d]}\abs{V_{ij,\theta}}^q\right]\right)^{1/q}\leq \max_{i\in\Ic_1}\left(\sum_{j=1}^{d}\Eb\left[\abs{V_{ij,\theta}}^q\right]\right)^{1/q}\leq d^{1/q}\max_{i\in\Ic_1}\max_{j\in[d]}\left(\Eb\left[\abs{V_{ij,\theta}}^q\right]\right)^{1/q}\nonumber\\
    &\qquad\leq d^{1/q}L_q(\theta).
\end{align} Hence, \eqref{eq:pf.prop.A.1:2} verifies condition \eqref{eq:propB.1:2} with $D_n=d^{1/q}L_4(\theta)$, respectively. Hence, Proposition~\ref{prop:B1} applies and yields
\begin{equation*}
    \varrho_n^\theta\lesssim\frac{L_4(\theta)\log^{5/4}(ed)}{n^{1/4}}+\frac{d^{1/q}L_q(\theta)\log^{3/2}(ed)}{n^{1/2-1/q}},
\end{equation*} where $\lesssim$ only hides universal constants. This proves \eqref{eq:prop:B1:1}.

In order to prove \eqref{eq:prop:B1:2}, we use Theorem~2.1 of \cite{chernozhukov2023nearly} which reads 
\begin{equation}\label{eq:cher2023}
    \varrho_n^\theta\lesssim\log(en_1)\left\{\sqrt{\Delta_1\log(ed)}+\frac{(\Mc\log(ed))^2}{n_1\lambda(\theta)}\right\}+\sqrt{\frac{\Lambda_1M(\psi)}{n_1\lambda(\theta)^2}}+\frac{\log^{3/2}(ed)\psi}{\sqrt{n_1\lambda(\theta)}},
\end{equation} for all $\psi>0$ and $\lambda(\theta)=\lambda_{\rm min}(\Gamma(\theta;Q))$. Moreover, the parameters in \eqref{eq:cher2023} are defined as follows:
\begin{eqnarray*}
    \Delta_1 = \frac{\log^2(ed)}{n_1^2\lambda(\theta)^2}\max_{j\in[d]}\sum_{i\in\Ic_1}\Eb[(V_{ij,\theta})^4],&\quad&\Mc =\left(\Eb\left[\max_{i\in\Ic_1}\norm{V_{i,\theta}}_\infty^4\right]\right)^{1/4},\\
    \Lambda_1 = \log^2(ed)\log(en_1)\log(nd),&\quad& M(\psi)=\max_{i\in\Ic_1}\Eb\left[\norm{V_{i,\theta}}_\infty^4\mathbbm{1}(\norm{V_{i,\theta}}_\infty>\psi)\right].
\end{eqnarray*} First, $\Delta_1$ can be controlled using the fourth Lyapunov ratio as
\begin{equation*}
    \max_{j\in[d]}\sum_{i\in\Ic_1}\Eb[(V_{ij,\theta})^4]\leq n_1L_4^4(\theta).
\end{equation*} This implies that
\begin{equation}\label{eq:control_Delta1}
    \Delta_1 \leq \frac{\log^2(ed)L_4^4(\theta)}{n_1\lambda(\theta)^2}.
\end{equation} 
To control $\Mc$, we observe
\begin{align}\label{eq:control_Mc}
    &\Mc\leq \left(\Eb\left[\max_{i\in\Ic_1}\max_{j\in[d]}\abs{V_{ij,\theta}}^q\right]\right)^{1/q}\leq\left(\Eb\left[\sum_{i\in\Ic_1}\sum_{j=1}^{d}\abs{V_{ij,\theta}}^q\right]\right)^{1/q}\nonumber\\
    &\quad\leq (n_1d)^{1/q}\max_{j\in[d]}\left(\frac{1}{n_1}\sum_{i\in\Ic_1}\Eb[\abs{V_{ij,\theta}}^q]\right)^{1/q}\leq (n_1d)^{1/q}L_q,
\end{align} for any $q\geq 2$. Finally, for $M(\psi)$, we note that $x^4\mathbbm{1}(x\geq t)\leq x^{4+\delta}/t^\delta$ for any $x,t,\delta\geq 0$. This implies that
\begin{align}\label{eq:control_Mpsi}
    M(\psi)\leq \psi^{4-q}\max_{i\in\Ic_1}\Eb\left[\max_{j\in[d]}\abs{V_{ij,\theta}}^q\right]\leq\frac{dL_q^q}{\psi^{q-4}},
\end{align} for any $q\geq 4$, where the last inequality follows from \eqref{eq:control_Mc}. Combining \eqref{eq:cher2023}, \eqref{eq:control_Delta1}, \eqref{eq:control_Mc} and \eqref{eq:control_Mpsi}, we get
\begin{align}
    &C^{-1}\varrho_n^\theta\leq \frac{L_4^2(\theta)\log^{3/2}(ed)\log(en)}{n^{1/2}\lambda(\theta)}+\frac{L_q^2(\theta)\log^{2}(ed)d^{2/q}}{n^{1-2/q}\lambda(\theta)}\nonumber\\
    &\qquad\qquad+\frac{L_q^{q/2}(\theta)\Lambda_1d^{1/2}}{n^{1/2}\lambda(\theta)\psi^{q/2-2}}+\frac{\log^{3/2}(ed)\psi}{n^{1/2}\lambda^{1/2}(\theta)},
\end{align} for some universal constant $C>0$ and any $\psi>0$. The optimal choice of $\psi$ is
\begin{equation*}
    \psi = \left(\frac{L_q^{q/2}(\theta)\Lambda_1d^{1/2}}{\log^{3/2}(ed)\lambda^{1/2}(\theta)}\right)^{2/(q-2)}.
\end{equation*} This leads to
\begin{align}\label{eq:prop.A.1_res}
    &C^{-1}\varrho_n^\theta\leq \frac{L_4^2(\theta)\log^{3/2}(ed)\log(en)}{n^{1/2}\lambda(\theta)}+\frac{L_q^2(\theta)\log^{2}(ed)d^{2/q}}{n^{1-2/q}\lambda(\theta)}\nonumber\\
    &\qquad+\frac{L_q^{q/(q-2)}(\theta)\Lambda_{n,q}d^{1/(q-2)}}{\lambda^{q/(2q-4)}(\theta)n^{1/2}},
\end{align}where $\Lambda_{n,q}=[\log^{\frac{3q-4}{4}}(ed)\log(en)\log(nd)]^{\frac{2}{q-2}}.$ Therefore, the right-hand side of \eqref{eq:prop.A.1_res} can be bounded by
\begin{eqnarray}\label{eq:propC1:final}
    C^{-1}\varrho_n^\theta\leq \frac{L_q^2(\theta)\log^{2}(ed)d^{2/q}}{n^{1-2/q}\lambda(\theta)}+\frac{L_q^{q/(q-2)}(\theta)\Lambda'_{n,q}d^{1/(q-2)}}{\lambda^{q/(2q-4)}(\theta)n^{1/2}},
\end{eqnarray} with $\Lambda'_{n,q}=\Lambda_{n,q}\vee\log^{3/2}(ed)\log(en)\leq \log^{2}(ed)\log(en)\log(nd)$. Finally, since we may assume that
\begin{equation*}
    \frac{L_q^{q/(q-2)}(\theta)\Lambda'_{n,q}d^{1/(q-2)}}{\lambda^{q/(2q-4)}(\theta)n^{1/2}}\leq 1,
\end{equation*}otherwise, Proposition~\ref{prop:B1} holds by simply taking a large constant, the rightmost term in \eqref{eq:propC1:final} dominates the bound since $q\geq 4$. This proves \eqref{eq:prop:B1:2}.
\end{proof}

\begin{proof}[proof of Theorem~\ref{thm:B2}.] For $g\sim\Nc(0_d,I_d)$, let $\widehat S(t;\Gamma^\dagger)=\Pb(\norm{[\Gamma^\dagger]^{1/2}g}_\infty>t|\Gamma^\dagger)$. The one-sided DKW inequality (Corollary 1 of \cite{massart1990}) applies conditionally on $\Gamma^\dagger$, and yields that $$\Pb\left(\sup_{t\geq 0}[\hat S(t;\Gamma^\dagger)-\frac{1}{B}\sum_{i=1}^B\mathbbm{1}(\norm{[\Gamma^\dagger]^{1/2}W_i^{\rm b}}_\infty>t)]_+\leq\sqrt{\frac{\log(en)}{2B}}|\Gamma^\dagger\right)\geq 1-\frac{1}{n}.$$ Therefore, $\Pb[\Ec|\Gamma^\dagger]\geq 1-1/n$ where $\Ec =\set{\sup_{\alpha\in(0,1)}[\hat S( K_{\alpha,B}(\Gamma^\dagger);\Gamma^\dagger)-\alpha]_+\leq\sqrt{\frac{\log(en)}{2B}}}.$ Applications of Proposition~2.1 of \cite{chernozhuokov2022improved} and Theorem 2.1 of \cite{lopes2022central} yields
\begin{equation}\label{eq:pfthm16.1}
    \sup_{t\in\Real}\abs{\hat S(t;\Gamma^\dagger)-S(t)}\leq \begin{cases}
        c\log(en)\Upsilon(\Gamma^\dagger)^{1/2}&\mbox{w.p.1},\\
        c\lambda^{-1}\log(en)\Upsilon(\Gamma^\dagger)[1\vee\log(1/\Upsilon(\Gamma^\dagger))&\mbox{on}~\Ec_\lambda;
    \end{cases}.
\end{equation}Consider an event $\Mc_\delta:=\set{\Upsilon(\Gamma^\dagger)\leq \delta}$. To prove \eqref{eq:thm:B2:1}, we note that on $\Ec\cap\Mc_\delta$, $$\sup_{\alpha\in(0,1)}[S( K_{\alpha,B}(\Gamma^\dagger))-\alpha]_+\leq\sqrt{\frac{\log(en)}{2B}}+c\log(en)\delta^{1/2}.$$ Since $S(\cdot)$ is non-increasing, this implies that $ K_{\alpha,B}(\Gamma^\dagger)\geq S^{-1}(\alpha+\sqrt{\frac{\log(en)}{2B}}+c\log(en)\delta^{1/2})$ on $\Ec\cap\Mc_\delta$, where the right-hand side is $\Dc_2$-measurable. Therefore, we conclude that, w.p.1, uniformly for $\alpha\in(0,1)$, $\Pb(T^{\rm SeN}(\theta_0)>K_{\alpha,B}(\Gamma^\dagger)|\Dc_2 )\leq\Pb(T^{\rm SeN}(\theta_0)>S^{-1}(\alpha+\sqrt{\log(en)/(2B)}+c\log(en)\delta^{1/2})|\Dc_2 ) +\Pb((\Ec\cap\Mc_\delta)^\complement|\Dc_2)\leq \varepsilon_{n,q}^{(1)}+\alpha+\sqrt{\log(en)/(2B)}+c\log(en)\delta^{1/2}+\Pb((\Ec\cap\Mc_\delta)^\complement|\Dc_2)$ where the last inequality follows from Theorem~\ref{thm:B1}. The expected value of the last probability term is bounded by $1/n+\Pb(\Upsilon(\Gamma^\dagger)>\delta)$ by the union bound. Hence, the law of total probability implies \eqref{eq:thm:B2:1}. The inequality \eqref{eq:thm:B2:2} can be proven similarly with a sharper inequality in \eqref{eq:pfthm16.1} on $\Ec_\lambda$. We omit this for brevity.
\end{proof}
\begin{proof}[proof of Lemma~\ref{lem:B1}] We begin with simple triangle inequality that $$\mnorm{\hat\Delta}_\infty\leq \mnorm{\Gamma(\theta_0,\widetilde Q_2)-\Gamma(\widetilde\theta_2;\widetilde Q_2)}_\infty+ \mnorm{\Gamma(\widetilde\theta_2;\widetilde Q_2)-\widehat\Gamma_1(\widetilde\theta_2;\widetilde Q_2)}_\infty.$$ To control the leading term of the right-hand side, we use Lemma~\ref{lem:C6} to get $\mnorm{\Gamma(\theta_0,\widetilde Q_2)-\Gamma(\widetilde\theta_2;\widetilde Q_2)}_\infty\leq 4\norm{V(\theta_0)^{-1/2}(V(\widetilde\theta_2)-V(\theta_0))V(\theta_0)^{-1/2}}_{\rm op}\leq 4\phi(\norm{\widetilde\theta_2-\theta_0}_2)$. Now, an application of the following lemma (conditional on $\Dc_2$) completes the proof.
\end{proof}
\begin{lemma}\label{lem:B3} For any $\theta\in\Theta$, let 
\begin{equation}\label{eq:def_Anmq}
    A_{n,q}(\theta) = L_q^2(\theta)\frac{\log^{1-2/q}(ed)d^{2/q}}{n^{1-2/q}}.
\end{equation} For a non-stochastic $Q\in\Real^{d\times d}$ and $q\geq4$, there a constant $C=C(q)>0$ such that for any $\varepsilon\in(0,1)$,
\begin{equation*}
    \Pb\left(\mnorm{\Gamma(\theta;Q)-\widehat\Gamma_1(\theta;Q)}_\infty\geq C A_{n,q}(\theta)\left(1+\frac{1}{\varepsilon^{2/q}\log^{1-2/q}(ed)}\right)\right)\leq \varepsilon,
\end{equation*} Therefore, $\mnorm{\Gamma(\theta;Q)-\widehat\Gamma_1(\theta;Q)}_\infty=O_\Pb( A_{n,q}(\theta))$.
\end{lemma}

\begin{proof}[proof of Lemma~\ref{lem:B3}] Let $a_j=Q^\top\ev_j$ for $j\in[d]$. Lemma~\ref{lem:C6} implies that
\begin{align*}
    &\mnorm{\Gamma(\theta;Q)-\widehat\Gamma_1(\theta;Q)}_\infty\leq 4\max_{j,k\in[d]}\Abs{\frac{a_j^\top(V(\theta)-\widehat{V}_1(\theta))a_k}{(a_j^\top V(\theta))a_j)^{1/2}(a_k^\top V(\theta))a_k)^{1/2}}}\\
    &\quad\leq4\max_{j,k\in[d]}\Abs{\frac{a_j^\top(V(\theta)-\widehat{V}_1(\theta))a_k}{(\sum_{i\in\Ic_1}{\rm Var}[a_j^\top\psi(Z_i,\theta)])^{1/2}(\sum_{i\in\Ic_1}{\rm Var}[a_k^\top\psi(Z_i,\theta)])^{1/2}}}=:4\Xi_n.
\end{align*} Fix $\theta\in\Theta$, and denote
\begin{equation*}
    y_{ij}=\frac{a_j^\top(\psi(Z_i,\theta)-\mu(\theta))}{\sqrt{\sum_{i\in\Ic_1}{\rm Var}[a_j^\top\psi(Z_i,\theta)]}},
\end{equation*} for all $i\in\Ic_1$ and $j\in[d]$. Then, we may write
\begin{align*}
    \Xi_n&\leq\max_{j,k\in[d]}\Abs{\sum_{i\in\Ic_1}\Set{y_{ij}y_{ik}-\Eb[y_{ij}y_{ik}]}}+\max_{j\in[d]}\frac{n_1(a_j^\top(\widebar\psi_1(\theta)-\mu(\theta)))^2}{\sum_{i\in\Ic_1}{\rm Var}[a_j^\top\psi(Z_i,\theta)]}\\
    &=\max_{j,k\in[d]}\Abs{\sum_{i\in\Ic_1}\Set{y_{ij}y_{ik}-\Eb[y_{ij}y_{ik}]}}+\frac{1}{n_1}\max_{j\in[d]}\Abs{\sum_{i\in\Ic_1}y_{ij}}^2\\
    &=:\Iv_1+n_1^{-1}\Iv_2^2.
\end{align*} The proof is complemented by Lemma~\ref{lem:B4}, which establishes tail inequalities for $\Iv_1$ and $\Iv_2$.

\end{proof}
\begin{lemma}\label{lem:B4}The following concentration inequalities hold. The constant $C_q$ only depends on $q$ and may differ by lines.

\begin{enumerate}[left=0cm]
    \item For any $q\geq4$ and $\varepsilon\in(0,1)$, 
    \begin{equation*}
        \Pb\left(C_q^{-1}\Iv_1\geq\frac{L_q^2(\theta)\log^{1-2/q}(ed)d^{2/q}}{n^{1-2/q}}+\frac{1}{\varepsilon^{2/q}}\left[\frac{L_4^4(\theta)d^{4/q}}{n}+\frac{L_q^2(\theta)d^{2/q}}{n^{1-2/q}}\right]\right)\leq\varepsilon.
    \end{equation*}\label{lem:B4.1}
    \item For any $q\geq2$and $\varepsilon\in(0,1)$, 
    \begin{equation*}
    \Pb\left(C_q^{-1}n_1^{-1}\Iv_2^2\geq\frac{L_q^2(\theta)\log^{2-2/q}(ed)d^{2/q}}{n^{2-2/q}}+\frac{1}{\varepsilon^{2/q}}\left[\frac{L_q^4(\theta)d^{4/q}}{n^3}+\frac{L_q^2(\theta)d^{2/q}}{n^{2-2/q}}\right]\right)\leq\varepsilon.
\end{equation*}\label{lem:B4.3}
\end{enumerate}    
\end{lemma}

Proof strategy follows. To control the maximal average of mean-zero independent random variables, we begin with Markov's inequality: for any $r\geq1$,
\begin{equation*}
    \Pb\left(\Iv_{m}\geq\Eb[\Iv_m]+t\right)=\Pb\left((\Iv_{m}-\Eb[\Iv_m])_+\geq t\right)\leq t^{-r}\Eb[(\Iv_{m}-\Eb[\Iv_m])_+^r],
\end{equation*} for all $t>0$. Hence, it suffices to control $\Eb[\Iv_m]$ and $\Eb[(\Iv_{m}-\Eb[\Iv_m])_+^r]$ for some $r\geq1$. We use Proposition~B.1 of \cite{kuchibhotla2022least} to control the expectation and Theorem~15.14 of \cite{Boucheron2013concentration} to control the deviation parts. For a general version of this result, see Proposition~\ref{prop:C3}.

To facilitate the analysis, we compute the average of moments of $y_{ik}$. For any $q\geq2$
\begin{equation}\label{eq:moment_of_y_and_z}
    \max_{j\in[d]}\frac{1}{n_1}\sum_{i\in\Ic_1}\Eb[\abs{y_{ij}}^q]\leq n_1^{-q/2}L_q^q(\theta).  
\end{equation}

\begin{proof}[proof of Lemma~\ref{lem:B4}]
    In order to prove Lemma~\ref{lem:B4}\eqref{lem:B4.1}, we let $\xi_{i,kl}=y_{ik}y_{il}-\Eb[y_{ik}y_{il}]$. To apply Proposition~B.1 of \cite{kuchibhotla2022least}, we need to control some moments of $\xi_{i,kl}$. We first observe that
\begin{align}\label{eq:pf_lem3:1}
    &\max_{k,l\in[d]}\sum_{i\in\Ic_1}{\rm Var}(\xi_{i,kl})\leq\max_{k,l\in[d]}\sum_{i\in\Ic_1}\Eb[y_{ik}^2y_{il}^2]\leq\max_{k,l\in[d]}\sum_{i\in\Ic_1}\sqrt{\Eb[y_{ik}^4]\Eb[y_{il}^4]}\nonumber\\
    &\leq\max_{k,l\in[d]}\sqrt{\sum_{i\in\Ic_1}\Eb[y_{ik}^4]\sum_{i\in\Ic_1}\Eb[y_{il}^4]}=\max_{k\in[d]}\sum_{i\in\Ic_1}\Eb[y_{ik}^4]=n_1^{-1}L^4_4(\theta).
\end{align} Furthermore, one has
\begin{align}\label{eq:pf_lem3:2}
    &\sum_{i\in\Ic_1}\Eb\left[\max_{k,l\in[d]}\abs{\xi_{i,kl}}^{q/2}\right]\overset{(a)}{\leq} 2\sum_{i\in\Ic_1}\Eb\left[\max_{k,l\in[d]}\abs{y_{ik}y_{il}}^{q/2}\right]=2\sum_{i\in\Ic_1}\Eb\left[\max_{k\in[d]}\abs{y_{ik}}^{q}\right]\nonumber\\
    &\qquad\leq2d\max_{j\in[d]}\sum_{i\in\Ic_1}\Eb[\abs{y_{ik}}^q]=2dn^{1-q/2}L^q_q(\theta).
\end{align} Here, the inequality $(a)$ follows from the standard symmetrization lemma (see, for instance, Lemma~2.3.1 of \cite{van1996weak}). With \eqref{eq:pf_lem3:1} and \eqref{eq:pf_lem3:2}, Proposition~B.1 of \cite{kuchibhotla2022least} applies and yields
\begin{eqnarray}\label{eq:lem.3_pf:result_1}
    \Eb[\Iv_1]\lesssim\sqrt{\frac{L_4^4(\theta)\log(ed)}{n}}+\frac{L_q^2(\theta)\log^{1-2/q}(ed)d^{2/q}}{n^{1-2/q}},
\end{eqnarray} for any $q>2$ where $\lesssim$ only hides universal constants. Let $\kappa = \sqrt{e}/(2(\sqrt{e}-1))<1.271$ where $e$ is the natural constant. Theorem 15.14 of \cite{Boucheron2013concentration} with $r=q/2$ implies that for $q\geq2$,
    \begin{align}\label{eq:pf_lem3:3}
        &\left(\Eb[(\Iv_1-\Eb[\Iv_1])_+^r]\right)^{1/r}\leq \sqrt{\kappa q}\left[\Eb\left(\max_{k,l\in[d]}\sum_{i\in\Ic_1} \xi_{i,kl}^2\right)+\max_{k,l\in[d]}\sum_{i\in\Ic_1}\Eb [\xi_{i,kl}^2]\right]\nonumber\\
        &\qquad+~\kappa q\left[\left\{\Eb\left(\max_{i\in\Ic_1}\max_{k,l\in[d]}\abs{\xi_{i,kl}}^{q/2} \right)\right\}^{2/q}+\max_{i\in\Ic_1}\max_{k,l\in[d]}\left(\Eb[\xi_{i,kl}^2]\right)^{1/2}\right].
    \end{align}
We now analyze the quantities on the right-hand side of the last display. First, we get from \eqref{eq:pf_lem3:2} that
\begin{eqnarray}\label{eq:pf_lem3:4}
    \Eb\left(\max_{k,l\in[d]}\sum_{i\in\Ic_1} \xi_{i,kl}^2\right)&\leq& \sum_{i\in\Ic_1} \Eb\left(\max_{k,l\in[d]}\xi_{i,kl}^2\right)\leq \sum_{i\in\Ic_1} \left[\Eb\left(\max_{k,l\in[d]}\abs{\xi_{i,kl}}^{q/2}\right)\right]^{4/q}\nonumber\\
    &\leq& n_1^{1-4/q}\left[\sum_{i\in\Ic_1} \Eb\left(\max_{k,l\in[d]}\abs{\xi_{i,kl}}^{q/2}\right)\right]^{4/q}\leq \frac{2^{4/q}L^4_q(\theta)d^{4/q}}{n_1},
\end{eqnarray} for $q\geq2$. Furthermore, we get from \eqref{eq:pf_lem3:2} that
\begin{equation}\label{eq:pf_lem3:5}
    \left\{\Eb\left(\max_{i\in\Ic_1}\max_{k,l\in[d]}\abs{\xi_{i,kl}}^{q/2} \right)\right\}^{2/q}\leq\left\{\sum_{i\in\Ic_1}\Eb\left(\max_{k,l\in[d]}\abs{\xi_{i,kl}}^{q/2} \right)\right\}^{2/q}\leq \frac{2^{2/q}L_q^2(\theta)d^{2/q}}{n_1^{1-2/q}}.
\end{equation} Finally, it follows from the definition of $L_{n,2}$ that
\begin{align}\label{eq:pf_lem3:6}
    &\max_{i\in\Ic_1}\max_{k,l\in[d]}\Eb[\xi_{i,kl}^2]\leq \max_{i\in\Ic_1}\max_{k,l\in[d]}\Eb[y_{ik}^2y_{il}^2]\leq\max_{i\in\Ic_1}\max_{k,l\in[d]}\sqrt{\Eb[y_{ik}^4]\Eb[y_{il}^4]}\nonumber\\
    &\qquad=\max_{j\in[d]}\max_{i\in\Ic_1}\Eb[y_{ik}^4]=n_1^{-2}L_4^4(\theta).
\end{align}Combining \eqref{eq:pf_lem3:3}---\eqref{eq:pf_lem3:6}, we get
\begin{align}\label{eq:lem.3_pf:result_2}
    &\left(\Eb[(\Iv_1-\Eb[\Iv_1])_+^{q/2}]\right)^{2/q}\nonumber\\
    &\quad\lesssim \sqrt{q}\left[\frac{2^{4/q}L_q^4(\theta)d^{4/q}}{n}+\frac{L_4^4(\theta)}{n}\right]+ q\left[\frac{2^{2/q}L_q^2(\theta)d^{2/q}}{n^{1-2/q}}+\frac{L_4^2(\theta)}{n}\right].
\end{align} Therefore, combining \eqref{eq:lem.3_pf:result_1} and \eqref{eq:lem.3_pf:result_2} leads to that there exists a universal constant $C>0$ such that for any $q>2$,
\begin{align*}
    &\Pb\Bigg(C_q^{-1}\Iv_1\geq\sqrt{\frac{L_4^4(\theta)\log(ed)}{n}}+\frac{L_q^2(\theta)\log^{1-2/q}(ed)d^{2/q}}{n^{1-2/q}}\\
    &+\frac{1}{\varepsilon^{2/q}}\Set{\sqrt{ q}\left[\frac{L_q^4(\theta)d^{4/q}}{n}+\frac{L_4^4(\theta)}{n}\right]+ q\left[\frac{L_q^2(\theta)d^{2/q}}{n^{1-2/q}}+\frac{L_4^2(\theta)}{n}\right]}\Bigg)\leq \varepsilon,
\end{align*} for any $\varepsilon>0$. If $q\geq4$, this reduces to Lemma~\ref{lem:B4}\eqref{lem:B4.1}.

To control $\Iv_2$, we apply Proposition~\ref{prop:C3} to control $\max_{j\in[d]}|n_1^{-1/2}\sum_{i\in\Ic_1}y_{ij}|$. Since $L_{n,2}(\theta)\leq n^{1/q}L_q(\theta)$, Proposition~\ref{prop:C3} leads to that for any for any $\varepsilon\in(0,1)$, with probability at least $1-\varepsilon$,
\begin{equation}\label{eq:C.3.3.final}
    C_q^{-1}\max_{j\in[d]}\Abs{\frac{1}{\sqrt{n_1}}\sum_{i\in\Ic_1}y_{ij}}\leq\frac{L_q(\theta)\log^{1-1/q}(ed)d^{1/q}}{n^{1-1/q}}+\frac{1}{\varepsilon^{1/q}}\left\{\frac{L_q^2(\theta)d^{2/q}}{n^{3/2}}+\frac{L_q(\theta)d^{1/q}}{n^{1-1/q}}\right\},
\end{equation}for a constant $C_q$, depending only on $q$. This proves Lemma~\ref{lem:B4}\eqref{lem:B4.3}.

\end{proof}

\subsection{Proof of Theorem~\ref{thm:2}.}

\begin{theorem}\label{thm:B4}Let $\hat r_{n} = \norm{\widetilde\theta_2-\theta_0}_2$ and $$g(\widetilde\theta_2)= L_q^2(\widetilde\theta_2)[(\frac{\log(en)}{2n})^{1-2/q}+\frac{\log^{1-1/q}(ed)d^{1/q}}{n^{1-1/q}}],$$ and $$\Fk_n =\frac{\phi(\hat r_{n})\log(en)}{(1-2\phi(\hat r_{n}))_+} +\frac{g(\widetilde\theta_2)\log(en)}{(1-2g(\widetilde\theta_2))_+}.$$ Define $\Ec^{\rm SpN}:=\sup_{t\geq 0}[\Pb(T^{\rm SpN}(\theta_0)>t| \Dc_2)-S(t)]_+$. For $q>2$, there exists a constant $C=C(q)>0$ such that $\Ec^{\rm SpN}\leq C\varepsilon_{n,q}^{(1)} + C\Fk_n $ w.p.1. Moreover, for $q\geq 4$, $\Ec^{\rm SpN}\leq C\set{\varepsilon_{n,q}^{(2)}(\lambda)\mathbbm{1}_{\Ec_\lambda} +\varepsilon_{n,q}^{(1)} \mathbbm{1}_{(\Ec_\lambda)^\complement}} +C\Fk_n$  w.p.1, for any $\lambda>0.$
\end{theorem}
The remaining part of the proof proceeds similarly to Theorem~\ref{thm:1}, by leveraging bootstrap consistency results in Theorem~\ref{thm:B2} and concentration inequality for correlation matrix estimator in Lemma~\ref{lem:B1}. Thus, we only prove Theorem~\ref{thm:B4}.

\paragraph{proof of Theorem~\ref{thm:B4}.} Note that $T^{\rm SpN}(\theta_0) \leq A_1\times A_2 \times A_3$ where $$A_1=\max_{j}\sqrt{n_1}\frac{|\ev_j^\top Q\bar\psi_1(\theta_0)|}{\norm{\ev_j}_{V(\theta_0)}}, A_2 = \max_j\frac{\norm{\ev_j}_{V(\theta_0)}}{\norm{\ev_j}_{V(\widetilde\theta_2)}}, \mbox{ and }~A_3= \max_j\frac{\norm{\ev_j}_{V(\widetilde\theta_2)}}{\norm{\ev_j}_{\hat V_1(\widetilde\theta_2)}}.$$ Definition of $\phi$ implies $A_2\leq(1-\phi(\hat r_{n}))_+^{-1/2}$. Regarding the last term, Lemma~\ref{lem:C4} proves $\Pb(A_3^{-2}\leq 1-g(\widetilde\theta_2)|\Dc_2)\leq 2d/n_1.$ Putting together, we get for any $t>0$, $$\Pb(T^{\rm SpN}(\theta_0)>t | \Dc_2)\leq \Pb(A_1>\widetilde\gamma t |\Dc_2) + 2d/n_1,$$ where $\widetilde\gamma^2 = (1-g(\widetilde\theta_2))(1-\phi(\hat r_{n})).$ We write $\rho_n=\sup_{t\geq0}\abs{\Pb(A_1\geq t)-S(t)}$ and $\Phi_{\rm AC}=\sup_{t\geq 0, \varepsilon>0}\frac{1}{\varepsilon}[S(t)-S(t+\varepsilon)]$. For $t\geq0$, $\Pb(A_1\geq \widetilde\gamma t)\leq \rho_n + S(\widetilde\gamma t)\leq \rho_n +\Phi_{\rm AC}(1-\widetilde\gamma)t+S(t)$. Fix $t\geq 0$ and consider the event that $\set{\widetilde\gamma t\geq \sqrt{(q-2)\log(en_1)}}$. On this event, Gaussian maximal inequality in Lemma~\ref{lem:C1} implies that $S(\widetilde\gamma t)\leq S(\sqrt{(q-2)\log(en_1)})\leq d/n_1^{q/2-1}$. Meanwhile, on the complement  $\set{\widetilde\gamma t< \sqrt{(q-2)\log(en_1)}}$, we have $\Phi_{\rm AC}(1-\widetilde\gamma)t\leq \Phi_{\rm AC}(1-\widetilde\gamma)/\widetilde\gamma\sqrt{(q-2)\log(en_1)}$. Combining these, we get $$\Pb(A_1\geq t)-S(t)\leq \rho_n + \Phi_{\rm AC}(1-\widetilde\gamma)/\widetilde\gamma\sqrt{(q-2)\log(en_1)} +  2d/n_1^{q/2-1}.$$ Nazarov's inequality (\cite{nazarov2003maximal} or \cite{chernozhukov2017detailed}) gives that $\Phi_{\rm AC}\leq C\sqrt{\log(en)}$ for universal constant $C>0$. The bounds for $\rho_n$ results from high-dimensional CLTs \citep{chernozhukov2023nearly,chernozhuokov2022improved}, and can be found in Proposition~\eqref{prop:B1}. Finally, since $1-\sqrt{(1-a)(1-b)}\leq a+b$ and $x\mapsto x/(1-x)_+$ is increasing and convex in $x$, we conclude the proof.

\subsection{Proof of Theorem~\ref{thm:3}.}
Let $\hat a\in\Real^d$ be $\Dc_2$-measurable, and let
$$ B_1=n_1^{1/2}\frac{\abs{\hat a^\top\bar\psi_1(\theta_0)}}{\norm{\hat a}_{V(\theta_0)}},~B_2= \frac{\norm{\hat a}_V}{\norm{\hat a}_{V(\widetilde\theta_2)}},~\mbox{and}~B_3= \frac{\norm{\hat a}_{V(\widetilde\theta_2)}}{\norm{\hat a}_{\hat V_1(\widetilde\theta_2)}}.$$ Then, $T^{\rm GaN}(\theta_0) = B_1B_2 B_3$ with the choice $\hat a = \widetilde Q_2^\top(\theta-\widetilde\theta_2)$. Let $\Pb_2(\cdot) = \Pb(\cdot|\Dc_2)$ and $$\Ec^{\rm GaN}:=\sup_{t\geq 0}|\Pb_2(T^{\rm GaN}(\theta_0)\geq t)-2\Phi(-t)|.$$ For any $\epsilon,\epsilon_3>0$, $\Pb_2(T^{\rm GaN}(\theta_0)\geq t)\leq \Pb_2(B_1\geq (1+\epsilon)^{-1}(1+\epsilon_3)^{-1} t) + \Pb_2(|B_2-1|\geq\epsilon) +  \Pb_2(|B_3-1|\geq \epsilon_3),$ and, $\Pb_2(T^{\rm GaN}(\theta_0)\leq t)\leq \Pb_2(B_1\leq (1-\epsilon)^{-1}(1-\epsilon_3)^{-1} t) + \Pb_2(|B_2-1|\geq\epsilon) +  \Pb_2(|B_3-1|\geq \epsilon_3).$ Letting $\eta = \max\set{1-(1+\epsilon)^{-1}(1+\epsilon_3)^{-1}, (1-\epsilon)^{-1}(1-\epsilon_3)^{-1}-1}$, we get
\begin{equation}\label{eq:pf.thm3:1}
    \Ec^{\rm GaN}\leq C_{q\wedge 3}L_{q\wedge 3}^{q\wedge 3}/n^{(q\wedge 3-2)/2} + C|\eta|+ \Pb_2(|B_2-1|\geq\epsilon) +  \Pb_2(|B_3-1|\geq \epsilon_3),
\end{equation} where the inequality follows from Berry Esséen bound under $q<3$ moments (see, e.g., Theorem~6 in Chapter~3 of \cite{petrov1975sums}) and univariate Gaussian comparison inequality. The definition of $\phi$ implies $|B_2^{-2}-1|\leq \phi(\hat r_{n})$ for $\hat r_{n}=\norm{\widetilde\theta_2-\theta_0}_2$, and thus, $|B_2-1|\leq \phi(\hat r_{n})/(1-\phi(\hat r_{n}))_+=:\epsilon$ w.p.1. To prove the concentration of $B_3$ around 1, we consider $m$ i.i.d r.v.s $\Eb[X_i]=0$, $\Eb[X_i^2]=1$, and $\Eb^{1/q}[|X_i|^q]=M_q$ ($q> 3$) for $i\in[m]$. Let $\hat\sigma_X^2 = m^{-1}\sum_{i=1}^m (X_i-\bar X)^2 = m^{-1}\sum_{i=1}^m (X_i^2-1) -(\bar X)^2+1$. Chebyshev inequality implies $\Pb[(\bar X)^2\geq t]\leq 1/(mt)$, and Markov inequality and von Bahr-Esséen inequality \citep{vonbahr1965inequalities} imply that $\Pb(|m^{-1}\sum_{i=1}^m (X_i^2-1)|\geq t)\leq t^{-q/2} \Eb[|m^{-1}\sum_{i=1}^m (X_i^2-1)|^{q/2}]\leq (mt)^{-q/2} C_q\sum_{i=1}^m\Eb[|X_i^2-1|^{q/2}]\leq  C_q'M_q^qt^{-q/2}m^{1-q/2}$ for some constants $C_q,C_q'$ depending only on $q\in[2,4]$. Combining these yields $\Pb(|\hat\sigma_Z^2-1|\geq c'_{q\wedge 4} M_{q\wedge 4}^2/(m^{1-2/q\wedge 4}\delta^{2/{q\wedge 4}}))\leq \delta$ where $q>4$ case follows from Chebyshev inequality.

Applying to $X_i=\hat a ^\top\psi(Z_i,\widetilde\theta_2)/(\hat a^\top V(\widetilde\theta_2)\hat a)^{1/2}$ ($i\in\Ic_1$) implies that $\Pb_2(|B_3^{-2}-1|\geq \eta_3(\delta))\leq \delta$ where
$$\eta_3(\delta)=C_{q\wedge 4}\frac{L_{q\wedge 4}^2\gamma^2_{q\wedge 4}(\hat r_{n})}{n^{1-2/q\wedge 4}\delta^{2/{q\wedge 4}}}.$$
Therefore, $\Pb_2(|B_3-1|\geq \epsilon_3(\delta))\leq \delta$ for $\epsilon_3(\delta)= \eta_3(\delta)/(1-\eta_3(\delta))_+.$ We take $$\delta_0 = (L_{q\wedge 4}^2\gamma^2_{q\wedge 4}(\hat r_{n})/n^{1-2/(q\wedge 4)})^{\frac{q\wedge 4}{2+q\wedge 4}},$$ then $\epsilon_3 = \epsilon_3(\delta_0)\asymp\delta_0$. We may assume $\epsilon,\epsilon_3\leq 1/3$, otherwise the bounds in Theorem~\ref{thm:3} holds by taking a large constant. As long as $\epsilon,\epsilon_3\leq 1/3$, $\eta\leq 3(\epsilon+\epsilon_3)$ follows from simple algebraic manipulation, and hence, \eqref{eq:pf.thm3:1} becomes $\Ec^{\rm GaN}\leq C(\phi(\hat r_{n})+(L_{q\wedge 4}^2\gamma^2_{q\wedge 4}(\hat r_{n})/n^{1-2/(q\wedge 4)})^{\frac{q\wedge 4}{2+q\wedge 4}}),$ where the BE bound has been absorbed in the last term. The law of total probability proves Theorem~\ref{thm:3}.

\subsection{Proof of Theorem~\ref{thm:4}.}

For any $\theta_1\in\Theta$ and $\tau\in[0,1]$, define an interpolating path between $\theta_1$ and $\theta_0$ as $\theta_\tau:=\theta(\tau) = \tau\theta_1 + (1-\tau)\theta_0.$ Then, $\theta_\tau$ lies on the boundary of $\tau\norm{\theta_1-\theta_0}_2$-neighborhood of $\theta_0$. We denote the symmetrized Bregman divergence as $\hat\Dc_n(\theta,\theta') = (\theta-\theta')^\top (\bar\psi_1(\theta)-\bar\psi_1(\theta'))$. To enhance localized analysis of the minimizer of the convex function, we use Lemma~2 of \cite{Sun2020adaptivehuber} (or Lemma~F.2 of \cite{Fan2018ILAMM}) to get $\hat\Dc_n(\theta_\tau,\theta_0)\leq \tau\hat\Dc_n(\theta_1,\theta_0).$ This implies that $(\theta_\tau-\theta_0)^\top(\bar\psi_1(\theta_\tau)-\bar\psi_1(\theta_0))\leq \tau(\theta_1-\theta_0)^\top(\bar\psi_1(\theta_1)-\bar\psi_1(\theta_0)) = (\theta_\tau-\theta_0)^\top(\bar\psi_1(\theta_1)-\bar\psi_1(\theta_0)).$ Consequently, we get $(\theta_\tau-\theta_0)^\top \bar\psi_1(\theta_\tau)\leq (\theta_\tau-\theta_0)^\top \bar\psi_1(\theta_1)$. Let $r$ be any positive real, no larger than $\norm{\theta_1-\theta_0}_2$, and taking $\tau = r/\norm{\theta_1-\theta_0}_2$ and $\eta_\tau = (\theta_\tau-\theta_0)/\norm{\theta_\tau-\theta_0}_2$ leads to that
\begin{equation}\label{eq:pf.thm:7and9:0}
    \inf_{\eta\in\Sb^{d-1}}\eta^\top \bar\psi_1(\theta_0+r\eta)\leq \eta_\tau^\top \bar\psi_1(\theta_\tau)\leq \eta_\tau^\top \bar\psi_1(\theta_1)\leq \norm{\bar\psi_1(\theta_1)}_2.
\end{equation} Now, let $r_n$ satisfy \eqref{eq:thm:4:1}. Assumption~\ref{asmp:C.1} and \ref{asmp:C.2} implies that $\Pb(\Yc_{1,n}(r_n))\geq 1-\varepsilon$ where $\Yc_{1,n}(r)=\{\inf_{\eta\in\Sb^{d-1}}\eta^\top \bar\psi_1(\theta_0+r\eta)\geq \underline{\mu}(r)-u_{1,n}(r,\varepsilon)-\norm{\bar\psi_1(\theta_0)}_2\}$. Moreover, Proposition~\ref{prop:C1} yields that $\Pb(\Yc_{2,n})\geq 1-\epsilon $ where $\Yc_{2,n}=\{\norm{\bar\psi_1(\theta_0)}_{V^{-1}}\leq c_q(\sqrt{\frac{d+\log(2/\varepsilon)}{n}}+\frac{L_q d^{1/2}}{\varepsilon^{1/q}n^{1-1/q}})\}$. Thus, on the event $\set{\norm{\hat\theta_1-\theta_0}_2>r_n}\cap\Yc_{1,n}(r_n)\cap\Yc_{2,n}$, taking $\theta_1 = \hat\theta_1$ and $r=r_n$ in \eqref{eq:pf.thm:7and9:0} leads to the contradiction. Therefore, $\Pb(\norm{\hat\theta_1-\theta_0}_2>r_n)\leq \Pb(\Yc_{1,n}^\complement\cup\Yc_{2,n}^\complement)\leq2\varepsilon$.
\subsection{Proof of Theorem~\ref{thm:5}.}
Let $s_n$ satisfy \eqref{eq:thm:5:1} for a large $C_q>0$ and let $\hat r_{n}=\norm{\hat\theta_2-\theta_0}_2$. Let $\Yc_0=\set{\exists \theta^\star\in \widehat{\mathrm{CI}}^{\rm SpN}(K_n):~\norm{\theta^\star-\theta_0}_2>s_n}$, and on event $\Yc_0\cap\Yc_{1,n}(s_n)\cap\Yc_{2,n}$, take $n=n$, $\theta_1=\theta^\star$, and $r=s_n$ in \eqref{eq:pf.thm:7and9:0} to get
\begin{equation}\label{eq:pf.thm5:1}
    \underline{\mu}(s_n)\leq \frac{u_{1,n}(s_n,\varepsilon) }{\sqrt{n}}+c_q\overline{\lambda}_V\bigg(\sqrt{\frac{d+\log(2/\varepsilon)}{n}}+\frac{L_q d^{1/2}}{\varepsilon^{1/q}n^{1-1/q}}\bigg) + \norm{\bar\psi_1(\theta^\star)}_2.
\end{equation}Since $\theta^\star\in\widehat{\mathrm{CI}}^{\rm SpN}(K_n)$, $\sqrt{n}\norm{\bar\psi_1(\theta^\star)}_2\leq \underline{\sigma}^{-2}\sqrt{n}\norm{\bar\psi_1(\theta^\star)}_{\widetilde Q^\top\widetilde Q}\leq \underline{\sigma}^{-2}K_n{\rm tr}^{1/2}(\widetilde Q\hat V_1(\hat\theta_2)\widetilde Q^\top)\leq (\overline{\sigma}/\underline{\sigma})^{2}K_n{\rm tr}^{1/2}(\hat V_1(\hat\theta_2)).$ Let $\widecheck V = V(\hat\theta_2)^{-1/2}\hat V_1(\hat\theta_2)V(\hat\theta_2)^{-1/2}$, then ${\rm tr}(\hat V_1(\hat\theta_2))\leq {\rm tr}(\widecheck V)\norm{V(\hat\theta_2)}_{\rm op}\leq {\rm tr}(\widecheck V)\overline{\lambda}_V(1+\phi(\hat r_{n}))$. Proposition~\ref{prop:C1} controls ${\rm tr}(\widecheck V)$, and implies $\Pb(\Yc_3|\Dc_2)\geq 1-\varepsilon$ where
\begin{equation}\label{eq:pf.thm5:2}
    \Yc_3 = \bigg\{\norm{\bar \psi_1(\theta^\star)}_2\leq c_q\frac{\overline{\sigma}^2\overline{\lambda}_V(1+\phi(\hat r_{n}))}{\underline{\sigma}^2}\bigg[\sqrt{\frac{K_n^2d}{n}}+\frac{\gamma_q(\hat r_{n})L_qK_nd^{1/2}}{\varepsilon^{1/q}n^{1-1/q}}\bigg]\bigg\}.
\end{equation} On $\Yc_0\cap\Yc_{1,n}(s_n)\cap\Yc_{2,n}\cap\set{\hat r_{n}\leq r_{n}}$, combining \eqref{eq:pf.thm5:1} and \eqref{eq:pf.thm5:2} contradicts to \eqref{eq:thm:5:1}. Therefore, $\Pb(\Yc_0)\leq \Pb(\Yc_{1,n}(s_n)^\complement)+\Pb(\Yc_{2,n}^\complement)+\Pb(\hat r_{n}>r_{n})\leq 4\varepsilon$, and this concludes the proof.

\subsection{Proof of Lemma~\ref{lem:1}.}
For $W_i^{\rm b}| Z_1,\ldots,Z_n\sim\Nc(0_d,\widecheck\Gamma)$, $i\in[n]$ and $\widecheck\Gamma=\widehat\Gamma_1(\hat\theta;\hat Q)$, write $\hat S_Z(t)=n^{-1}{\textstyle\sum}_{b=1}^n\mathbbm{1}(\norm{W_i^{\rm b}}_\infty>t)$ and recall $K^{\rm b}_\alpha(\hat\theta;\hat Q):=\inf\{t\geq0:\hat S_Z(t)\leq\alpha\}.$ For $g\sim\Nc(0_d,I_d)$, let $\widecheck S(t)=\Pb(\norm{[\widecheck\Gamma]^{1/2}g}_\infty>t|\widecheck\Gamma)$. The one-sided DKW inequality (Corollary 1 of \cite{massart1990}) applies conditionally on $\widecheck\Gamma$, and yields that $\Pb(\sup_{t\geq 0}[\hat S_Z(t)-\widecheck S(t)]_+\leq\sqrt{\log(en)/(2n)}|\Gamma^\dagger)\geq 1-1/n.$ Note that $\ev_j^\top\widecheck\Gamma^{1/2}g|\widecheck\Gamma\sim\Nc(0,1)$ for $j\in[d]$, and thus, the union bound implies $\widecheck S(K)\leq d(1-\Phi(K))\leq 2de^{-K^2/2}\leq\alpha/2$ for $K=\sqrt{2\log(4d/\alpha)}$. This implies $\Pb(K^{\rm b}_\alpha(\hat\theta;\hat Q)\geq K)\leq \Pb(\hat S_Z(K)>\alpha-1/n)\leq1/n+\Pb(\widecheck S(K)>\alpha/2)=1/n$.
\subsection{Proof of Theorem~\ref{thm:6}.}
By simple implication, $\Pb(\xi^\top\theta_0\notin\hat{\rm CI}^{\rm a}_\alpha)\leq \Pb(\theta_0\notin\hat{\mathscr{C}}^{\rm a}_\alpha) + \Pb(\theta_0\notin\widehat{\mathrm{CI}}^{\rm SpN}_{\alpha/n})$ for ${\rm a}\in\set{\rm SeN, SpN}$. The bound for $(\Pb(\theta_0\notin\widehat{\mathrm{CI}}^{\rm SpN}_{\alpha/n})-\alpha/n)_+$ is provided in Theorem~\ref{thm:2}. For SeN, Theorem~1.2 of \cite{bentkus1996berry} gives that $|\Pb(\theta_0\notin\hat{\mathscr{C}}_\alpha^{\rm SeN}|\Dc_2)-\alpha|\leq cL_{q\wedge 3}^{q\wedge 3}/n^{(q\wedge 3)/2-1}$, w.p.1, for a constant $c>0$ and $q>2$. For SpN, we use Theorem~\ref{thm:3} whereas a sharper one-side control of $B_3$ is possible; Lemma~\ref{lem:C4} proves $\Pb(B_3^{-2}\leq 1-cL_q^2(\widetilde\theta_2)(\log(2n)/n)^{1-2/q})\leq 2/n,$ for $q>2$. This implies a one-side analogue of \eqref{eq:pf.thm3:1}, $(\Pb(\theta_0\notin\hat{\mathscr{C}}_\alpha^{\rm SeN}|\Dc_2)-\alpha)_+\leq cL_{q\wedge 3}^{q\wedge 3}/n^{q\wedge3/2-1} + c\phi(\hat r_{n})+ c\gamma_q^2(\hat r_{n})L_q^2(\log(2n)/n)^{1-2/q}) + 2/n$, for $q>2$. This proves the theorem.

\subsection{Verification of Assumptions in Misspecified GLMs}
In this section, we verify assumptions in the context of misspecified GLMs. A constant $C>0$ may vary line by line, while we retain explicit constants in their proofs. We prove that 
\vspace{-0.2cm}
\begin{enumerate}[left=0cm]
\setlength{\itemsep}{-.3em}
    \item The $q$-th Lyapnouv ratio $L_q$ is uniformly bounded for $q\leq q_{xy}$ and \ref{asmp:m} holds.\label{_toglm:1}
    \item $\phi(\delta)=C\delta,~\delta\in(0,\delta_0)$ for some $\delta_0>0$.\label{_toglm:2}
    \item $\gamma_q(\delta)=1+C\delta,~\delta\in(0,\delta_0)$ for some $\delta_0>0$, for $q\leq \min\{q_x/2,q_{xy}\}$.\label{_toglm:3}
    \item \ref{asmp:C.1} holds with $\underline{\mu}(\delta)=C(\delta\wedge\delta_0)$ for all $\delta>0$ and some $\delta_0>0$.\label{_toglm:4}
    \item \ref{asmp:C.2} holds for all $\delta>0$, $\epsilon\in(0,1)$ with $\frac{u_{1,n}(\delta,\varepsilon)}{\sqrt{n}} = C\delta[ \left(\frac{d}{n}\right)^{\frac{1}{2}-\frac{1}{2(q_x-1)}}+\sqrt{\frac{\log(1/\varepsilon)}{n}}+\frac{\log(1/\varepsilon)}{(n^{q_x-2}d)^{1/(q_x-1)}}].$\label{_toglm:5}
\end{enumerate}
\vspace{-.5cm}
\paragraph{Verification of \ref{_toglm:1}.}
We first analyze the $q$-th Lyapunov ratio $L_q(\theta)$. At $\theta=\theta_0$, this can be controlled through \ref{asmp:glm2} and \ref{asmp:glm3}. For $q\leq q_{xy}:=(1/q_x+1/q_y)^{-1}$,
\begin{align}\label{eq:pf.glm:0}
    &\sup_{u\in\Sb^{d-1}}(\Eb[\abs{u^\top X(Y-\dot{\ell}(X^\top\theta_0))}^q])^{1/q}\leq \sup_{u\in\Sb^{d-1}}(\Eb[\abs{u^\top X(Y-\dot{\ell}(X^\top\theta_0))}^{q_{xy}}])^{1/q_{xy}}\nonumber\\
    &\quad\leq\sup_{u\in\Sb^{d-1}}\norm{\Sigma^{1/2}u}_2\sup_{v\in\Sb^{d-1}}(\Eb[\abs{v^\top \Sigma^{-1/2}X(Y-\dot{\ell}(X^\top\theta_0))}^{q_{xy}}])^{1/q_{xy}}\leq\overline{\lambda}_\Sigma^{1/2}K_xK_y.
\end{align} The first inequality is the Jensen inequality, and the last is H\"older inequality. In particular, $q=2$ verifies \ref{asmp:m}. Meanwhile, we observe that $\inf_{u\in\Sb^{d-1}}(\Eb[|u^\top X(Y-\dot{\ell}(X^\top\theta_0))|^2])^{1/2}\geq\underline{\lambda}_V^{1/2}.$ Combining these implies that $L_q(\theta_0)\leq (\overline\lambda_\Sigma/\underline\lambda_V)^{1/2}K_xK_y$ for $q\leq q_{xy}$.
\paragraph{Verification of \ref{_toglm:2}.}
For a general $\theta\in\Real^d$, we note the decomposition: $\psi_1(\theta)-\mu(\theta) = \psi_1(\theta_0)+\psi_1(\theta)-\psi_1(\theta_0)-\mu(\theta)$. Hence, it follows from Minkowski inequality that $ \abs{\Eb^{1/q}[\abs{u^\top(\psi_1(\theta)-\mu(\theta))}^q]-\Eb^{1/q}[\abs{u^\top\psi_1(\theta_0)}^q]}\leq \Eb^{1/q}[\abs{u^\top(\psi_1(\theta)-\psi_1(\theta_0)}^q] + |u^\top\mu(\theta)|^q.$ Meanwhile, \ref{asmp:glm1} implies $\abs{\dot{\ell}(X^\top\theta_0)-\dot{\ell}(X^\top\theta)}\leq b_1\abs{X^\top(\theta-\theta_0)}$, and thus, for $q\leq q_x/2$ and $u\in\Sb^{d-1}$,
\begin{equation*}
    \Eb^{1/q}[\abs{u^\top(\psi_1(\theta)-\psi_1(\theta_0)}^q]\leq b_1\Eb[|u^\top X|^{q}|X^\top(\theta-\theta_0)|^{q}]^{1/q}\leq \overline b_1\overline\lambda_\Sigma K_x^2\norm{\theta_0-\theta}_2.
\end{equation*} For $q\leq \min\set{q_{xy},q_x/2}$, we have $\abs{\Eb^{1/q}[\abs{u^\top(\psi_1(\theta)-\mu(\theta))}^q]-\Eb^{1/q}[\abs{u^\top\psi_1(\theta_0)}^q]}\leq2b_1\overline\lambda_\Sigma K_x^2\norm{\theta-\theta_0}_2$. Combining this with \eqref{eq:pf.glm:0}, we get
\begin{equation}\label{eq:pf.glm:1}
    \abs{\Eb^{2/q}[\abs{u^\top(\psi_1(\theta)-\mu(\theta))}^q]-\Eb^{2/q}[\abs{u^\top\psi_1(\theta_0)}^q]}\leq2b_1\overline\lambda_\Sigma K_x^2\norm{\theta-\theta_0}_2(2b_1\overline\lambda_\Sigma K_x^2\norm{\theta-\theta_0}_2 + 2\overline\lambda^{1/2}_\Sigma K_xK_y).
\end{equation} Together with \ref{asmp:glm3} verifies $\phi(\delta) = 4b_1\overline\lambda^{3/2}_\Sigma\underline{\lambda}_V^{-1} K_x^3K_y\delta$ for $\delta\in(0,b_1\overline{\lambda}_\Sigma^{1/2}K_x/K_y)$.

\paragraph{Verification of \ref{_toglm:3}}

For any $A,B,A',B'$ such that $|A-A'|\vee|B-B'|<\varepsilon$, we observe that
\begin{align*}
    \Abs{\frac{B}{A}-\frac{B'}{A'}}\leq\Abs{\frac{B-B'}{A}}+\Abs{\frac{(A'-A)B'}{AA'}}\leq \frac{\varepsilon}{A} + \frac{\varepsilon(B+\varepsilon)}{A(A-\varepsilon)}=\frac{\varepsilon(A+B)}{A(A-\varepsilon)},\mbox{ if }A>\varepsilon.
\end{align*} We put $A = \Eb[|u^\top \psi_1(\theta_0)|^2]$, $B=(\Eb[|u^\top \psi_1(\theta_0)|^q])^{2/q}$, $A' = \Eb[|u^\top(\psi_1(\theta)-\mu(\theta))|^2]$, and $B'=(\Eb[|u^\top (\psi_1(\theta)-\mu(\theta))|^q])^{2/q}$ to get from \eqref{eq:pf.glm:1} that $|A-A'|\vee|B-B'|\leq \varepsilon:= 4b_1\overline\lambda^{3/2}_\Sigma K_x^3K_y\norm{\theta_0-\theta}_2+4b_1^2\overline\lambda_\Sigma^2 K_x^4\norm{\theta_0-\theta}_2^2.$ Also, $A,B\leq \overline{\lambda}_\Sigma K_x^2K_y^2$ from \eqref{eq:pf.glm:0} and $ A\geq \underline\lambda_V$ from \ref{asmp:glm3}. Thus,
\begin{equation*}
    \Abs{L_q^2(\theta)-L_q^2}\leq \frac{2\overline{\lambda}_\Sigma K_x^2K_y^2\set{4b_1\overline\lambda^{3/2}_\Sigma K_x^3K_y\norm{\theta_0-\theta}_2+4b_1^2\overline\lambda_\Sigma^2 K_x^4\norm{\theta_0-\theta}_2^2}}{\underline\lambda_V(\underline\lambda_V-\set{4b_1\overline\lambda^{3/2}_\Sigma K_x^3K_y\norm{\theta_0-\theta}_2+4b_1^2\overline\lambda_\Sigma^2 K_x^4\norm{\theta_0-\theta}_2^2})_+}
\end{equation*} If $\norm{\theta-\theta_0}_2\leq \varepsilon_\circ:=\min\set{b_1\overline\lambda_\Sigma^{1/2}K_x/K_y,\underline{\lambda}_V/(9b_1\overline\lambda_\Sigma^{3/2}K_x^2K_y)}$, then for $q\leq \min\set{q_{xy},q_x/2}$, $
    |L_q^2(\theta)-L_q^2|\leq C\norm{\theta-\theta_0}_2$ where $C=144b_1\underline\lambda_V^{-2}\overline\lambda_\Sigma^{5/2}K_x^5K_y^3$, i.e., $\gamma_q(\delta)=1+C\delta$ for $\delta\in(0,\varepsilon_\circ)$.

\paragraph{Verification of \ref{_toglm:4}} For $\eta\in\Sb^{d-1},$ Taylor expansion of $\mu(\theta)$ at $\theta=\theta_0$ gives $\eta^\top [\mu(\theta_0+\delta\eta)-\mu(\theta_0)]=\delta(\eta^\top\nabla_\theta\mu(\theta_0)\eta)+\delta^2\int_0^1(1-s)\ip{\nabla_\theta^2\mu(\theta_0+s\delta\eta),\eta^{\otimes 3}}~ds\geq\delta(\lambda_{\rm min}(\Eb[XX^\top\ell''(X^\top\theta_0)]))-(\delta^2/2)\sup_{\theta\in \Bb(\theta_0,\delta)}\abs{\Eb[\ell'''(X^\top\theta)(\eta^\top X)^3]} \geq \underline\lambda_V\delta -b_1K_x^3\overline\lambda_\Sigma^{3/2}\delta^2/2$. Here, the last inequality follows from \ref{asmp:glm3} and that $\sup_{\theta\in \Bb(\theta_0,\delta)}\abs{\Eb[\ell'''(X^\top\theta)(\eta^\top X)^3]}\leq b_1\Eb[|\eta^\top X|^3]\leq b_1K_x^3\overline\lambda_\Sigma^{3/2}$. Therefore, $\eta^\top [\mu(\theta_0+\delta\eta)-\mu(\theta_0)]\geq\underline\lambda_V\delta/2$ for $\delta\in(0,b_2/(b_1K_x^3\overline{\lambda}_\Sigma^{3/2}))$. Otherwise, since $\partial [\eta^\top (\mu(\theta_0+\delta\eta)-\mu(\theta_0))]/\partial\delta = \eta^\top \Eb[XX^\top\ell''(X^\top(\theta_0+\delta\eta))]\eta\geq0$, followed by convexity of $\ell$, $\eta^\top [\mu(\theta_0+\delta\eta)-\mu(\theta_0)]$ is non-decreasing in $\delta$. This verifies \ref{_toglm:4}.
\paragraph{Verification of \ref{_toglm:5}} We shall prove the lemma below.

\begin{lemma}\label{lem:momentbound_G} Assume \ref{asmp:glm1}, \ref{asmp:glm2}, and \ref{asmp:glm3}. Then, \ref{asmp:C.2} is satisfied with
\begin{equation*}
    \frac{u_{1,n}(\delta,\varepsilon)}{\sqrt{n}} = Cb_1\delta\bigg[ K_x^2\overline\lambda_\Sigma^{1/2} \left(\frac{d}{n}\right)^{\frac{1}{2}-\frac{1}{2(q-1)}}+\sqrt{\frac{\overline{\lambda}_\Sigma\log(1/\varepsilon)}{n}}+\frac{K_x^4\overline{\lambda}_\Sigma\log(1/\varepsilon)}{(n^{q-2}d)^{1/(q-1)}}\bigg],
\end{equation*} for $q\in[2,q_x]$ and some universal constant $C>0$.
\end{lemma}
\begin{proof}
    For a given $\delta>0$, denote $g_{\eta,\delta}(x) = -(x^\top\eta)(\ell'(x^\top\theta_0 +\delta x^\top\eta)-\ell'(x^\top\theta_0)).$ It follows from the construction that $-\eta^\top(\bar\psi_1(\theta_0+\delta\eta)-\bar\psi_1(\theta_0)-\mu(\theta_0+\delta\eta))=(\Pb_n-P)g_{\eta,\delta}(X)$, where $\Pb_n$ denotes the uniform distribution on $\set{X_1,\ldots,X_{n/2}}.$ Hence, it suffices to (stochastically) upper-bound the quantity, $\sup_{\eta\in\Sb^{d-1}}(\Pb_n-P)g_{\eta,\delta}(X).$ We introduce a smooth truncation: Let $\varphi:\Real\to[0,1]$ be a mollified cosine cutoff function: $\psi(u)=1$ for $u\in[-1,1]$, $\frac{1+\cos(\pi\cdot(|u|-1))}{2}$ for $|u|\in(1,2)$, and 0 otherwise. From the construction, $\varphi$ is supported on $[-2,2]$ and $\norm{\varphi'}_\infty\leq \pi/2$, and thus, $\varphi(\cdot)$ is $\pi/2-$Lipschitz, that is, $|\varphi(u_1)-\varphi(u_2)|\leq \pi/2\cdot|u_1-u_2|$ for $u_1,u_2\in\Real$. For $B>0$, define $\varphi_B(u)=\varphi(u/B)$, then $\varphi_B$ is compactly supported on $[-2B, 2B]$, equals 1 on $[-B,B]$, and $\pi/(2B)-$Lipschitz. Equipped with this, we define $g^{B}_{\eta,\delta}(x) := g_{\eta,\delta}(x)\varphi_B(|\eta^\top x|).$ Since $g_{\eta,\delta}(\cdot)\leq 0$ due to the convexity of $\ell$, we have $g_{\eta,\delta}(x)-g^{B}_{\eta,\delta}(x) = g_{\eta,\delta}(x)(1-\varphi_B(\eta^\top x))\leq 0$, that is, $g_{\eta,\delta}(x)\leq g^B_{\eta,\delta}(x)\leq 0$ for all $x$. This implies that $(\Pb_n-P)g_{\eta,\delta}(X)\leq \Pb_n g^B_{\eta,\delta}(X)-Pg_{\eta,\delta}(X)= (\Pb_n-P) g^B_{\eta,\delta}(X)+P(g^B_{\eta,\delta}-g_{\eta,\delta})(X).$ Consequently, we get
    \begin{equation}\label{eq:B.93}
        \sup_{\eta\in\Sb^{d-1}}(\Pb_n-P)g_{\eta,\delta}(X)\leq \sup_{\eta\in\Sb^{d-1}}(\Pb_n-P) g^B_{\eta,\delta}(X) + \sup_{\eta\in\Sb^{d-1}}P(g^B_{\eta,\delta}-g_{\eta,\delta})(X).
    \end{equation} We first control the rightmost term by $0\leq (g^B_{\eta,\delta}-g_{\eta,\delta})(x) = -g_{\eta,\delta}(x)(1-\varphi_B(\eta^\top x))\leq -g_{\eta,\delta}(x)\mathbbm{1}(|\eta^\top x|>B).$ Since
    $\abs{g_{\eta,\delta}(x)}\leq \norm{\ell''}_\infty\delta(\eta^\top x)^2\leq b_1\delta(\eta^\top x)^2$, we get, for any $q\in[2,q_x],$
    \begin{equation}\label{eq:B.92}
        \sup_{\eta\in\Sb^{d-1}}\abs{P(g_{\eta,\delta}-g^B_{\eta,\delta})} \leq b_1\delta\sup_{\eta\in\Sb^{d-1}}\Eb[(\eta^\top X)^2\mathbbm{1}_{|\eta^\top X|>B}]\leq b_1\delta \sup_{\eta\in\Sb^{d-1}}\Eb\bigg[\frac{|\eta^\top X|^q}{B^{q-2}}\bigg]\leq\frac{b_1 K_x^q\overline\lambda_\Sigma^{q/2} \delta}{B^{q-2}}.
    \end{equation} To control the leading term of the right-hand side of \eqref{eq:B.93}, we write $Z=\sup_{\eta\in\Sb^{d-1}}\abs{(\Pb_n-P) g^B_{\eta,\delta}(X)}$ and we shall prove that
    \begin{align}
        &\qquad\Eb^*[Z]\leq (8+4\pi)b_1\overline\lambda_\Sigma^{1/2}\delta B\sqrt{2d/n},\quad\mbox{and}\label{eq:B.94}\\
        &\Pb^*\bigg(Z\geq 3\Eb[Z] + 3b_1\delta \sqrt{\frac{2\overline{\lambda}_\Sigma\log(1/\varepsilon)}{n}}+38b_1\delta\frac{B^2\log(1/\varepsilon)}{n}\bigg)\leq \varepsilon,\label{eq:B.95}
    \end{align} for all $\varepsilon\in(0,1)$.
    To prove \eqref{eq:B.94}, we let $f(u):=f(u;t,\delta,B) = -u(\ell'(t+\delta u)-\ell'(t))\varphi_B(u).$ The first derivative is given by $f'(u) = -(\ell'(t+\delta u)-\ell'(t))\varphi_B(u) -u\delta\ell''(t+\delta u)\varphi_B(u)-u(\ell'(t+\delta u)-\ell'(t))\varphi_B'(u).$ Hence, \ref{asmp:glm1} and the construction of $\varphi_B$ leads to that $\abs{f'(u)}\leq 2b_1\delta |u|\varphi_B(u) + b_1\delta u^2\varphi'_B(u)\leq (4+2\pi)b_1\delta B.$ Thus, for any $\eta_1,\eta_2\in\Sb^{d-1}$ and $x\in\Real^d$, $\abs{g_{\eta_1,\delta}^B(x)-g_{\eta_2,\delta}^B(x)}\leq (4+2\pi)b_1\delta B \abs{\eta_1^\top x - \eta_2^\top x}.$ By a standard symmetrization argument and Ledoux–Talagrand contraction inequality (see Equation~(4.20) in \cite{Ledoux2011Isoperimetry} or Theorem~11.6 in \cite{Boucheron2013concentration}), 
    \begin{align*}
        &\Eb^*\bigg[\sup_{\eta\in\Sb^{d-1}}(\Pb_n-P) g^B_{\eta,\delta}(X)\bigg]\leq 2\Eb^*\bigg[\sup_{\eta\in\Sb^{d-1}}\frac{1}{n/2}\sum_{i=1}^{n/2}\varepsilon_i^\circ g_{\eta,\delta}^B(X_i)\bigg]\\
        &\quad\leq(8+4\pi)b_1\delta B\Eb\bigg[\sup_{\eta\in\Sb^{d-1}}\frac{1}{n/2}\sum_{i=1}^{n/2}\varepsilon_i^\circ \eta^\top X_i\bigg]=(8+4\pi)b_1\delta B\Eb\bigg[\bigg\|{\frac{1}{n/2}\sum_{i=1}^{n/2}\varepsilon_i^\circ X_i}\bigg\|_2\bigg].
    \end{align*} The last term can be bounded by  
    \begin{equation*}
        \Eb^{1/2}\bigg[\Norm{\frac{1}{n/2}\sum_{i=1}^{n/2}\varepsilon_i^\circ X_i}_2^2\bigg] = \frac{{\rm tr}^{1/2}\left({\rm Var}\left(\varepsilon^\circ X_1\right)\right)}{\sqrt{n/2}}=\frac{{\rm tr}^{1/2}\left({\rm Var}\left( X_1\right)\right)}{\sqrt{n/2}}\leq\overline{\lambda}^{1/2}_\Sigma\sqrt{\frac{2d}{n}}.
    \end{equation*}
    This proves \eqref{eq:B.94}. To prove \eqref{eq:B.95}, we note that $g_{\eta,\delta}^B$ is uniformly bounded as $|g_{\eta,\delta}^B(x)|\leq |g_{\eta,\delta}(x)\mathbbm{1}(|\eta^\top x|\leq B)|\leq b_1\delta (\eta^\top x)^2\mathbbm{1}(|\eta^\top x|\leq B)|\leq b_1\delta B^2$. Moreover, we have ${\rm Var}(g_{\eta,\delta}^B(X_i))\leq \Eb[g_{\eta,\delta}^B(X_i)^2]\leq \Eb[g_{\eta,\delta}(X_i)^2]\leq b_1^2\delta^2\Eb[(\eta^\top X_i)^2]\leq b_1^2\delta^2\overline{\lambda}_\Sigma$ uniformly in $\eta\in\Sb^{d-1}$. Hence, Bennett's inequality for bounded empirical processes (see, for instance, Equation~(11) of \cite{massart1990} with $\kappa=4$ and $\epsilon=2$) yields \eqref{eq:B.95}. Finally, combining \eqref{eq:B.93}---\eqref{eq:B.95} leads to that there exists a universal constant $C>0$ such that for any $\varepsilon>0$, w.p. at least $1-\varepsilon$, 
    \begin{equation*}
        \sup_{\eta\in\Sb^{d-1}}(\Pb_n-P)g_{\eta,\delta}(X)\leq Cb_1\delta\bigg[\frac{ K_x^q\overline\lambda_\Sigma^{q/2}}{B^{q-2}} + \overline\lambda_\Sigma^{1/2} B\sqrt{\frac{d}{n}}+\sqrt{\frac{\overline{\lambda}_\Sigma\log(1/\varepsilon)}{n}}+\frac{B^2\log(1/\varepsilon)}{n}\bigg],
    \end{equation*} for any $q\in [2,q_x]$. We take $B = K_x^{\frac{q}{q-1}}\overline{\lambda}_\Sigma^{1/2}(n/d)^{\frac{1}{2(q-1)}}$, and this proves the lemma.
    

\end{proof}

\subsection{Proof of Theorem~\ref{thm:7} and \ref{thm:12}.} We may assume $d\leq cn$ for sufficiently small constant; otherwise the results follow by taking a large constant. Theorem~\ref{thm:11} implies that $\Pb(\norm{\hat\beta_2-\beta_0}_2\geq C\sqrt{(d+\log(en))/n})\leq 1/n$. (See also, Theorem~4.2 of \cite{oliveira2016lower} or Proposition~26 of \cite{kuchibhotla2020berry} for OLSE results). Hence, \eqref{eq:thm:7:1}, \eqref{eq:thm:7:3}, and Theorem~\ref{thm:12} directly follow from the application of Theorem~\ref{thm:1}--\ref{thm:3} through verifications \ref{_toglm:1}--\ref{_toglm:5} of assumptions. Hence, we prove \eqref{eq:thm:7:2} only. We may assume $d\log^2(en)\leq n$, otherwise, the bounds in Theorem~\ref{thm:12} hold immediately. We use \eqref{eq:pf_of27.1} for precise control; in particular, for ${\rm a}\in\set{\rm SN, SSW}$, $$\sup_{\alpha\in(0,1)}[\Pb(\theta_0\notin\widehat{\mathrm{CI}}_{\alpha}^{\rm a})-\alpha]_+\lesssim \log^2(en)\sqrt{\frac{d}{n}} + p_\lambda +\left(\frac{\log(en)\log^2(e d) d^{2/q}}{\lambda n^{1-2/q}}\right)^{\frac{q}{q+2}},$$ for any $\lambda>0$ where $q=\min\set{q_x/2,q_{cy}}$ and $p_\lambda = \Pb(\lambda_{\rm min}({\rm Corr}(\hat\Sigma_2^{-1}V\hat \Sigma_{\Dc_2}^{-1}))\leq \lambda)$. Note that $\lambda_{\rm min}({\rm Corr}(\hat\Sigma_2^{-1}V\hat\Sigma_2^{-1}))\geq \underline{\lambda}_V/\overline{\lambda}_V(\lambda_{\rm min}(\hat\Sigma_2)/\lambda_{\rm max}(\hat\Sigma_2))^2$. We take $$\lambda^{-1} \asymp 1 + \frac{d\log(ed\delta)}{n} + \frac{d^2\log^3(ed/\delta)}{\delta^{4/q_x}n^{2-4/q_x}},$$ which leads to $p_\lambda\lesssim\delta$ from Proposition~\ref{prop:C4}. As long as $\delta\leq {\rm poly}(n)$ and $q\geq 6$, $$p_\lambda+\frac{\log^{3/4}(en)\log^{3/2}(ed)d^{1/4}}{\lambda^{3/4} n^{1/2}}\lesssim \delta + \frac{\log^{3/4}(en)\log^{15/4}(ed)d^{7/4}}{\delta^{3/q_x}n^{2-3/q_x}}.$$ We take the following $\delta$, which can be bounded as long as $q_x\geq 12$ and proves the theorem.
\begin{align*}
    &\delta = \left(\frac{\log^{3/4}(en)\log^{15/4}(ed)d^{7/4}}{n^{2-3/q_x}}\right)^{q_x/(q_x+3)}\lesssim \frac{\log^{3/5}(en)\log^{3}(ed)d^{7/5}}{n^{7/5}}\lesssim\log^2(en)\sqrt{\frac{d}{n}}.
\end{align*} 

\subsection{Proof of Theorem~\ref{thm:8}.} We first control the widths of the SeN and SpN confidence sets. We define some events: \begin{align*}
    \Sc_1&=\Set{\lambda_{\rm min}(\Sigma^{-\frac{1}{2}}\hat\Sigma_i\Sigma^{-\frac{1}{2}})\geq 1/2;~i=1,2}\\
    \Sc_2 &= \Set{\frac{K_n}{\sqrt{n_1}}\max_{i\in\Ic_1}\norm{X_i}_{\Sigma^{-1}}\leq 1\wedge \frac{\underline{\lambda}_\Sigma^{1/2}{\underline\sigma}_a}{4\overline{\lambda}_\Sigma^{1/2}{\overline\sigma_a}}}\\
    \Sc_3 &= \Set{\lambda_{\rm max}(\Sigma^{-\frac{1}{2}}\hat\Sigma_1\Sigma^{-\frac{1}{2}})\leq cK_x^2}~\mbox{for a large universal constant}~c>0.\\
    \Sc_4&=\Set{\frac{K_n^2}{n_1^2}\sum_{i\in\Ic_1}\norm{\Sigma^{-1/2}X_i}_2^2(Y_i-X_i^\top\beta_0)^2\leq 2K_x^2K_y^2\frac{K_n^2d}{n}},\\
    \Sc_5&=\Set{\norm{\hat\beta_i-\beta_0}_{\hat\Sigma_i}^2\leq\left(\frac{4\overline{\lambda}_V}{\underline{\lambda}_\Sigma}+ C_{q_{xy}}K_x^2K_y^2\right)\left(\frac{d+2\log(2n)}{n}\right);i=1,2}~\mbox{for some } C_{q_{xy}}>0.
\end{align*}Throughout, the constant $C$, independent of $n$, $\overline{\sigma}$, and $\underline{\sigma}$, may differ in line by line.

Suppose that $\Sc_1$---$\Sc_5$ hold true. From the definition of $\hat\beta_1$, the confidence set $\widehat{\mathrm{CI}}^{\rm SeN}(K_n)$ can be expressed as
\begin{equation*}
    \widehat{\mathrm{CI}}^{\rm SeN}(K_n)=\Set{\theta: \max_{1\leq j\leq d}\Abs{\frac{\ev_j^\top\widetilde Q_2\hat\Sigma_1(\hat\beta_1-\theta)}{\sqrt{n_1^{-1}\sum_{i\in\Ic_1}(\ev_j^\top\widetilde Q_2X_i)^2(Y_i-X_i^\top\theta)^2}}}\leq \frac{K_n}{\sqrt{n_1}}\frac{1}{\sqrt{1+K_n^2/n_1}}}.
\end{equation*} Therefore, every $\theta\in\widehat{\mathrm{CI}}^{\rm SeN}(K_n)$ satisfies that
\begin{align}\label{eq:thm.10:1}
    &\underline{\sigma}_a^{2}\lambda_{\rm min}(\hat\Sigma_1)\norm{\hat\Sigma_1^{1/2}(\hat\beta_1-\theta)}_2^2\leq \norm{\widetilde Q_2\hat\Sigma_1(\hat\beta_1-\theta)}_2^2\nonumber\\
    &\quad\leq\frac{K_n^2}{n_1}\sum_{j=1}^d\frac{1}{n_1}\sum_{i\in\Ic_1}(\ev_j^\top\widetilde Q_2X_i)^2(Y_i-X_i^\top\theta)^2=\frac{K_n^2}{n_1^2}\sum_{i\in\Ic_1}\norm{\widetilde Q_2X_i}_2^2(Y_i-X_i^\top\theta)^2\nonumber\\
    &\quad \leq \lambda_{\rm max}(\Sigma^{1/2}\widetilde Q_2^\top\widetilde Q_2\Sigma^{1/2})\frac{K_n^2}{n_1^2}\sum_{i\in\Ic_1}\norm{\Sigma^{-1/2}X_i}_2^2(Y_i-X_i^\top\theta)^2,
\end{align}for all $\theta\in\widehat{\mathrm{CI}}^{\rm SeN}(K_n)$. This implies that
\begin{align}\label{eq:thm.10:2}
    &\underline{\sigma}_a\lambda_{\rm min}^{1/2}(\hat\Sigma_1)\norm{\hat\Sigma_1^{1/2}(\hat\beta_1-\theta)}_2\leq\\
    &\overline{\lambda}^{1/2}_\Sigma\overline{\sigma}_a\frac{K_n}{n_1}\left\{\left(\sum_{i\in\Ic_1}\norm{\Sigma^{-1/2}X_i}_2^2(Y_i-X_i^\top\hat\beta_1)^2\right)^{1/2}+\left(\sum_{i\in\Ic_1}\norm{\Sigma^{-1/2}X_i}_2^2(X_i^\top(\hat\beta_1-\theta))^2\right)^{1/2}\right\}\nonumber.
\end{align} To analyze the last display, we note that
\begin{eqnarray}\label{eq:thm.10:3}
    \sum_{i\in\Ic_1}\norm{\Sigma^{-1/2}X_i}_2^2(X_i^\top(\hat\beta_1-\theta))^2 &\leq& \max_{i\in\Ic_1}\norm{\Sigma^{-1/2}X_i}_2^2\sum_{i\in\Ic_1}(X_i^\top(\hat\beta_1-\theta))^2 \nonumber\\
    &=& n_1\norm{\hat\Sigma_1^{1/2}(\hat\beta_1-\theta)}_2^2\max_{i\in\Ic_1}\norm{\Sigma^{-1/2}X_i}_2^2.
\end{eqnarray} On $\Sc_1\cap\Sc_2$, we can deduce from \eqref{eq:thm.10:2} and \eqref{eq:thm.10:3} that
\begin{align*}
    \frac{1}{4}\underline{\sigma}_a\underline{\lambda}_\Sigma^{1/2}\norm{\hat\Sigma^{1/2}(\hat\beta_1-\theta)}_2\leq \overline{\lambda}^{1/2}_\Sigma\overline{\sigma}_a\frac{K_n}{n_1}\left(\sum_{i\in\Ic_1}\norm{\Sigma^{-1/2}X_i}_2^2(Y_i-X_i^\top\hat\beta_1)^2\right)^{1/2},
\end{align*} which implies that
\begin{equation*}
    \norm{\hat\beta_1-\theta}_2\leq \frac{8\overline{\lambda}_\Sigma^{1/2}{\overline{\sigma}_a}}{\underline{\lambda}_\Sigma^{1/2}{\underline{\sigma}_a}}\frac{K_n}{n_1}\left(\sum_{i\in\Ic_1}\norm{\Sigma^{-1/2}X_i}_2^2(Y_i-X_i^\top\hat\beta_1)^2\right)^{1/2},
\end{equation*} for all $\theta\in \widehat{\mathrm{CI}}^{\rm SeN}(K_n)$. Note that the right-hand side no longer depends on $\theta$, and therefore, taking supremum over $\theta\in \widehat{\mathrm{CI}}^{\rm SeN}(K_n)$, we get from triangle inequality that
\begin{equation}\label{eq:thm.10:4}
    {\rm diam}\left(\widehat{\mathrm{CI}}^{\rm SeN}(K_n)\right)\leq \frac{16\overline{\lambda}_\Sigma^{1/2}{\overline{\sigma}_a}}{\underline{\lambda}_\Sigma^{1/2}{\underline{\sigma}_a}}\frac{K_n}{n_1}\left(\sum_{i\in\Ic_1}\norm{\Sigma^{-1/2}X_i}_2^2(Y_i-X_i^\top\hat\beta_1)^2\right)^{1/2},
\end{equation} as long as the event $\Sc_1\cap\Sc_2$ holds true. We now analyze the quantity in \eqref{eq:thm.10:4}. Observe that
\begin{align}
    \left(\frac{K_n^2}{n_1^2}\sum_{i\in\Ic_1}\norm{\Sigma^{-1/2}X_i}_2^2(Y_i-X_i^\top\hat\beta_1)^2\right)^{1/2}\leq \Av_1^{1/2}+\Av_2^{1/2}~\mbox{ where }\\
    \Av_1 = \frac{K_n^2}{n_1^2}\sum_{i\in\Ic_1}\norm{\Sigma^{-1/2}X_i}_2^2(Y_i-X_i^\top\beta_0)^2,~\Av_2 = \frac{K_n^2}{n_1^2}\sum_{i\in\Ic_1}\norm{\Sigma^{-1/2}X_i}_2^2(X_i^\top(\beta_0-\hat\beta_1))^2.\nonumber
\end{align}
On $\Sc_2\cap\Sc_5$, one has
\begin{eqnarray*}
    \Av_2&\leq& \frac{K_n^2}{n_1^2}\left(\max_{i\in\Ic_1}\norm{\Sigma^{-1/2}X_i}_2^2\right)\sum_{i\in\Ic_1}(X_i^\top(\beta_0-\hat\beta_1))^2\\
    &\leq&\frac{\underline{\lambda}_\Sigma{\underline{\sigma}_a}^2}{8\overline{\lambda}_\Sigma{\overline{\sigma}_a}^2}\frac{1}{n_1}\sum_{i\in\Ic_1}(X_i^\top(\beta_0-\hat\beta_1))^2\leq \frac{\underline{\lambda}_\Sigma{\underline{\sigma}_a}^2}{8\overline{\lambda}_\Sigma{\overline{\sigma}_a}^2}(\beta_0-\hat\beta_1)^\top\hat\Sigma_1(\beta_0-\hat\beta_1)\\
    &\leq&\left(\frac{4\overline{\lambda}_V}{\underline{\lambda}_\Sigma}+ C_{q_{xy}}K_x^2K_y^2\right)\frac{\underline{\lambda}_\Sigma{\underline{\sigma}_a}^2}{8\overline{\lambda}_\Sigma{\overline{\sigma}_a}^2}\left(\frac{d+2\log(2n)}{n}\right).
\end{eqnarray*}
Combining all, on $\Sc_{1}\cap\Sc_{2}\cap\Sc_{4}\cap \Sc_{5}$, one has
\begin{equation*}
    {\rm diam}\left(\widehat{\mathrm{CI}}^{\rm SeN}(K_n)\right)\lesssim \frac{{\overline{\sigma}_a}}{{\underline{\sigma}_a}}\left(\Av_1^{1/2}+\Av_2^{1/2}\right)\lesssim \frac{{\overline{\sigma}_a}}{{\underline{\sigma}_a}}\left(\frac{K_n^2d+\log(en)}{n}\right)^{1/2}.
\end{equation*}

The confidence set $\widehat{\mathrm{CI}}^{\rm SpN}(K_n)$ can be similarly analyzed. Following the derivation of \eqref{eq:thm.10:1}, we get that every $\theta\in\widehat{\mathrm{CI}}^{\rm SpN}(K_n)$ satisfies that 
\begin{equation*}
    \underline{\sigma}_a^{2}\lambda^2_{\rm min}(\hat\Sigma_1)\norm{\hat\beta_1-\theta}_2^2\leq \overline\lambda_\Sigma \overline\sigma^2\frac{K_n^2}{n_1^2}\sum_{i\in\Ic_1}\norm{\Sigma^{-1/2}X_i}_2^2(Y_i-X_i^\top\hat\beta_2)^2.
\end{equation*} Moreover, we note that $(\sum_{i\in\Ic_1}\norm{\Sigma^{-1/2}X_i}_2^2(Y_i-X_i^\top\hat\beta_2)^2)^{1/2}\leq(\sum_{i\in\Ic_1}\norm{\Sigma^{-1/2}X_i}_2^2(Y_i-X_i^\top\beta_0)^2)^{1/2}+(\sum_{i\in\Ic_1}\norm{\Sigma^{-1/2}X_i}_2^4\norm{\hat\beta_2-\beta_0}_\Sigma^2)^{1/2}\leq(\sum_{i\in\Ic_1}\norm{\Sigma^{-1/2}X_i}_2^2(Y_i-X_i^\top\beta_0)^2)^{1/2}+n_1^{1/2}\max_{i\in\Ic_1}\norm{\Sigma^{-1/2}X_i}_2^2\norm{\hat\beta_2-\beta_0}_\Sigma.$ Therefore, we get $\Sc_{1}\cap\Sc_{2}\cap\Sc_{4}\cap \Sc_{5}$, 
\begin{align*}
    \norm{\theta-\hat\beta_1}_2\lesssim\frac{\overline{\sigma}}{\underline{\sigma}}\left[\norm{\hat\beta_2-\beta_0}_2 + \sqrt{\frac{K_n^2d+\log(en)}{n}}\right],
\end{align*} uniformly over $\theta\in\widehat{\mathrm{CI}}^{\rm SpN}(K_n)$.

Finally, we analyze the GaN confidence set. Since $\widetilde Q_2$ is positive definite with \ref{asmp:S}, on $\Sc_1\cap\Sc_2\cap\Sc_3$, we have for $\widehat{\mathrm{CI}}^{\rm GaN}(K_n)$,
$$\abs{(\theta-\hat\beta_2)^\top\widetilde Q_2\hat\Sigma_1(\hat\beta_2-\theta)}\geq C\underline{\sigma}\norm{\theta-\hat\beta_2}_\Sigma^2,$$ and $$\abs{(\theta-\hat\beta_2)^\top\widetilde Q_2\hat\Sigma_1(\hat\beta_1-\hat\beta_2)}\leq C\overline{\sigma}\norm{\theta-\hat\beta_2}_\Sigma\norm{\hat\beta_1-\hat\beta_2}_\Sigma.$$ Also, $\norm{\Sigma^{1/2}\widetilde Q_2^\top(\theta-\hat\beta_2)}_2\leq C\overline{\sigma}\norm{\theta-\hat\beta_2}_\Sigma$, and on $\Sc_2$,  $K_n^2n_1^{-1}\sum_{i\in\Ic_1}\norm{\Sigma^{-\frac{1}{2}}X_i}^2(Y_i-X_i^\top\hat\beta_2)^2\leq n_1^{-1}\sum_{i\in\Ic_1}(Y_i-X_i^\top\hat\beta_2)^2$. Hence, for any $\theta\in\widehat{\mathrm{CI}}^{\rm GaN}(K_n)$, on $\Sc_1\cap\Sc_2\cap\Sc_3\cap\Sc_5$,
\begin{align*}
    \underline{\sigma}\sqrt{n}\norm{\theta-\hat\beta_2}_\Sigma\leq C\overline{\sigma}\left[\sqrt{n}\norm{\hat\beta_1-\hat\beta_2}_\Sigma+ \Set{n_1^{-1}\sum_{i\in\Ic_1}(Y_i-X_i^\top\hat\beta_2)^2}^{1/2}\right]\\
    \lesssim \overline{\sigma}\left[\sqrt{n}\norm{\hat\beta_1-\hat\beta_2}_\Sigma+ \Set{n_1^{-1}\sum_{i\in\Ic_1}(Y_i-X_i^\top\beta_0)^2}^{1/2}+\norm{\hat\beta_2-\beta_0}_{\hat\Sigma_1}\right]\\
    \lesssim \overline{\sigma}\left[\sqrt{d+\log(en)}+ \Set{n_1^{-1}\sum_{i\in\Ic_1}(Y_i-X_i^\top\beta_0)^2}^{1/2}\right]
\end{align*} From Chebyshev inequality, one further gets $\Pb(n_1^{-1}\sum_{i\in\Ic_1}(Y_i-X_i^\top\beta_0)^2\geq 2K_y^2)\leq 1/n$.

Now, we analyze the probabilities of events. Provided that $d+2\log(4n)\leq(n/18K_x^2)^2$, an application of Theorem~4.1 of \cite{oliveira2016lower} results in $\Pb(\Sc_1)\geq1-1/n$. Furthermore, we note from Markov's inequality that for $A>0$,
\begin{align}\label{eq:prob_control_max_norm}
    &\Pb\left(\frac{K_n^2}{n_1}\max_{i\in\Ic_1}\norm{\Sigma^{-1/2}X_i}_2^2>A\right)\leq \frac{K_n^{q_x}}{A^{q_x/2}n_1^{q_x/2}}\Eb\left[\max_{i\in\Ic_1}\norm{\Sigma^{-1/2}X_i}_2^{q_x}\right]\nonumber\\
    &\quad\leq\frac{K_n^{q_x}}{A^{q_x/2}n_1^{q_x/2}}\Eb\left[\sum_{i\in\Ic_1}\norm{\Sigma^{-1/2}X_i}_2^{q_x}\right]\leq\frac{K_n^{q_x}}{A^{q_x/2}n^{q_x/2}}\left(K_x^{q_x}d^{q_x/2}n\right)=\left(\frac{K_xK_n^2d}{An^{1-2/q_x}}\right)^{q_x/2}.
\end{align} Proposition~\ref{prop:C4} (reproduced in Proposition~\ref{prop:C4}) implies $\Pb(\Sc_3)\geq 1-(d\log^{3/2}(en)/n^{1-2/q_x})^{q_x/2}$. Finally, Proposition~\ref{prop:C5} controls the probabilities of $\Sc_4$ and $\Sc_5$. This completes the proof.

\subsection{Proof of Theorem~\ref{thm:9}.}
The miscoverage rates of $\hat{\mathscr{C}_\alpha^{\rm a}}$ for ${\rm a}\in\{\rm SeN, SpN\}$ are provided in the proof of Theorem~\ref{thm:6}. To prove the validity of SeN with Bonferroni correction, we apply Theorem~D.1 of \cite{Chernozhukov2019inference} with union bounds which gives $\Pb(T^{\rm SeN}(\theta,\widetilde Q_2)\leq z_{1/(2nd)}|\Dc_2)\leq \alpha/n +CL_{q\wedge 3}^{q\wedge 3}(1+z_{1/(2nd)})^{q\wedge 3}/n^{q\wedge 3}\leq 1/n+CL_{q\wedge 3}^{q\wedge 3}\log^{{q\wedge 3}/2}(end)/n^{q\wedge 3},$ for some universal constant $c>0$. Algebraic simplifications leads to the desired results.

\subsection{Proof of Theorems~\ref{thm:A1} and \ref{thm:10}.} We first prove Thereom~\ref{thm:A1} for $\widehat{\rm CI}_\alpha^{\rm SeN}$ whereas the results for SpN follow similarly. Denote $\widehat{\rm CI}_\alpha^{\rm SeN}=\widehat{\rm CI}_\alpha$ and let $\hat b_n = \norm{(J\widetilde Q_2^\top-I_d)\xi}_{\rm op}$, $\hat A_n(r)=\sup_{\theta\in\Bb(\theta_0,r)}\norm{\xi\widetilde Q_2(\hat J_\theta-J)\widetilde Q_2^\top\xi}_{\rm op}$ where $\hat J_\theta = \int_0^1\nabla\bar\psi_1(t\theta+(1-t)\hat\theta_1)\,dt$, and $\hat C_n(r)=\sup_{\theta\in\Bb(\theta_0,r)}\norm{\hat V_1(\theta)-V(\theta)}_{\rm op}$. Define $s_n^{\rm W} = z_{\alpha/2}\sqrt{\xi^\top J^{-1}VJ^{-1}/n}$, so that $\widehat{\rm CI}_\alpha^{\rm Wald}=[\xi^\top\hat\theta_1\pm s_n^{\rm W}]$, and define $\hat{\rm CI}_1 = [\xi^\top\hat\theta_1\pm \set{s_n^{\rm W}+\hat b_ns_n+\overline{\sigma}s_n\hat A_n(s_n)+ z_{\alpha/2}/\sqrt{n}(\overline{\sigma}_J\hat b_n+\phi(s_n)+\hat C_n(s_n))^{1/2}}].$ Suppose that $\widehat{\rm CI}_\alpha\subseteq\widehat{\rm CI}_1$, then ${\rm d}_{\rm H}(\widehat{\rm CI}_\alpha\to \widehat{\rm CI}_\alpha^{\rm Wald})\leq{\rm d}_{\rm H}(\widehat{\rm CI}_1\to \widehat{\rm CI}_\alpha^{\rm Wald})$, which proves Thereom~\ref{thm:A1}. Now, Taylor expansion of $\bar\psi_1(\theta)$ at $\hat\theta_1$ yields $\bar\psi_1(\theta)=\hat J_\theta(\theta-\hat\theta_1)$, and note that $$\xi^\top \widetilde Q_2\bar\psi_1(\theta) = \xi^\top(\theta-\hat\theta_1) + \xi^\top(\widetilde Q_2J-I_d)(\theta-\hat\theta_1)+\xi^\top\widetilde Q_2(\hat J_\theta-J)(\theta-\hat\theta_1).$$ Regarding the variance, we note that
$\norm{\widetilde Q_2^\top\xi}_{\hat V_1(\theta)}=\norm{J^{-1} \xi}_{V(\theta_0)} + (\norm{\widetilde Q_2^\top\xi}_{V(\theta_0)}-\norm{J^{-1} \xi}_{V(\theta_0)})+(\norm{\widetilde Q_2^\top\xi}_{V(\theta)}-\norm{\widetilde Q_2^\top\xi}_{V(\theta_0)})+(\norm{\widetilde Q_2^\top\xi}_{\hat V_1(\theta)}-\norm{\widetilde Q_2^\top\xi}_{V(\theta)})$. Therefore, $\widehat{\rm CI}_\alpha\subseteq\widehat{\rm CI}_1$. For SpN, one can prove $\widehat{\rm CI}_\alpha^{\rm SpN}\subseteq\widehat{\rm CI}_2$ where $\widehat{\rm CI}_3=[\xi^\top\hat\theta_2\pm \set{s_n^{\rm W}+\hat b_ns_n+\overline{\sigma}s_n\hat A_n(s_n)+ z_{\alpha/2}/\sqrt{n}(\overline{\sigma}_J\hat b_n+\phi(r_{n})+|\norm{\widetilde Q_2^\top\xi}_{\hat V_1(\hat \theta_2)}-\norm{\widetilde Q_2^\top\xi}_{V(\hat \theta_2)}|)^{1/2}}]$.

Now we prove Theorem~\ref{thm:10}. Write $\hat a = \hat\Sigma_2^{-1}\xi$. In the context of linear regression, $\hat b_n = \norm{(\Sigma\hat\Sigma_2^{-1}-I_d)\xi}_2$, $\hat A_n(r) \equiv |\hat a^\top(\hat\Sigma_1-\Sigma)\hat a^\top|$ for all $r>0$. In event $\set{\norm{\Sigma^{-1/2}\hat\Sigma_i\hat\Sigma^{-1/2}-I_d}_{\rm op}\leq C\sqrt{d/n}, i=1,2}$, we have $\hat b_n,\hat A_n(r)\lesssim \sqrt{d/n}$ where this event happened with probability at least $1-C((d\log^3(en)/n^{1-4/q_x})^{q_x/4})$ by Proposition~\ref{prop:C4}. Also, according to Theorem~\ref{thm:10}, we have $s_n\lesssim\sqrt{d/n}$. Moreover, $\hat a^\top\hat V_1(\theta)\hat a=\hat a^\top\hat V_1(\theta_0)\hat a +2 \hat M_1(\theta-\beta_0) + (\theta-\beta_0)^\top \hat M_2(\theta-\beta_0)$ where $\hat M_1 = n_1^{-1}\sum_{i\in\Ic_1}(\hat a^\top X_i)^2(Y_i-X_i^\top\beta_0)X_i$, $\hat M_2 = n_1^{-1}\sum_{i\in\Ic_1}(\hat a^\top X_i)^2X_iX_i^\top$. Due to Chebyshev inequality, $\Pb(\hat a^\top\hat V_1(\theta_0)\hat a-\hat a^\top V\hat a\geq t|\Dc_2)\leq (nt^2)^{-1}\norm{\hat a}_\Sigma^2K_x^4K_y^4$ and $\Pb(\norm{\hat M_1-\Eb[\hat M_1|\Dc_2]}_\Sigma\geq t)\leq (nt^2/d)^{-1}\norm{\hat a}_\Sigma^2K_x^6K_y^2$. Proposition~\ref{prop:C4} implies that $$\Pb(\norm{\hat M_2-\Eb[\hat M_2|\Dc_2]}_{\rm op}\geq \sqrt{d\log(ed/\delta)/n} + d\log^{3/2}(ed/\delta)/(\delta^{4/q_x}n^{1-4/q_x}))\leq\delta.$$ Therefore, $$\Pb(\hat C_n(r)\geq 1/\sqrt{n}+r+r^2)\leq C((d\log^{3/2}(en)/n^{1-4/q_x})^{q_x/4}).$$ Now an application of Theorem~\ref{thm:A1} proves the concentration of SeN. The results for SpN follows similarly.
\subsection{Proof of Theorem~\ref{thm:11}.}
As per \ref{_toglm:4}, we can choose $\underline\mu(r_n)=C_1 (r_n\wedge \delta_0)$, and \ref{_toglm:5} implies that there exists $c_0>0$ such that if $d/n\leq c_0$, then $u_{1,n}(r_n,\varepsilon )/\sqrt{n}\leq C_0C_1r_n/2.$ Set $$C^{-1}r_n = \sqrt{\frac{d+\log(2/\varepsilon)}{n}}+\frac{d^{1/2}}{n^{1-1/q_{xy}}\varepsilon^{1/q_{xy}}},$$ where $C=2/(C_0C_1)$. As long as the RHS is no larger than $\min\set{\sqrt{c_0}, C_0C_1\delta_0/2}$, $r_n$ satisfies \eqref{eq:thm:4:1}, and thus, Theorem~\ref{thm:4} applies and yields the result.
\subsection{Proof of Theorem~\ref{thm:13}.}
Let $r_{n}$ be the convergence rate of $\hat\theta_2$ in Theorem~\ref{thm:11}, satisfying $\Pb(\norm{\hat\theta_2-\theta_0}_2\geq r_n)\leq\varepsilon$. Then, $\phi(r_n)$ and $\gamma_q(r_n)$ are bounded by constant due to \ref{_toglm:2} and \ref{_toglm:3}, respectively. From \ref{_toglm:4}, we can choose  $\underline\mu(s_n)=C_1 (s_n\wedge \delta_0)$, and \ref{_toglm:5} implies that $u_{1,n}(r_{n},\varepsilon )/\sqrt{n}\leq C_0C_1s_n/2$, provided that $K_n^2d/n\leq c$ for a sufficiently small $c_0>0$. Set $$
   C^{-1}s_n = \sqrt{\frac{K_n^2d+\log(2/\varepsilon)}{n}}+\frac{d^{1/2}}{n^{1-1/q_{xy}}\varepsilon^{1/q_{xy}}},$$ where $C=3(C_0C_1)^{-1}$. As long as the RHS is no larger than $c = \min\set{\sqrt{c_0}, C_0C_1\delta_0/2}$, $s_n$ satisfies \eqref{eq:thm:5:1}, and Theorem~\ref{thm:5} yields the result.

\section{Useful Lemmas and Propositions}
\begin{lemma}\label{lem:C1}
    Let $X_1,\ldots,X_N$ be not necessarily independent centered Gaussian random variables. Let $\sigma_{\rm max}^2=\max_{i\in[N]}{\rm Var}(X_i)$, then
    \begin{equation*}
        \Pb\left(\max_{i\in[N]}\abs{X_i}>t\right)\leq Ne^{-t^2/(2\sigma_{\rm max}^2)},
    \end{equation*} for all $t>0$.
\end{lemma}
\begin{proof}[proof of Lemma~\ref{lem:C1}] Let $\sigma_i^2={\rm Var}(X_i)$. It follows from the union bound and that $1-\Phi(t/\sigma_i)\leq e^{-t^2/(2\sigma_i^2)}/2$ for $t>0$ that
\begin{eqnarray*}
    \Pb\left(\max_{i\in[n]}\abs{X_i}>t\right)\leq \sum_{i=1}^N\Pb\left(\abs{X_i}>t\right)=2\sum_{i=1}^N\left(1-\Phi(t/\sigma_i)\right)\leq \sum_{i=1}^Ne^{-t^2/(2\sigma_i^2)}\leq Ne^{-t^2/(2\sigma_{\rm max}^2)}.
\end{eqnarray*}
    
\end{proof}
\begin{lemma}\label{lem:C2} Let $X_1,\ldots,X_n$ be independent non-negative random variables and let $\Eb[X_i^{1+\varepsilon}]<\infty$, $i=1,\ldots,n$, for some $\varepsilon>0$. Let $\mu_n=\sum_{i=1}^n\Eb[X_i]$ and $B_{n,\varepsilon}=(\sum_{i=1}^n\Eb[X_i^{1+\varepsilon}])^{1/(1+\varepsilon)}$. Then,
for $0\leq x\leq\mu_n$,
\begin{equation}\label{eq:lem:C2:1}
    \Pb\left(\sum_{i=1}^nX_i\leq x\right)\leq\exp\left(-\frac{\varepsilon}{1+\varepsilon}\left(\frac{\mu_n-x}{B_{n,\varepsilon}}\right)^{1+1/\varepsilon}\right).
\end{equation}
\end{lemma}
\begin{proof}[proof of Lemma~\ref{lem:C2}] We follow the proof of lemma 2.19 of \cite{pena2009self}. Note that $e^{-a}\leq 1-a+(1+\varepsilon)^{-1}a^{1+\varepsilon}$ for all $a>0$. For any $t\geq 0$ and $0\leq x\leq \mu_n$, we have from Chernoff inequality that
    \begin{eqnarray*}
        \Pb\left(\sum_{i=1}^nX_i\leq x\right)&\leq& e^{tx}\prod_{i=1}^n\Eb\left[\exp\left(-tX_i\right)\right]\\
        &\leq& e^{tx}\prod_{i=1}^n\Eb\left[1-tX_i+\frac{(tX_i)^{1+\varepsilon}}{1+\varepsilon}\right]\\
        &\leq&\exp\left(t(x-\mu_n)+\frac{(B_{n,\varepsilon}t)^{1+\varepsilon}}{1+\varepsilon}\right).
    \end{eqnarray*} Taking $t=((\mu_n-x)/B_{n,\varepsilon}^{1+\varepsilon})^{1/\varepsilon}$ yields \eqref{eq:lem:C2:1}.
\end{proof}
\begin{lemma}\label{lem:C3}
    Let $X_1,\ldots,X_n$ be independent random vectors in $\Real^p$ where $X_i=(X_{i1},\ldots,X_{ip})^\top$ and $X_{ij}\geq0$ for all $i\in\set{1,\ldots,n}$ and $j\in\set{1,\ldots,p}$. For $\varepsilon>0$, define
    \begin{equation*}
        L_{1+\varepsilon} = \max_{1\leq j\leq p}\frac{(n^{-1}\sum_{i=1}^n\Eb[X_{ij}^{1+\varepsilon}])^{1/(1+\varepsilon)}}{(n^{-1}\sum_{i=1}^n\Eb[X_{ij}])}.
    \end{equation*} For $c\in(0,1)$, one has
    \begin{equation*}
        \Pb\left(\min_{1\leq j\leq p}\frac{\sum_{i=1}^nX_{ij}}{\sum_{i=1}^n\Eb[X_{ij}]}\geq 1-c\right)\geq 1-p\exp\left(-\frac{n\varepsilon}{1+\varepsilon}\left(\frac{c}{L_{1+\varepsilon}}\right)^{1+1/\varepsilon}\right).
    \end{equation*} Hence, if $c_n:=L_{1+\varepsilon}((1+1/\varepsilon)\log(en)/n)^{\varepsilon/(1+\varepsilon)}\leq1$, then one has
    \begin{equation*}
         \Pb\left(\min_{1\leq j\leq p}\frac{\sum_{i=1}^nX_{ij}}{\sum_{i=1}^n\Eb[X_{ij}]}\geq 1-c_n\right)\geq 1-\frac{p}{n}.
    \end{equation*}
\end{lemma}
\begin{proof}[proof of Lemma~\ref{lem:C3}]
    For any $c\in(0,1)$, Lemma~\ref{lem:C2} implies that 
    \begin{eqnarray*}
        \Pb\left(\sum_{i=1}^nX_{ij}\leq(1-c)\sum_{i=1}^n\Eb[X_{ij}]\right)&\leq&\exp\left(-\frac{\varepsilon c^{1+1/\varepsilon}(\sum_{i=1}^n\Eb[X_{ij}])^{1+1/\varepsilon}}{(1+\varepsilon)(\sum_{i=1}^n\Eb[X_{ij}^{1+\varepsilon}])^{1/\varepsilon}}\right)\\
        &\leq&\exp\left(-\frac{n\varepsilon}{1+\varepsilon}\left(\frac{c}{L_{1+\varepsilon}}\right)^{1+1/\varepsilon}\right),
    \end{eqnarray*} for all $j=1,\ldots,p$. Hence, we get
    \begin{eqnarray*}
        \Pb\left(\min_{1\leq j\leq p}\frac{\sum_{i=1}^nX_{ij}}{\sum_{i=1}^n\Eb[X_{ij}]}\geq 1-c\right)&=&1-\Pb\left(\min_{1\leq j\leq p}\frac{\sum_{i=1}^nX_{ij}}{\sum_{i=1}^n\Eb[X_{ij}]}< 1-c\right)\\
        &\geq&1-\sum_{j=1}^p\Pb\left(\frac{\sum_{i=1}^nX_{ij}}{\sum_{i=1}^n\Eb[X_{ij}]}< 1-c\right)\\
        &\geq&1-p\exp\left(-\frac{n\varepsilon}{1+\varepsilon}\left(\frac{c}{L_{1+\varepsilon}}\right)^{1+1/\varepsilon}\right).
    \end{eqnarray*} 
\end{proof}

\begin{proposition}\label{prop:C1}
    Let $X_1,\ldots,X_n\in\Real^d$ be i.i.d with $\Eb[X]=0_d$ and $\Eb[X_1X_1^\top]=\Sigma$. Then,
    \begin{align*}
        \Pb\left(\Norm{\bar X_n}_{\Sigma^{-1}}\geq \bigg[\sqrt{\frac{2d+3\log(2/\delta)}{n}}\mathbbm{1}_{q>2}+C_q\frac{\Eb^{1/q}[\norm{X}^q_{\Sigma^{-1}}]}{\delta^{1/q}n^{1-1/q}}\bigg]\right)\leq\delta,~\mbox{for}~q\geq2\\
        \Pb\left(\frac{1}{n}\sum_{i\in[n]]}\norm{X_i}^2_{\Sigma^{-1}}\geq d + \bigg[\sqrt{\frac{2\Eb[\norm{X}^4_{\Sigma^{-1}}]\log(1/\delta)}{n}}\mathbbm{1}_{q>4}+C_q\frac{\Eb^{2/q}[\norm{X}^q_{\Sigma^{-1}}]}{\delta^{2/q}n^{1-2/q}}\bigg]\right)\leq \delta~\mbox{for}~q\geq 4.
    \end{align*}
\end{proposition}
Moment conditions on norms are generally weaker than those on projections. For e.g., one can show that $\Eb^{1/q}[\norm{X}^q_{\Sigma^{-1}}]\leq \sqrt{d}\sup_{j\in[d]}\Eb^{1/q}[|\ev_j^\top\Sigma^{-1/2}X|^q]$ while the reverse does not hold in general. See \cite{mourtada2022exact} for related discussions on norm kurtosis assumption.
\begin{proof} The $q=2$ case follows from Chebyshev inequality. For $q>2$, Theorem~4 of \cite{einmahl2008characterization} applies with $\delta=\eta=1$ and yields $$\Pb(\norm{\bar X_n}_{\Sigma^{-1}}\geq 2\Eb[\norm{\bar X_n}_{\Sigma^{-1}}^2]+t)\leq e^{-\frac{t^2}{3\Lambda_n}}+C_q\frac{\Eb\norm{X_1}_{\Sigma^{-1}}^q}{t^qn^{q-1}},$$ for any $q>2$ where $\Lambda_n = n^{-1}\lambda_{\rm max}(\Eb[\Sigma^{-1/2}X_1X_1^\top\Sigma^{-1/2}])=1/n$. Moreover, $\Eb[\norm{\bar X_n}_{\Sigma^{-1}}]\leq \Eb^{1/2}[\norm{\bar X_n}_{\Sigma^{-1}}^2]=\sqrt{d/n}$. Hence, taking the following $t=\sqrt{\frac{3\log(2/\delta)}{n}} + (2C_q)^{1/q}\frac{\Eb^{1/q}[\norm{X}^q_{\Sigma^{-1}}]}{\delta^{1/q}n^{1-1/q}}$ proves the first part of the proposition. Now, let $Y_i = \norm{X_i}_{\Sigma^{-1}}^2-1$, and apply Equation (1.9) of \cite{rio2017constants} to get, for $q>4$, $$\Pb(\bar Y_n\geq \sigma\sqrt{2\log(1/\delta)}+C_q \frac{\Eb^{q/2}[|Y_1|^{q/2}]}{n^{1-2/q}\delta^{2/q}})\leq\delta,$$ where $\sigma^2 = \sum_{i\in[n]}\Eb[Y_i^2]/n=\Eb[\norm{X_1}_{\Sigma^{-1}}^4]-1$, and $\Eb^{q/2}[|Y_1|^{q/2}]\leq \Eb^{q/2}[\norm{X_1}_{\Sigma^{-1}}^{q}]+1$. This proves the second part where $q=4$ case follows from Chebyshev's ineqaulity.
\end{proof}
\begin{lemma}\label{lem:C4}
    Let $X_1,\ldots,X_N$ be not necessarily independent centered Gaussian random variables. Let $\sigma_{\rm max}^2=\max_{i\in[N]}{\rm Var}(X_i)$, then $\Pb(\max_{i\in[N]}\abs{X_i}>t)\leq Ne^{-t^2/(2\sigma_{\rm max}^2)},$ for $t>0$.
\end{lemma}
\begin{proof}[proof of Lemma~\ref{lem:C4}] Let $\sigma_i^2={\rm Var}(X_i)$. It follows from the union bound and that $1-\Phi(t/\sigma_i)\leq e^{-t^2/(2\sigma_i^2)}/2$ for $t>0$ that $\Pb(\max_{i\in[N]}\abs{X_i}>t)\leq \sum_{i\in[N]}\Pb(\abs{X_i}>t)=2\sum_{i\in[N]}(1-\Phi(t/\sigma_i))\leq \sum_{i\in[N]}e^{-t^2/(2\sigma_i^2)}\leq Ne^{-t^2/(2\sigma_{\rm max}^2)}.$
    
\end{proof}
\begin{lemma}\label{lem:C5} Let $\Ac = \set{a_j:j\in[m]}$ and let $$R_{n,\Ac} = \min_{j\in[m]}\sqrt{\frac{\sum_{i\in[n]}(a_j^\top(\psi(Z_i,\theta)-\bar\psi_1(\theta))^2}{na_j V(\theta)a_j}}.$$ For any $q>2$, $\Pb\left(R^2_{n,\Ac}\geq 1-c_{n,q}-e_{n,q}\right)\geq 1-2m/n^{q/2-1},$ where $c_{n,q}=L_q^2(\theta)(q\log(en)/(2n))^{1-2/q}$ and $e_{n,q}=cL_q^2(\theta)(\log^{1-1/q}(ed)m^{1/q}/n^{1-1/q})$ provided that $c_{n,q}\in(0,1)$.

\end{lemma}
\begin{proof}[proof of Lemma~\ref{lem:C5}.] It holds that $R_{n,\Ac}\geq D_1-D_2$ where
\begin{equation*}
    D_1 = \min_{j\in[m]}\sqrt{\frac{\sum_{i\in[n]}(a_j^\top(\psi(Z_i,\theta)-\mu(\theta))^2}{na_j^\top V(\theta)a_j}}, ~D_2 = \max_{j\in[m]}\frac{|a_j^\top(\bar\psi_1(\theta)-\mu(\theta))|}{(a_j^\top V(\theta)a_j)^{1/2}}.
\end{equation*}Let $ X_{ij}=(a_j^\top(\psi(Z_i,\theta)-\mu(\theta))^2$ for $i\in[n]$ and $j\in[m]$. An Lemma~\ref{lem:C3} with $\varepsilon=q/2-1$ implies that $\Pb(D_1\geq 1-c_{n,q})\geq 1-m/n^{q/2-1}$. Meanwhile, $\Dc_2$ can be controlled via Lemma~\ref{lem:B4.3}. 
\end{proof}

\begin{lemma}\label{lem:C6} Let $\Sigma =(\sigma_{jk})$ and $\Sigma' =(\sigma'_{jk})$ be $p\times p$ covariance matrices. Then,
\begin{equation}\label{eq:lem:C5.1}
    \max_{1\leq j\leq k\leq p}\Abs{\frac{\sigma'_{jk}}{\sqrt{\sigma'_{jj}\sigma'_{kk}}}-\frac{\sigma_{jk}}{\sqrt{\sigma_{jj}\sigma_{kk}}}}\leq 4\max_{1\leq j\leq k\leq p}\Abs{\frac{\sigma'_{jk}-\sigma_{jk}}{\sqrt{\sigma_{jj}\sigma_{kk}}}}.
\end{equation}
Suppose $X_1,\ldots,X_N\in\Real^p$ are independent random vectors where $X_i=(X_{i1},\ldots,X_{ip})^\top$ for all $i=1,\ldots,p$. Set
\begin{equation*}
    \hat\sigma_{jk}:=\frac{1}{N}\sum_{i=1}^NX_{ij}X_{i,k},\quad \sigma_{jk}=\Eb[\sigma'_{jk}],
\end{equation*} for all $j,k=1,\ldots,p$. Further, denote $\hat\rho_{jk}=\hat\sigma_{jk}/\sqrt{\hat\sigma_{jj}\hat\sigma_{kk}}$ and $\rho_{jk}=\sigma_{jk}/\sqrt{\sigma_{jj}\sigma_{kk}}$. Then, it follows from \eqref{eq:lem:C5.1} that
\begin{equation*}
    \max_{1\leq j\leq k\leq p}\abs{\hat\rho_{jk}-\rho_{jk}}\leq 4\max_{1\leq j\leq k\leq p}\Abs{\frac{\hat\sigma_{jk}-\sigma_{jk}}{\sqrt{\sigma_{jj}\sigma_{kk}}}}.
\end{equation*}
\end{lemma}
\begin{proof} We follow the proof of Lemma~19 of \cite{kuchibhotla2020berry}. Define
\begin{equation*}
    \Delta = \max_{1\leq j\leq k\leq p}\Abs{\frac{\sigma'_{jk}-\sigma_{jk}}{\sqrt{\sigma_{jj}\sigma_{kk}}}}.
\end{equation*} We may assume that $\Delta\leq 1/2$, otherwise \eqref{eq:lem:C5.1} follows immediately. From the definition, we get \begin{equation}\label{eq:lem:C5.2}
    \frac{1}{1+\Delta}\leq\frac{\sigma_{jj}}{\sigma'_{jj}}\leq\frac{1}{1-\Delta},
\end{equation} for all $1\leq j\leq p$. We note that
\begin{eqnarray*}
    \max_{1\leq j\leq k\leq p}\Abs{\frac{\sigma'_{jk}}{\sqrt{\sigma'_{jj}\sigma'_{kk}}}-\frac{\sigma_{jk}}{\sqrt{\sigma_{jj}\sigma_{kk}}}}&\leq&\frac{\abs{\sigma'_{jk}-\sigma_{jk}}}{\sqrt{\sigma'_{jj}\sigma'_{kk}}}+\frac{\abs{\sigma_{jk}}}{\sqrt{\sigma_{jj}\sigma_{kk}}}\Abs{\frac{\sqrt{\sigma_{jj}\sigma_{kk}}}{\sqrt{\sigma'_{jj}\sigma'_{kk}}}-1}\\
    &\leq&\Delta +\Abs{\frac{\sqrt{\sigma_{jj}\sigma_{kk}}}{\sqrt{\sigma'_{jj}\sigma'_{kk}}}-1},
\end{eqnarray*} where the first inequality follows from the triangle inequality, and for the second inequality, we used the definition of $\Delta$ and the fact that $\abs{\sigma_{jk}}\leq\sqrt{\sigma_{jj}\sigma_{kk}}$, which is a consequence of the Cauchy-Schwarz inequality. To bound the rightmost term, we use \eqref{eq:lem:C5.2} to obtain
\begin{equation*}
    \frac{1}{(1+\Delta)^2}\leq\frac{\sigma_{jj}\sigma_{kk}}{\sigma'_{jj}\sigma'_{kk}}\leq\frac{1}{(1-\Delta)^2},
\end{equation*} for all $j,k=1,\ldots,p$. This implies that
\begin{equation*}
    \Abs{\frac{\sqrt{\sigma_{jj}\sigma_{kk}}}{\sqrt{\sigma'_{jj}\sigma'_{kk}}}-1}\leq\frac{\Delta}{1-\Delta}.
\end{equation*} Therefore, as long as $\Delta_n\leq 1/2$,
\begin{equation*}
    \max_{1\leq j\leq k\leq p}\Abs{\frac{\sigma'_{jk}}{\sqrt{\sigma'_{jj}\sigma'_{kk}}}-\frac{\sigma_{jk}}{\sqrt{\sigma_{jj}\sigma_{kk}}}}\leq \frac{2\Delta}{1-\Delta}\leq 4\Delta.
\end{equation*}
    
\end{proof}

\begin{proposition}\label{prop:C2}
    Let $X_1,\ldots,X_n$ be centered independent random vectors in $\Real^p$ where $X_i=(X_{i1},\ldots,X_{ip})^\top$ for all $i=1,\ldots,n$. Suppose that there exists positive constants $b_1$ and $b_2$, and a sequence of positive reals $\set{B_n\geq1}$ such that
    \begin{equation}\label{eq:propB.1:1}
        b_1\leq \min_{1\leq j\leq p}\left(\frac{1}{n}\sum_{i=1}^n\Eb[X_{ij}^2]\right)^{1/2}\quad\mbox{and}\quad\max_{1\leq j\leq p}\left(\frac{1}{n}\sum_{i=1}^n\Eb[X_{ij}^4]\right)^{1/4}\leq b_2 B_n.
    \end{equation} Furthermore, suppose there exists a sequence of reals $\set{D_n\geq1}$ such that
    \begin{equation}\label{eq:propB.1:2}
        \max_{1\leq i\leq n}\left(\Eb[\norm{X_i}_\infty^q]\right)^{1/q}\leq D_n,
    \end{equation} for some $q>2$. Let $\Ac$ be a class of hyperrectangles in $\Real^p$ and let $G\sim\Nc(0_P,\Sigma_n)$ where $\Sigma_n={\rm Var}(n^{-1/2}\sum_{i=1}^nX_i)$. Then, there exists a constant $C=C(b_1,b_2,q)>0$ such that
    \begin{equation*}
        \sup_{A\in\Ac}\Abs{\Pb\left(\frac{1}{\sqrt{n}}\sum_{i=1}^nX_i\in A\right)-\Pb\left(G\in A\right)}\leq C\left[\left(\frac{B_n^4\log^5(p)}{n}\right)^{1/4}+\frac{D_n\log^{3/2}(p)}{n^{1/2-1/q}}\right].
    \end{equation*}
\end{proposition}

This proposition refines Theorem~2.5 of \cite{chernozhuokov2022improved} where the author assumed $D_n=B_n^2$. The current version is more informative since $B_n^2$ and $D_n$ generally scale differently.

\begin{proof}[proof of Proposition~\ref{prop:C2}] We shall show that
\begin{equation*}
        \chi_n:=\sup_{y\in|\Real^p}\Abs{\Pb\left(\frac{1}{\sqrt{n}}\sum_{i=1}^nX_i\preceq y\right)-\Pb\left(G\preceq y\right)}\leq C\left[\left(\frac{B_n^4\log^5(p)}{n}\right)^{1/4}\right],
\end{equation*} which implies the desired result. Here, $x\preceq y$ denotes $\ev_j^\top x\leq \ev_j^\top y$ for all $j=1,\ldots,p$. We use Lemma~A.1 of \cite{chernozhuokov2022improved} which states that
\begin{eqnarray}\label{eq:pf_propA.1:1}
    \chi_n\leq C\bigg[\phi\log^2(p)\sqrt{\Delta_1}+\phi^3\log^{7/2}(p)\Delta_1+\log(p)\sqrt{\Delta_2(\phi)}\nonumber\\
    +\phi\log^{3/2}(p)\Delta_2(\phi)+\log^{3/2}(p)\sqrt{\Delta_3(\phi)}+\frac{\log^{1/2}(p)}{\phi}\bigg],
\end{eqnarray} for all $\phi>0$ where $C>0$ only depends on $b_1$ and
\begin{eqnarray*}
    \Delta_1 &=& \frac{1}{n^2}\max_{1\leq j\leq p}\sum_{i=1}^n\Eb[X_{ij}^4],\quad \Delta_2(\phi)=\max_{1\leq j\leq p}\sum_{i=1}^n\Eb\left[Y_{ij}^2\mathbbm{1}\set{\norm{Y_i}_\infty>(\phi\log(p))^{-1}}\right],\\
    \Delta_3(\phi) &=& \Eb\left[\max_{1\leq i\leq n}\norm{Y_i}_\infty^2\mathbbm{1}\set{\norm{Y_i}_\infty>(\phi\log(p))^{-1}}\right],
\end{eqnarray*} and $Y_i = (X_i-\tilde{X}_i)/\sqrt{n}$ for all $i=1,\ldots,n$ and $\set{\tilde X_1,\ldots,\tilde X_n}$ is an independent copy of $\set{X_1,\ldots,X_n}$. First, it is immediate that
\begin{equation}\label{eq:pf_propA.1:3}
    \Delta_1\leq n^{-1}B_n^2.
\end{equation} To analyze $\Delta_2(\phi)$, we first note that
\begin{equation}\label{eq:pf_propA.1:2}
    \max_{1\leq j\leq p}\frac{1}{n}\sum_{i=1}^n\Eb[Y_{ij}^4]\lesssim n^{-2}B_n^4,\quad\mbox{and}\quad \max_{1\leq i\leq n}\Eb\left[\norm{Y_i}_\infty^q\right]\lesssim n^{-q/2}D_n^q,
\end{equation} where $\lesssim$ denotes an inequality that holds up to constant depending only on $b_1,b_2$ and $q$. We observe that
\begin{eqnarray*}
    Y_{ij}^2\mathbbm{1}\set{\norm{Y_i}_\infty>A}&=&Y_{ij}^2\mathbbm{1}\set{\abs{Y_{ij}}>A,\norm{Y_i}_\infty>A}+Y_{ij}^2\mathbbm{1}\set{\abs{Y_{ij}}\leq A,\norm{Y_i}_\infty>A}\\
    &=&Y_{ij}^2\mathbbm{1}\set{\abs{Y_{ij}}>A}+Y_{ij}^2\mathbbm{1}\set{\abs{Y_{ij}}\leq A,\norm{Y_i}_\infty>A}\\
    &\leq&Y_{ij}^2\mathbbm{1}\set{\abs{Y_{ij}}>A}+A^2\mathbbm{1}\set{\norm{Y_i}_\infty>A},
\end{eqnarray*} for any $A>0$. Consequently, we have
\begin{eqnarray*}
    \Eb\left[Y_{ij}^2\mathbbm{1}\set{\norm{Y_i}_\infty>A}\right]&\leq& \Eb\left[Y_{ij}^2\mathbbm{1}\set{\abs{Y_{ij}}>A}\right]+A^2\Pb\left(\norm{Y_i}_\infty>A\right)\\
    &\leq&\Eb\left[\frac{Y_{ij}^4}{A^2}\mathbbm{1}\set{\abs{Y_{ij}}>A}\right]+A^2\frac{\Eb[\norm{Y_i}_\infty^q]}{A^q}\\
    &\leq&\frac{\Eb[Y_{ij}^4]}{A^2}+\frac{\Eb[\norm{Y_i}_\infty^q]}{A^{q-2}}.
\end{eqnarray*} We take $A=(\phi\log(p))^{-1}$, then \eqref{eq:pf_propA.1:2} with leads to
\begin{equation}\label{eq:pf_propA.1:4}
    \Delta_2(\phi)\lesssim \frac{(\phi\log(p))^2B_n^4}{n}+\frac{(\phi\log(p))^{q-2}D_n^q}{n^{q/2-1}}.
\end{equation} The quantity $\Delta_3(\phi)$ can be controlled as
\begin{eqnarray}\label{eq:pf_propA.1:5}
    \Delta_3(\phi)\leq \sum_{i=1}^n \Eb\left[\norm{Y_i}_\infty^2\mathbbm{1}\set{\norm{Y_i}_\infty>(\phi\log(p))^{-1}}\right]\leq \sum_{i=1}^n (\phi\log(p))^{q-2}\Eb\left[\norm{Y_i}_\infty^q\right]\lesssim\frac{(\phi\log(p))^{q-2}D_n^q}{n^{q/2-1}},
\end{eqnarray} where the last inequality follows from \eqref{eq:pf_propA.1:2}. Hence, combining \eqref{eq:pf_propA.1:3}, \eqref{eq:pf_propA.1:4}, \eqref{eq:pf_propA.1:5}, and \eqref{eq:pf_propA.1:1} yields
\begin{eqnarray*}
    \chi_n&\lesssim& \frac{\phi\log^2(p)B_n}{\sqrt{n}}+\frac{\phi^3\log^{7/2}(p)B_n^2}{n}+\frac{\phi\log^2(p)B_n^2}{\sqrt{n}}+\frac{\phi^{q/2-1}\log^{q/2}(p)D_n^{q/2}}{n^{q/4-1/2}}\\
    &&+\frac{\phi^3\log^{7/2}(p)B_n^4}{n}+\frac{\phi^{q-1}\log^{q-1/2}(p)D_n^{q}}{n^{q/2-1}}+\frac{\phi^{q/2-1}\log^{q/2+1/2}(p)D_n^{q/2}}{n^{q/4-1/2}}+\frac{\log^{1/2}(p)}{\phi},
\end{eqnarray*} for all $\phi>0$. Here, we used that $\sqrt{a+b}\leq\sqrt{a}+\sqrt{b}$. Since $B_n\geq1$ and $\log(p)\geq1$, we can further deduce that
\begin{equation*}
    \chi_n\lesssim\frac{\phi\log^2(p)B_n^2}{\sqrt{n}}+\frac{\phi^3\log^{7/2}(p)B_n^4}{n}+\frac{\phi^{q-1}\log^{q-1/2}(p)D_n^{q}}{n^{q/2-1}}+\frac{\phi^{q/2-1}\log^{q/2+1/2}(p)D_n^{q/2}}{n^{q/4-1/2}}+\frac{\log^{1/2}(p)}{\phi}.
\end{equation*} We choose
\begin{equation*}
    \phi^{-1}=\frac{B_n\log^{3/4}(p)}{n^{1/4}}+\frac{D_n\log(p)}{n^{1/2-1/q}}.
\end{equation*} This leads to
\begin{eqnarray*}
    \frac{\phi\log^2(p)B_n^2}{\sqrt{n}}+\frac{\phi^3\log^{7/2}(p)B_n^4}{n}&\leq&\left(\frac{B_n\log^{3/4}(p)}{n^{1/4}}\right)^{-1}\frac{\log^2(p)B_n^2}{\sqrt{n}}+\left(\frac{B_n\log^{3/4}(p)}{n^{1/4}}\right)^{-3}\frac{\phi\log^{7/2}(p)B_n^4}{n}\\
    &=&\frac{2B_n\log^{5/4}(p)}{n^{1/4}}.
\end{eqnarray*} Furthermore, 
\begin{eqnarray*}
    &&\frac{\phi^{q-1}\log^{q-1/2}(p)D_n^{q}}{n^{q/2-1}}+\frac{\phi^{q/2-1}\log^{q/2+1/2}(p)D_n^{q/2}}{n^{q/4-1/2}}\\
    &\leq&\left(\frac{D_n\log(p)}{n^{1/2-1/q}}\right)^{-(q-1)}\frac{\log^{q-1/2}(p)D_n^{q}}{n^{q/2-1}}+\left(\frac{D_n\log(p)}{n^{1/2-1/q}}\right)^{-(q/2-1)}\frac{\log^{q/2+1/2}(p)D_n^{q/2}}{n^{q/4-1/2}}\\
    &=&\frac{D_n\log^{1/2}(p)}{n^{1/2-1/q}}+\frac{D_n\log^{3/2}(p)}{n^{1/2-1/q}}\leq\frac{D_n\log^{3/2}(p)}{n^{1/2-1/q}}.
\end{eqnarray*} Finally, one has
\begin{equation*}
    \frac{\log^{1/2}(p)}{\phi}=\frac{B_n\log^{5/4}(p)}{n^{1/4}}+\frac{D_n\log^{3/2}(p)}{n^{1/2-1/q}}.
\end{equation*} This proves the proposition.
\end{proof}

\begin{proposition}[Concentration Inequality for Maximal Independent Sum] \label{prop:C3}Suppose that $X_1,\ldots,X_n$ are independent centered random vectors in $\Real^p$ where $X_i=(X_{i1},\ldots,X_{ip})^\top$ for all $i=1,\ldots,n$. Let
\begin{equation*}
    B_{n,q}:=\max_{1\leq j\leq p}\left(\frac{1}{n}\sum_{i=1}^n\Eb[\abs{X_{ij}}^q]\right)^{1/q} \mbox{ for }q\geq1, \mbox{ and } v_n=\max_{1\leq j\leq p}\max_{1\leq i\leq n}\left(\Eb[X_{ij}^2]\right)^{1/2}.
\end{equation*}Let $A=\max_{1\leq j\leq p}\abs{n^{-1}\sum_{i=1}^nX_{ij}}$, then
\begin{equation}\label{eq:propC.2:1}
    \Eb[A]\leq B_{n,2}\sqrt{\frac{6\log(1+p)}{n}}+\sqrt{2}B_{n,q}\left(\frac{3\log(1+p)(2p)^{1/(q-1)}}{n}\right)^{1-1/q},
\end{equation} for $q\geq1$, and
\begin{equation}\label{eq:propC.2:2}
    \left(\Eb[A-\Eb[A])_+^q]\right)^{1/q}\leq \sqrt{\kappa q}\frac{B_{n,q}^2(p^{2/q}+1)}{n}+\kappa q\left(\frac{p^{1/q}B_{n,q}}{n^{1-1/q}}+\frac{v_n}{n}\right),
\end{equation} for $q\geq2$ where $\kappa=\sqrt{e}/(\sqrt{e}-1)$. Combining \eqref{eq:propC.2:1} and \eqref{eq:propC.2:2} together with Markov's inequality readily gives the concentration inequality for $A$.
    
\end{proposition}
\begin{proof} From Markov's inequality, we get
\begin{equation*}
    \Pb\left(A\geq \Eb[A]+t\right)\leq t^{-q}\Eb\left[(A-\Eb[A])_+^q\right],
\end{equation*} for any $t>0$.
We use Proposition~B.1 of \cite{kuchibhotla2022least} to control $\Eb[A]$, and use Theorem~15.14 of \cite{Boucheron2013concentration} to control $\Eb[(A-\Eb[A])_+]^q$. First, we note that
\begin{eqnarray*}
    \sum_{i=1}^n\Eb\left[\max_{1\leq j\leq p}\Abs{\frac{X_{ij}}{n}}^q\right]\leq \frac{p}{n^q}\max_{1\leq j\leq p}\sum_{i=1}^n\Eb\abs{X_{ij}}^q=\frac{pB_{n,q}^q}{n^{q-1}}.\\
\end{eqnarray*} Hence, an application of Proposition~B.1 of \cite{kuchibhotla2022least} yields that
\end{proof}
\begin{eqnarray*}
    \Eb[A]\leq B_{n,2}\sqrt{\frac{6\log(1+p)}{n}}+\sqrt{2}B_{n,q}\left(\frac{3\log(1+p)(2p)^{1/(q-1)}}{n}\right)^{1-1/q},
\end{eqnarray*} for any $q\geq1$. Meanwhile, an application of Theorem~15.14 of \cite{Boucheron2013concentration} leads to that for any $q\geq 2$,
\begin{eqnarray*}
        \left(\Eb[A-\Eb[A])_+^q]\right)^{1/q}&\leq& \sqrt{\kappa q}\left[\Eb\left(\max_{1\leq j\leq p}\frac{1}{n^2}\sum_{i=1}^n X_{ij}^2\right)+\max_{1\leq j\leq p}\frac{1}{n^2}\sum_{i=1}^n\Eb [X_{ij}^2]\right]\\
        &&+\,\kappa q\left[\frac{1}{n}\left\{\Eb\left(\max_{1\leq i\leq n}\max_{1\leq j\leq p}\abs{X_{ij}}^q \right)\right\}^{1/q}+\frac{1}{n}\max_{1\leq i\leq n}\max_{1\leq j\leq p}\left(\Eb[X_{ij}^2]\right)^{1/2}\right],
    \end{eqnarray*} where $\kappa=\sqrt{e}/(\sqrt{e}-1)<1.542$. We further note that
    \begin{equation*}
        \Eb\left(\max_{1\leq j\leq p}\frac{1}{n^2}\sum_{i=1}^n X_{ij}^2\right)\leq \frac{1}{n^2}\sum_{i=1}^n\left(\Eb\max_{1\leq j\leq p}\abs{X_{ij}}^q\right)^{2/q}\leq \frac{p^{2/q}B_{n,q}^2}{n},\quad \max_{1\leq j\leq p}\frac{1}{n^2}\sum_{i=1}^n\Eb [X_{ij}^2]=\frac{B_{n,2}^2}{n}.
    \end{equation*} Moreover,
    \begin{eqnarray*}
        \frac{1}{n}\left\{\Eb\left(\max_{1\leq i\leq n}\max_{1\leq j\leq p}\abs{X_{ij}}^q \right)\right\}^{1/q}\leq\frac{1}{n}\left\{\sum_{j=1}^p\sum_{i=1}^n\Eb\abs{X_{ij}}^q\right\}^{1/q}\leq\frac{p^{1/q}B_{n,q}}{n^{1-1/q}},\quad\frac{1}{n}\max_{1\leq i\leq n}\max_{1\leq j\leq p}\left(\Eb[X_{ij}^2]\right)^{1/2}=\frac{v_n}{n}.
    \end{eqnarray*} This gives the desired result.

\begin{proposition}\label{prop:C4} For $d$-dimensional random vector $X_1,\ldots,X_n\sim X$, suppose there exists $q\geq4$ and $K\geq1$ such that
\begin{equation*}
    \sup_{u\in\Sb^{d-1}}\Eb[|u^\top\Sigma^{-1}X|^q]\leq K^q,
\end{equation*} where $\Sigma = \Eb[XX^\top]$. Define $\Dc^\Sigma=\norm{\Sigma^{-1/2}\hat\Sigma\Sigma^{-1/2}-I_d}_{\rm op}$ where $\hat\Sigma = n^{-1}\sum_{i=1}^nX_iX_i^\top$, then there exists a universal constant $C>0$ such that
    \begin{equation}\label{eq:conc_D_sigma}
        \Pb\left(\Dc^\Sigma\leq CK_x^2\left[\sqrt{\frac{d\log(2d/\delta)}{n}}+\left(\log^{3/2}(2d/\delta)+\frac{\log(2d/\delta)}{\delta^{2/q}}\right)\frac{d}{n^{1-2/q}}\right]\right)\geq 1-\delta,
    \end{equation} for all $\delta\in(0,1)$.
\end{proposition}
This proposition was previously established in Proposition~6 of \citet{chang2023inference}. For completeness, we provide a self-contained proof.

\begin{proof}[proof of Proposition~\ref{prop:C4}.]
This proposition applies Theorem 2.8 in \cite{brailovskaya2023universality}, and we first need to define the following matrix parameters. Let $\Omega = \Sigma^{-1/2}\hat\Sigma\Sigma^{-1/2}$, and define
\begin{align}\label{eq:matrixparameters}
    \sigma(\Omega):=\norm{\Eb[(\Omega -I)^2]}^{1/2}_{\rm op},\quad
    \sigma_*(\Omega):=\sup_{u,v\in\Sb^{d-1}}\Eb[|u^\top (\Omega-I)v|^2]^{1/2},\\
    Z_i = \frac{\Sigma^{-1/2}X_iX_i^\top\Sigma^{-1/2}-I_d}{n}~~\mbox{for}~i\in[n],\quad \bar R(\Omega) = \Eb\left[\max_{i\in[n]}\norm{Z_i}_{\rm op}^2\right]^{1/2}.
\end{align}
Now, we define a moment-matching $d\times d$ random matrix $G$ with jointly Gaussian entries such that $\Eb [G] = \Eb[\Omega]=I_d$ and ${\rm Cov}(G)= {\rm Cov}(\Omega)$, that is, for $i,j,k,l\in[d]$,
\begin{equation}\label{eq:cov_of_mat}
    \Eb[(G-I)_{ij}(G-I)_{kl}]=\Eb[(\Omega-I)_{ij}(\Omega-I)_{kl}].
\end{equation} The first two moments uniquely define the the distribution of $G$; see Section 2.1 of \cite{brailovskaya2023universality}. Denote the collection of eigenvalues of matrix $M$ as ${\rm sp}(M)$. Finally, Theorem~2.8 of \cite{brailovskaya2023universality} applies and yields that there exists a universal constant $C>0$ such that
\begin{align}\label{eq:universal_ineq_Brail}
    &\Pb\Bigg(d_H\left({\rm sp}(\Omega),{\rm sp}(G)\right)\geq C e_R(t),~~\max_{i\in[n]}\norm{Z_i}_{\rm op}\leq R\Bigg)\leq de^{-t},
\end{align} for all $t\geq0$ and $R\geq\bar R(\Omega)^{1/2}\sigma_*(\Omega)^{1/2}+\sqrt{2}\bar R(\Omega)$ where $e_R(t) = \sigma_*(\Omega)t^{1/2}+R^{1/3}\sigma(X)^{2/3}t^{2/3}+Rt$. Here $d_H$ denotes two-sided Hausdorff distance.
    
Before carefully analyzing \eqref{eq:universal_ineq_Brail}, we note that
\begin{align*}
    &\Dc^\Sigma = \max\set{|\lambda_{\rm max}(\Omega)-1|,|\lambda_{\rm min}(\Omega)-1|}\\
    &\quad\leq \max\set{|\lambda_{\rm max}(\Omega)-\lambda_{\rm max}(G)|,|\lambda_{\rm min}(\Omega)-\lambda_{\rm min}(G)|} + \max\set{|\lambda_{\rm max}(G)-1|,|\lambda_{\rm min}(G)-1|}\\
    &\quad\leq d_H({\rm sp}(\Omega),{\rm sp}(G)) + \norm{G-I_d}_{\rm op}.
\end{align*} Since $G$ has jointly Gaussian entries, it can be expressed as the linear combination of independent standard Gaussian variables, $G = I_d + \sum_{i=1}^N \Av_i\gv_i,$ for some deterministic symmetric matrices $\Av_1,\ldots,\Av_N$, and i.i.d. standard Gaussian variables $\gv_1,\ldots,\gv_N$. Therefore, Theorem 4.1.1 of \cite{tropp2015introduction} applies and yields
    \begin{equation}\label{eq:tropp2015}
        \Pb\left(\norm{G-I_d}_{\rm op}\geq t\right)\leq 2d\exp\left(-\frac{t^2}{2\sigma_*^2}\right)\mbox{ for all }t\geq0,
    \end{equation} where $\sigma_*^2=\norm{\Eb[(G-I)(G-I)^\top]}_{\rm op}=\sigma_*^2(\Omega)$. Combining this with \eqref{eq:universal_ineq_Brail}, we get
\begin{align}\label{eq:universal_ineq_Brail_new}
    &\Pb\Bigg(\Dc^\Sigma\geq C \varepsilon_R(t),~~\max_{i\in[n]}\norm{Z_i}_{\rm op}\leq R\Bigg)\leq 3de^{-t},
\end{align} for all $t\geq0$ and $R\geq\bar R(\Omega)^{1/2}\sigma_*(\Omega)^{1/2}+\sqrt{2}\bar R(\Omega)$ where $\varepsilon_R(t) = (1+\sqrt{2})\sigma_*(\Omega)t^{1/2}+R^{1/3}\sigma(X)^{2/3}t^{2/3}+Rt$.
Now, we will analyze the matrix parameters in \eqref{eq:matrixparameters}:
\begin{align*}
    &\sigma(\Omega)^2 = \norm{\Eb[(\Omega-I_d)^2]}_{\rm op}=\sup_{u\in\Sb^{d-1}}\Abs{\Eb[u^\top (\Omega-I_d)^2u]}\\
    &~~=\sup_{u\in \Sb^{d-1}}\Abs{\Eb\left[u^\top\left({\textstyle\sum}_{i=1}^nZ_i\right)^2u\right]} = \sup_{u\in \Sb^{d-1}}\Abs{\Eb\left[u^\top\left({\textstyle\sum}_{i=1}^nZ_i^2\right)u\right]}\\
    &\quad=n^{-1}\sup_{u\in\Sb^{d-1}}\Abs{\Eb[u^\top(\Sigma^{-1/2}XX^\top\Sigma^{-1/2}-I_d)^2u]}\\
    &\qquad\leq n^{-1}\sup_{u\in\Sb^{d-1}}\Abs{\Eb[(u^\top\Sigma^{-1/2}X)^2 \norm{X}_{\Sigma^{-1}}^2]-1}\leq \frac{K^4d}{n},
\end{align*} where we have used $Z_i$'s are mean-zero and independent across $i\in[n]$. On the other hand,
\begin{align*}
    &\sigma_*(\Omega)^2 = \sup_{u\in\Sb^{d-1}}\Eb[\abs{u^\top ({\textstyle\sum}_{i=1}^n Z_i) u}^2]=\sup_{u\in\Sb^{d-1}}\Eb[{\textstyle\sum}_{i=1}^n (u^\top Z_i u)^2]\\
    &\quad= n^{-1}\sup_{u\in\Sb^{d-1}}\Eb[(u^\top\Sigma^{-1/2}X)^2-1)^2]=n^{-1}\sup_{u\in\Sb^{d-1}}\Eb[(u^\top\Sigma^{-1/2}X)^4]-1\leq \frac{K^4}{n}.
\end{align*}In order to analyze $\bar R(\Omega)$, we note that $\norm{Z_i}_{\rm op} = n^{-1}\abs{\norm{X_i}_{\Sigma^{-1}}^2-1}$. Therefore,
\begin{equation*}
    \bar R(\Omega)^2\leq\frac{1}{n^2} + \frac{1}{n^2}\Eb\left[\max_{i\in[n]}\norm{X_i}^4_{\Sigma^{-1}}\right].
\end{equation*} Moreover, we get from Jensen's inequality and union bound that
\begin{align*}
    \Eb\left[\max_{i\in[n]}\norm{X_i}^4_{\Sigma^{-1}}\right] \leq \Eb^{4/q}\left[\max_{i\in[n]}\norm{X_i}^q_{\Sigma^{-1}}\right]\leq n^{4/q}\Eb^{4/q}\left[\norm{X_i}^q_{\Sigma^{-1}}\right]\leq K^4n^{4/q}d^2.
\end{align*} Hence, $\bar R(\Omega)^2\leq 1/n^2 + K^4d^2/n^{2-4/q}\leq 2K^4d^2/n^{2-4/q}$. Following the proof of Proposition B.6 of \cite{chang2023inference}, we let
\begin{equation*}
    R(\delta)=2K^2\left[\frac{d}{\delta^{2/q}n^{1-2/q}}+\frac{d^{3/4}}{n^{3/4-1/q}}\right],
\end{equation*} for $\delta\in(0,1]$. This ensures that $R(\delta)\geq \bar R(\Omega)^{1/2}\sigma_*(\Omega)^{1/2}+\sqrt{2}\bar R(\Omega)$ for all $\delta\in(0,1]$. Moreover, it follows from the union bound and Markov inequality that
\begin{align}\label{eq:propD4:1}
    &\Pb\left(\max_{i\in[n]}\norm{Z_i}_{\rm op}\geq R(\delta)\right)=\Pb\left(\max_{1\leq i\leq n}\Abs{\norm{ X_i}_{\Sigma^{-1}}^2-1}\geq nR(\delta)\right)\leq n\Pb\left(\Abs{\norm{ X}_{\Sigma^{-1}}^2-1}\geq nR(\delta)\right)\nonumber\\
    &\leq n(nR(\delta))^{-q/2}\Eb[\abs{\norm{ X}_{\Sigma^{-1}}^2-1}^{q/2}]\leq\frac{2^{q/2-1}K^{q}d^{q/2}}{n^{q/2-1}R(\delta)^{q/2}}\leq\frac{\delta}{2}.
\end{align} From Young's inequality, it follows that  $3R(\delta)^{1/3}\sigma(\Omega)^{2/3}t^{2/3}\leq 2\sigma(\Omega) t^{1/2}+R(\delta)t$, and this implies that for some universal constant $C>0$,
\begin{align*}
    &\varepsilon_{R(\delta)}(t) = (1+\sqrt{2})\sigma_*(\Omega)t^{1/2}+R(\delta)^{1/3}\sigma(\Omega)^{2/3}t^{2/3}+R(\delta)t\\
    &\quad \leq CK^2\left(\sqrt{\frac{dt}{n}}+t\left[\frac{d}{\delta^{2/q}n^{1-2/q}}+\frac{d^{3/4}}{n^{3/4-1/q}}\right]\right).
\end{align*} Moreover, we note that
\begin{equation*}
    \frac{d^{3/4}t}{n^{3/4-1/q}} = \left(\sqrt{\frac{dt}{n}}\right)^{1/2}\left(\frac{t^{3/2}d}{n^{1-q/2}}\right)^{1/2}\leq \sqrt{\frac{dt}{n}} + \frac{t^{3/2}d}{n^{1-2/q}}.
\end{equation*} Therefore, the inequality in \eqref{eq:universal_ineq_Brail_new} becomes
\begin{equation*}
    \Pb\left(\Dc^\Sigma\geq CK^2\left(\sqrt{\frac{dt}{n}}+\frac{t^{3/2}d}{n^{1-2/q}} +\frac{td}{\delta^{2/q}n^{1-2/q}}\right),~\max_{i\in[n]}\norm{Z_i}_{\rm op}\leq R(\delta) \right)\leq de^{-t},
\end{equation*} for all $t\geq 0$. Taking $t = \log(2d/\delta)$ and combining this with \eqref{eq:propD4:1} leads to the desired results.
\end{proof}
\begin{lemma}\label{lem:C7} Assume \ref{asmp:glm2}. For any $0<q_1,q_2$ such that $q_1/q_x+q_2/q_y\leq1$, it holds that
    \begin{equation*}
        \sup_{u\in\Sb^{p-1}}\Eb\left[\frac{1}{n}\sum_{i=1}^n\abs{u^\top\Sigma^{-1/2}X_i}^{q_1}\abs{Y_i-\ell'(X_i^\top\theta_0)}^{q_2}\right]\leq K_x^{q_1}K_y^{q_2},
    \end{equation*} and
    \begin{equation*}
        \Eb\left[\frac{1}{n}\sum_{i=1}^n\norm{\Sigma^{-1/2}X_i}_2^{q_1/2}\abs{Y_i-\ell'(X_i^\top\theta_0)}^{q_2}\right]\leq K_x^{q_1}K_y^{q_2}d^{q_1/2}.
    \end{equation*}
\end{lemma}
\begin{proof}
    Let $s=(q_1/q_x+q_2/q_y)^{-1}\geq1$. It follows from Jensen's inequality and H\"older's inequality that
    \begin{eqnarray*}
        \Eb\left[\norm{\Sigma^{-1/2}X_i}_2^{q_1}\abs{Y_i-\ell'(X_i^\top\theta_0)}^{q_2}\right]&\leq& \left(\Eb\left[\norm{\Sigma^{-1/2}X_i}_2^{sq_1}\abs{Y_i-\ell'(X_i^\top\theta_0)}^{sq_2}\right]\right)^{1/s}\\
        &\leq&\left(\Eb\left[\norm{\Sigma^{-1/2}X_i}_2^{q_x}\right]\right)^{q_1/q_x}\left(\Eb\left[\abs{Y_i-\ell'(X_i^\top\theta_0)}^{q_y}\right]\right)^{q_2/q_y}.
    \end{eqnarray*} Applications of Jensen's inequality and H\"older's inequality once more leads to
    \begin{eqnarray*}
        &&\frac{1}{n}\sum_{i=1}^n\left(\Eb\left[\norm{\Sigma^{-1/2}X_i}_2^{q_x}\right]\right)^{q_1/q_x}\left(\Eb\left[\abs{Y_i-\ell'(X_i^\top\theta_0)}^{q_y}\right]\right)^{q_2/q_y}\\
        &\leq& \left(\frac{1}{n}\sum_{i=1}^n\left(\Eb\left[\norm{\Sigma^{-1/2}X_i}_2^{q_x}\right]\right)^{sq_1/q_x}\left(\Eb\left[\abs{Y_i-\ell'(X_i^\top\theta_0)}^{q_y}\right]\right)^{sq_2/q_y}\right)^{1/s}\\
        &\leq&\left(\left\{\frac{1}{n}\sum_{i=1}^n\Eb\left[\norm{\Sigma^{-1/2}X_i}_2^{q_x}\right]\right\}^{sq_1/q_x}\left\{\frac{1}{n}\sum_{i=1}^n\Eb\left[\abs{Y_i-\ell'(X_i^\top\theta_0)}^{q_y}\right]\right\}^{sq_2/q_y}\right)^{1/s}\\
        &=&\left\{\frac{1}{n}\sum_{i=1}^n\Eb\left[\norm{\Sigma^{-1/2}X_i}_2^{q_x}\right]\right\}^{q_1/q_x}\left\{\frac{1}{n}\sum_{i=1}^n\Eb\left[\abs{Y_i-\ell'(X_i^\top\theta_0)}^{q_y}\right]\right\}^{q_2/q_y}\leq K_x^{q_1}K_y^{q_2}d^{q_1/2},
    \end{eqnarray*} where the last inequality follows from \ref{asmp:glm2}. In particular, one has from Jensen's inequality that
    \begin{eqnarray*}
        &&\frac{1}{n}\sum_{i=1}^n\Eb\left[\norm{\Sigma^{-1/2}X_i}_2^{q_x}\right]\\
        &=&\frac{1}{n}\sum_{i=1}^n\Eb\left[\left(\sum_{j=1}^d\abs{\ev_j^\top\Sigma^{-1/2}X_i}^2\right)^{q_x/2}\right]=\frac{d^{q_x/2}}{n}\sum_{i=1}^n\Eb\left[\left(\frac{1}{d}\sum_{j=1}^d\abs{\ev_j^\top\Sigma^{-1/2}X_i}^2\right)^{q_x/2}\right]\\
        &\leq&\frac{d^{q_x/2}}{n}\sum_{i=1}^n\Eb\left[\frac{1}{d}\sum_{j=1}^d\abs{\ev_j^\top\Sigma^{-1/2}X_i}^{q_x}\right]\leq K_x^{q_x}d^{q_x/2}.
    \end{eqnarray*}
\end{proof}

\begin{proposition}\label{prop:C5}Grant \ref{asmp:lm.1} and \ref{asmp:lm.2}. Then,
\begin{equation}\label{eq:prop.11.1}
    \Pb\left(\frac{1}{n^2}\sum_{i=1}^n\norm{\Sigma^{-1/2}X_i}_2^2(Y_i-X_i^\top\beta_0)^2\geq K_x^2K_y^2\frac{d}{n}\left(1+\frac{1}{\sqrt{n\delta}}\right)\right)\leq \delta,
\end{equation} for all $\delta>0$. Moreover, with probability at least $1-\delta$,
\begin{align}\label{eq:prop.11.2}
    (\beta_0-\hat\beta)^\top \hat\Sigma_1(\beta_0-\hat\beta)&\leq\left(1-9K_x\sqrt{\frac{d+2\log(4/\delta)}{n}}\right)_+^{-1}\nonumber\\
    &\quad\times\left(8\frac{\overline{\lambda}_V}{\underline{\lambda}_\Sigma }\frac{d+2\log(4/\delta)}{n}+C_{q_{xy}}K_x^2K_y^2\frac{d}{\delta^{2/q_{xy}}n^{2-2/q_{xy}}}\right)
\end{align} where the constant $C_{q_{xy}}$ only depends on $q_{xy}$ and $\hat\beta$ is the OLSE computed from $\set{(Y_i,X_i):i\in[n]}$.
    
\end{proposition}
\begin{proof}
    Note that \eqref{eq:prop.11.1} follows from Chebyshev's inequality. In particular, an application of Lemma~\ref{lem:C7} with $q_1=2$ and $q_2=2$ leads to
    \begin{equation*}
        \Eb\left[\frac{1}{n^2}\sum_{i=1}^n\norm{\Sigma^{-1/2}X_i}_2^2(Y_i-X_i^\top\beta_0)^2\right]\leq K_x^2K_y^2\frac{d}{n}.
    \end{equation*} In addition, an application of Lemma~\ref{lem:C7} with $q_1=4$ and $q_2=4$ gives
    \begin{eqnarray*}
        {\rm Var}\left[\frac{1}{n^2}\sum_{i=1}^n\norm{\Sigma^{-1/2}X_i}_2^2(Y_i-X_i^\top\beta_0)^2\right]&=&\frac{1}{n^4}\sum_{i=1}^n{\rm Var}\left[\norm{\Sigma^{-1/2}X_i}_2^2(Y_i-X_i^\top\beta_0)^2\right]\\
        &\leq& \frac{1}{n^4}\sum_{i=1}^n\Eb\left[\norm{\Sigma^{-1/2}X_i}_2^4(Y_i-X_i^\top\beta_0)^4\right]\\
        &\leq&K_x^4K_y^4\frac{d^2}{n^3}.
    \end{eqnarray*} This proves \eqref{eq:prop.11.1}. To prove \eqref{eq:prop.11.2}, we note from the definition of $\hat\beta$ that
    \begin{equation*}
        \hat\beta-\beta_0=\hat\Sigma^{-1}\Sigma^{1/2}\frac{1}{n}\sum_{i=1}^n\Sigma^{-1/2}X_i(Y_i-X_i^\top\beta_0).
    \end{equation*} Therefore, we can write
    \begin{eqnarray*}
        &&(\beta_0-\hat\beta)^\top \hat\Sigma_1(\beta_0-\hat\beta)\\
        &=&\left(\frac{1}{n}\sum_{i=1}^n\Sigma^{-1/2}X_i(Y_i-X_i^\top\beta_0)\right)^\top\Sigma^{1/2}\hat\Sigma^{-1}\Sigma^{1/2}\left(\frac{1}{n}\sum_{i=1}^n\Sigma^{-1/2}X_i(Y_i-X_i^\top\beta_0)\right)\\
        &\leq&\frac{1}{\lambda_{\rm min}\left(\Sigma^{-1/2}\hat\Sigma\Sigma^{-1/2}\right)}\Norm{\frac{1}{n}\sum_{i=1}^n\Sigma^{-1/2}X_i(Y_i-X_i^\top\beta_0)}_2^2.
    \end{eqnarray*} Hence, \eqref{eq:prop.11.2} follows from Theorem~4.1 of \cite{oliveira2016lower} and Proposition~\ref{prop:C4}.
\end{proof}

\end{document}